\documentclass[a4paper, leqno, 12pt]{amsart}
 \usepackage[usenames]{color}
 \usepackage{amssymb, amsthm, amsmath, amscd}
\usepackage{ascmac}
 \usepackage{fancybox}
 \usepackage{ulem}
 \usepackage[all]{xy}%xy-matrix
 \usepackage{comment}%コメント
 \usepackage{here}%図に必要
 \usepackage{float}%図の位置を並べる
 \input amssym.def
 \input amssym.tex
 \newcommand{\Z}{\mathbb{Z}}
 
 \newcommand{\R}{\mathbb{R}}
 
  \newcommand{\Ov}[1]{\overline{#1}}
 \newcommand{\Ti}[1]{\widetilde{#1}}
  \newcommand{\Ch}[1]{\check{#1}}
  \newcommand{\Bb}[1]{\mathbb{#1}}
  
   \newcommand{\Rm}[1]{\mathrm{#1}}
    
\newtheorem{prop}{Proposition}[section]
 \newtheorem{thm}[prop]{Theorem}
 
 \newtheorem{cor}[prop]{Corollary}

\theoremstyle{remark}
\newtheorem{define}[prop]{Definition}

\numberwithin{equation}{section} % 式番号を「(3.5)」のように印刷
 
 \makeatletter
\@namedef{subjclassname@2020}{2020 Mathematical Subject Classification}
\makeatother

\begin{document}
\baselineskip=18pt
 \title{Local delta invariants of weak del Pezzo surfaces with the anti-canonical degree $\geq 5$}
\author{Hiroto Akaike}
%\date{}

 
 
\maketitle
%\tableofcontents

\begin{abstract}
The delta invariant interprets the criterion for the K-(poly)stability of log terminal Fano varieties.
In this paper, we determine the whole local delta invariant for all weak del Pezzo surfaces with the anti-canonical degree $\geq 5$.
\end{abstract}

\section*{Introduction}
Throughout the paper, we work out over the complex number field $\Bb{C}$. 
It is an important problem whether does a log terminal Fano variety $X$ admits a weak K\"ahler Einstein metric or not.
Recently, 
the equivalence for the K-polystability of $X$ and the existence of weak K\"ahler Einstein metric on $X$ was proved. 
%(\cite{Ber16},\cite{CDS},\cite{Ti15}).
In order to interpret the criterion for the K-(poly)stability of $X$,
the delta invariant $\delta(X)$ is introduced in \cite{BJ},\cite{FO}.
%It is also important to serve as a criterion for the K-polystability of $X$.
In fact, it is known by \cite{BJ},\cite{BX},\cite{FO},\cite{Fjt19},\cite{Li17},\cite{LXZ} that 
\[
\delta(X) > 1 \iff 
X \text{ is K-stable }%\sharp\Rm{Aut}(X)< \infty.
\iff
X \text{ is K-polystable and }\sharp\Rm{Aut}(X)< \infty.
\]

The delta invariants of smooth del Pezzo surfaces are known in \cite[\S 2]{FAND}.
On the other hand, it is hard to estimate the delta invariant for higher dimensional Fano varieties.
In order to estimate the delta invariant, 
Abban and Zhuang introduced the local delta invariant instead of the delta invariant in \cite{AZ}.
They gave an important idea for reducing the estimation of local delta invariants to that of lower dimensional cases.

We recall the definition of the local delta invariant.
Let $X$ be a $n$-dimensional weak log Fano variety with at worst log terminal singularities
and let $p \in X$ be a closed point.
Take a projective birational morphism $\sigma: Y \to X$ with smooth variety $Y$ and a prime divisor $E$ on $Y$.
We call $E$ a divisor over $X$.
Let 
\[
A_{X}(E):=1+ \Rm{ord}_{E}(K_{Y}-\sigma^{\ast}K_{X}),
\]
and we let
\[
S(E):=\frac{1}{(-K_{X})^{n}}\int_{0}^{\tau}\Rm{vol}(\sigma^{\ast}(-K_{X})-uE)du,
\]
where $\tau$ is the pseudo effective threshold of $E$ with respect to $-K_{X}$, that is, 
\[
\tau := \Rm{sup} \{ u \in \Bb{Q}_{\geq 0}\mid \sigma^{\ast}(-K_{X})-uE \text{ is big } \}.
\]
The local delta invariant $\delta_{p}(X)$ of $X$ at $p \in X$ is defined as follows:
\[
\delta_{p}(X):=\Rm{inf}\left\{\frac{A_{X}(E)}{S(E)}\mid \text{$E$ is a prime divisor over $X$ with $p \in C_{X}(E)$} \right\}.
\]
Moreover, the delta invariant $\delta(X)$ of $X$ is given by 
\[
\delta(X):=\Rm{inf}_{p \in X}\delta_{p}(X).
\]

In this paper, for every weak del Pezzo surface $S$ with the anti-canonical degree $\geq 5$,
we calculate the local delta invariant at each $p \in S$.
These results are important in the following aspects:
\begin{itemize}
\item[(1)] As a directed corollary, we give the delta invariant for a smooth weak del Pezzo surface with the anti-canonical degree $\geq 5$.
Moreover, if $S \to \Ov{S}$ is the anti-canonical model of $S$,
then we can immediately give $\delta_{\Ov{p}}(\Ov{S})$ for each $\Ov{p} \in \Ov{S}$.
In particular, we get the value $\delta(\Ov{S})$ for any du Val del Pezzo surface with the anti-canonical degree $\geq 5$.
 
\item[(2)] The estimation of the local delta invariant of weak del Pezzo surfaces is useful for the K-stability of higher dimensional Fano varieties.
In fact, the estimation of the local delta invariant of the quintic del Pezzo surfaces plays a crucial role in determining the K-stability of certain  Fano $3$-folds in \cite[Lemma~24, 25]{CFKO}.
Our results seem to be useful for determining the K-stability of other higher dimensional Fano varieties.
\end{itemize}
Now, the main results can be stated as follows.

We state the local delta invariants of weak del Pezzo surfaces with the anti-canonical degree $5$.
It is known that there exist $7$ types of weak del Pezzo surfaces of the anti-canonical degree $5$
in terms of the configuration of negative curves (\cite{CT}).
\begin{thm}
Let $S$ be a weak del Pezzo surface with the anti-canonical degree $5$.
The symbols $(E_i, \bullet)$ and $(F_j, \circ)$ denote $(-1)$-curve and $(-2)$-curve, respectively.
The local delta invariants $\delta_{p}(S)$ of $S$ at $p \in S$ are as follows.

\noindent
$(1)$ If the configuration of negative curves of $S$ is
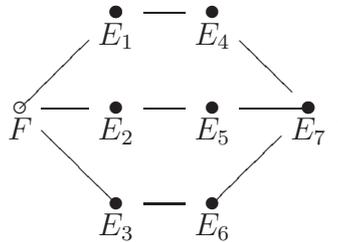
\[
\xygraph{
   {\circ}* +!D{ }*+!U{F}
      (- [ru] \bullet 
    {\bullet}* +!U{E_{1}} - [r]
    {\bullet}*+!U{E_{4}}  
   %     - [r] \cdots
   %     - [r] \bullet 
   %     - [r] \bullet 
       - [rd] \bullet
       ,
                -[r]
    {\bullet}*+!U{E_{2}}  - [r]
    {\bullet}*+!U{E_{5}}
   %     - [r] \cdots
   %     - [r] \bullet 
    %    - [r] \bullet 
         - [r] 
         ,
                -[rd]
    {\bullet}*+!U{E_{3}} - [r]
    {\bullet}*+!U{E_{6}}
   %     - [r] \cdots
  %      - [r] \bullet 
  %      - [r] \bullet 
         - [ru] {\bullet}*+!U{E_{7}}
         },
        \]
then the local delta invariants $\delta_{p}(S)$ of $S$ at $p \in S$ are as follows.
\begin{table}[H]
{\renewcommand\arraystretch{1.5}
 \begin{tabular}{|c||c|c|c|c|c|} \hline
     $p \in S$                     &  $F$   & $E_{i}\setminus F$ $(i=1,2,3)$ & $E_{i+3}\setminus E_{i}$ $(i=1,2,3)$ & $E_7$ &$S \setminus (F \cup \bigcup_{i=1}^{7} E_{i})$ \\ \hline 
  $\delta_{p}(S)$ & $\frac{15}{17}$ & $ 1 $  &  $\frac{15}{13}$                                                        & $\frac{15}{13}$  & $\frac{4}{3}$ \\ \hline
  \end{tabular}
  }
\end{table}
%%%%%%%%%%%%%

\noindent
$(2)$ If the configuration of negative curves of $S$ is
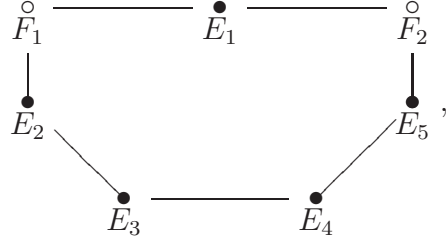
\[
\xygraph{
   {\bullet}* +!D{ }*+!U{E_{2}}
      (- [u] 
    {\circ}* +!U{F_{1}} - [rr]
    {\bullet}*+!U{E_{1}}  -[rr]
       {\circ}*+!U{F_{2}}  
   %     - [r] \cdots
   %     - [r] \bullet 
   %     - [r] \bullet 
       - [d] {\bullet}*+!U{E_{5}}
    ,
                -[rd]
    {\bullet}*+!U{E_{3}} - [rr]
    {\bullet}*+!U{E_{4}}
   %     - [r] \cdots
  %      - [r] \bullet 
  %      - [r] \bullet 
         - [ru] \bullet
         },
        \]
then the local delta invariants $\delta_{p}(S)$ of $S$ at $p \in S$ are as follows. 
\begin{table}[H]
{\renewcommand\arraystretch{1.5}
 \begin{tabular}{|c||c|c|c|c|c|} \hline
      $p \in S$                    &  $E_1$   & $F_{1}\setminus E_1$,  $F_{2}\setminus E_1$ & $E_{2}\setminus F_{1}$, $E_{5}\setminus F_{2}$ & 
  $E_{3}\setminus E_{2}$, $E_{4}\setminus E_{5}$ &$S \setminus  \bigcup_{i,j} (E_{i}\cup F_j)$ \\ \hline 
  $\delta_{p}(S)$ & $\frac{15}{19}$ &  $\frac{15}{17}$   &  $1$    & $\frac{15}{13}$  & $\frac{4}{3}$ \\ \hline
  \end{tabular}
  }
\end{table}
%%%%%%%%%%%%%%

\noindent
$(3)$ If the configuration of negative curves of $S$ is
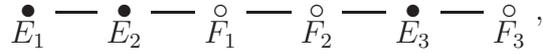
\[
\xygraph{
   {\bullet}* +!D{ }*+!U{E_1}
                    -[r]
    {\bullet}*+!U{E_{2}}  - [r]
    {\circ}*+!U{F_{1}}
   %     - [r] \cdots
       - [r] {\circ}* +!U{F_{2}}
       - [r] {\bullet}* +!U{E_{3}}
         - [r] {\circ}*+!U{F_{3}}
              },
        \]
then the local delta invariants $\delta_{p}(S)$ of $S$ at $p \in S$ are as follows.
\begin{table}[H]
{\renewcommand\arraystretch{1.5}
 \begin{tabular}{|c||c|c|c|c|c|c|} \hline
        $p \in S$                  &  $E_1 \setminus E_2$   & $E_{2}\setminus F_1$,  $F_{3}\setminus E_3$ & $F_{1}\setminus F_{2}$ & 
  $F_{2}\setminus E_{3}$ & $E_3$  & $S \setminus  \bigcup_{i,j} (E_{i}\cup F_j)$ \\ \hline 
  $\delta_{p}(S)$ & $\frac{15}{13}$ &  $\frac{15}{17}$   &  $\frac{15}{19}$    & $\frac{5}{7}$  & $\frac{15}{23}$ & $\frac{30}{23}$ \\ \hline
  \end{tabular}
  }
\end{table}
%%%%%%%%%%%%%%

\noindent
$(4)$ If the configuration of negative curves of $S$ is
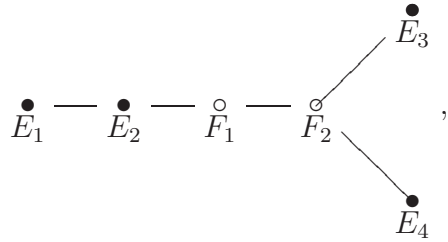
\[
\xygraph{
    {\bullet}* +!U{E_{1}} - [r]
    {\bullet}*+!U{E_{2}}  -[r]
    {\circ}* +!U{F_{1}} - [r]
   {\circ}* +!D{ }*+!U{F_{2}}
      (- [ru] \bullet 
    {\bullet}* +!U{E_{3}} 
   % {\bullet}*+!U{E_{4}}  
   %     - [r] \cdots
   %     - [r] \bullet 
   %     - [r] \bullet 
     %  - [rd] \bullet
              ,
                -[rd]
    {\bullet}*+!U{E_{4}} 
   % {\bullet}*+!U{E_{6}}
   %     - [r] \cdots
  %      - [r] \bullet 
  %      - [r] \bullet 
    %     - [ru] {\bullet}*+!U{E_{7}}
         },
        \]
then the local delta invariants $\delta_{p}(S)$ of $S$ at $p \in S$ are as follows. 
\begin{table}[H]
{\renewcommand\arraystretch{1.5}
 \begin{tabular}{|c||c|c|c|c|c|c|} \hline
      $p \in S$                    &  $E_1 \setminus E_2$   & $E_{2}\setminus F_1$ &  $F_{1}\setminus F_2$ & $ F_{2}$ & 
  $E_{i}\setminus F_{2}$ $(i=3,4)$  & $S \setminus  \bigcup_{i,j} (E_{i}\cup F_j)$ \\ \hline 
  $\delta_{p}(S)$ & $\frac{15}{13}$ & $\frac{15}{17}$ &   $\frac{15}{19}$ &  $\frac{5}{7}$    & $\frac{30}{31}$  &  $\frac{30}{23}$ \\ \hline
  \end{tabular}
  }
\end{table}
%%%%%%%%%%%%%%

\noindent
$(5)$ If the configuration of negative curves of $S$ is
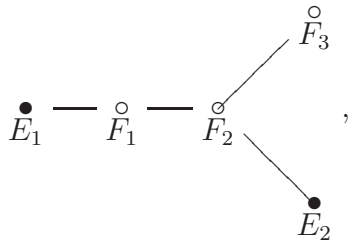
\[
\xygraph{
    {\bullet}* +!U{E_{1}} - [r]
    {\circ}*+!U{F_{1}}  -[r]
 %   {\circ}* +!U{F_{1}} - [r]
   {\circ}* +!D{ }*+!U{F_{2}}
      (- [ru] 
    {\circ}* +!U{F_{3}} 
   % {\bullet}*+!U{E_{4}}  
   %     - [r] \cdots
   %     - [r] \bullet 
   %     - [r] \bullet 
     %  - [rd] \bullet
              ,
                -[rd]
    {\bullet}*+!U{E_{2}} 
   % {\bullet}*+!U{E_{6}}
   %     - [r] \cdots
  %      - [r] \bullet 
  %      - [r] \bullet 
    %     - [ru] {\bullet}*+!U{E_{7}}
         },
        \]
then the local delta invariants $\delta_{p}(S)$ of $S$ at $p \in S$ are as follows.
\begin{table}[H]
{\renewcommand\arraystretch{1.5}
 \begin{tabular}{|c||c|c|c|c|c|c|} \hline
        $p \in S$                  &  $E_1 \setminus F_1$   & $F_{1}\setminus F_2$ &  $F_2$ & $ F_{3}\setminus F_2 $ & 
      $E_2 \setminus F_2 $ & $S \setminus  \bigcup_{i,j} (E_{i}\cup F_j)$ \\ \hline 
  $\delta_{p}(S)$ & $\frac{15}{16}$ & $\frac{30}{43}$ &   $\frac{5}{9}$ &  $\frac{15}{19}$    & $\frac{10}{13}$  & $\frac{5}{4}$  \\ \hline
  \end{tabular}
  }
\end{table}
%%%%%%%%%%%%%%

\noindent
$(6)$ If the configuration of negative curves of $S$ is
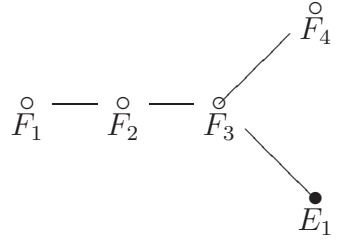
\[
\xygraph{
    {\circ}* +!U{F_{1}} - [r]
    {\circ}*+!U{F_{2}}  -[r]
 %   {\circ}* +!U{F_{1}} - [r]
   {\circ}* +!D{ }*+!U{F_{3}}
      (- [ru] 
    {\circ}* +!U{F_{4}} 
   % {\bullet}*+!U{E_{4}}  
   %     - [r] \cdots
   %     - [r] \bullet 
   %     - [r] \bullet 
     %  - [rd] \bullet
              ,
                -[rd]
    {\bullet}*+!U{E_{1}} 
   % {\bullet}*+!U{E_{6}}
   %     - [r] \cdots
  %      - [r] \bullet 
  %      - [r] \bullet 
    %     - [ru] {\bullet}*+!U{E_{7}}
         },
        \]
then the local delta invariants $\delta_{p}(S)$ of $S$ at $p \in S$ are as follows.
\begin{table}[H]
{\renewcommand\arraystretch{1.5}
 \begin{tabular}{|c||c|c|c|c|c|c|} \hline
      $p \in S$                    &  $F_1 \setminus F_2$   & $F_{2}\setminus F_3$ &  $F_3$ & $ F_{4}\setminus F_3 $ & 
      $E_1 \setminus F_3 $ & $S \setminus  (E_{1}\cup \bigcup_{j}  F_j)$ \\ \hline 
  $\delta_{p}(S)$ & $\frac{3}{4}$ & $\frac{6}{11}$ &   $\frac{3}{7}$ &  $\frac{9}{13}$    & $\frac{3}{5}$  & $\frac{6}{5}$  \\ \hline
  \end{tabular}
  }
\end{table}
%%%%%%%%%%%%%%%

\noindent
$(7)$ If $S$ is a del Pezzo surface with the anti-canonical degree $5$, then the local delta invariants $\delta_{p}(S)$ of $S$ at $p \in S$ are as follows. 
\begin{table}[H]
{\renewcommand\arraystretch{1.5}
 \begin{tabular}{|c||c|c|} \hline
        $p \in S$                  & $p$ lies on a $(-1)$ curve  & $p$ does NOT lies on any $(-1)$ curve  \\ \hline 
$\delta_{p}(S)$ &  $\frac{15}{13}$ & $\frac{40}{31}$  \\ \hline
  \end{tabular}
  }
\end{table}
\end{thm}

We state the local delta invariants of weak del Pezzo surfaces with the anti-canonical degree $6$.
It is known that there exist $6$ types of weak del Pezzo surfaces of the anti-canonical degree $6$
in terms of the configuration of negative curves (\cite{CT}).

%%%%%%%%%%%%%%%  degree 6
\begin{thm}
Let $S$ be a weak del Pezzo surface with the anti-canonical degree $6$.
The symbols $(E_i, \bullet)$ and $(F_j, \circ)$ denote $(-1)$-curve and $(-2)$-curve, respectively.
The local delta invariants at $p \in S$ are as follows.

\noindent
$(1)$ If the configuration of negative curves of $S$ is
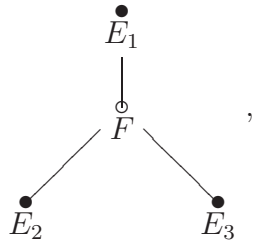
\[
\xygraph{
 %   {\circ}* +!U{F_{1}} - [r]
  %  {\circ}*+!U{F_{2}}  -[r]
 %   {\circ}* +!U{F_{1}} - [r]
   {\circ}* +!D{ }*+!U{F}
      (%- [ru] 
 %   {\circ}* +!U{F_{4}} 
   % {\bullet}*+!U{E_{4}}  
   %     - [r] \cdots
   %     - [r] \bullet 
   %     - [r] \bullet 
     %  - [rd] \bullet
              ,
                -[rd]
    {\bullet}*+!U{E_{3}} 
       ,
                -[ld]
    {\bullet}*+!U{E_{2}} 
     ,
                -[u]
    {\bullet}*+!U{E_{1}}
   % {\bullet}*+!U{E_{6}}
   %     - [r] \cdots
  %      - [r] \bullet 
  %      - [r] \bullet 
    %     - [ru] {\bullet}*+!U{E_{7}}
         },
        \]
then the local delta invariants $\delta_{p}(S)$ of $S$ at $p \in S$ are as follows.
\begin{table}[H]
{\renewcommand\arraystretch{1.5}
 \begin{tabular}{|c||c|c|c|} \hline
              $p \in S$           & $E_i \setminus F$ $(i=1,2,3)$ &  $F$  & $S \setminus  ( \bigcup_{i} E_{i}\cup  F)$  \\ \hline 
$\delta_{p}(S)$ &  $\frac{9}{10}$ & $\frac{3}{4}$  & $\frac{6}{5}$ \\ \hline
  \end{tabular}
  }
\end{table}
%%%%%%%%%%%%%%

\noindent
$(2)$ If the configuration of negative curves of $S$ is
\[
\xygraph{
   {\bullet}* +!D{ }*+!U{E_1}
                    -[r]
    {\bullet}*+!U{E_{2}}  - [r]
    {\circ}*+!U{F}
   %     - [r] \cdots
%       - [r] {\circ}* +!U{F_{2}}
       - [r] {\bullet}* +!U{E_{3}}
         - [r] {\bullet}*+!U{E_{4}}
              },
        \]
then the local delta invariants $\delta_{p}(S)$ of $S$ at $p \in S$ are as follows.
\begin{table}[H]
{\renewcommand\arraystretch{1.5}
 \begin{tabular}{|c||c|c|c|c|} \hline
            $p \in S$              &  $E_1 \setminus E_2$,  $E_4 \setminus E_3$  & $E_2$, $E_3$ &  $F \setminus (E_2 \cup E_3)$  
       & $S \setminus ( \bigcup_{i}E_{i}\cup F)$ \\ \hline 
  $\delta_{p}(S)$ & $\frac{9}{10}$ & $\frac{9}{11}$ &   $\frac{9}{11}$   & $\frac{9}{8}$  \\ \hline
  \end{tabular}
  }
\end{table}
%%%%%%%%%%%%%%%

\noindent
$(3)$ If the configuration of negative curves of $S$ is
\[
\xygraph{
   {\circ}* +!D{ }*+!U{F_1}
                    -[r]
    {\bullet}*+!U{E_{1}}  - [r]
    {\circ}*+!U{F_{2}}
   %     - [r] \cdots
%       - [r] {\circ}* +!U{F_{2}}
       - [r] {\bullet}* +!U{E_{2}}
       %  - [r] {\bullet}*+!U{E_{4}}
              },
        \]
 then the local delta invariants $\delta_{p}(S)$ of $S$ at $p \in S$ are as follows.
\begin{table}[H]
{\renewcommand\arraystretch{1.5}
 \begin{tabular}{|c||c|c|c|c|c|} \hline
         $p \in S$                 &  $F_{1}\setminus E_1$   & $E_1$ & $F_{2}\setminus E_{1}$ & 
  $E_{2}\setminus F_{2}$ &  $S \setminus  \bigcup_{i,j} (E_{i}\cup F_j)$ \\ \hline 
  $\delta_{p}(S)$ & $\frac{9}{11}$ &  $\frac{9}{14}$   &  $\frac{3}{4}$   & $\frac{9}{10}$  & $\frac{9}{8}$ \\ \hline
  \end{tabular}
  }
\end{table}
%%%%%%%%%%%%%%%

\noindent
$(4)$ If the configuration of negative curves of $S$ is
\[
\xygraph{
 %   {\circ}* +!U{F_{1}} - [r]
  %  {\circ}*+!U{F_{2}}  -[r]
 %   {\circ}* +!U{F_{1}} - [r]
   {\circ}* +!D{ }*+!U{F_2}
      (%- [ru] 
 %   {\circ}* +!U{F_{4}} 
   % {\bullet}*+!U{E_{4}}  
   %     - [r] \cdots
   %     - [r] \bullet 
   %     - [r] \bullet 
     %  - [rd] \bullet
              ,
                -[rd]
    {\bullet}*+!U{E_{2}} 
       ,
                -[l]
    {\circ}*+!U{F_{1}} 
     ,
                -[ru]
    {\bullet}*+!U{E_{1}}
   % {\bullet}*+!U{E_{6}}
   %     - [r] \cdots
  %      - [r] \bullet 
  %      - [r] \bullet 
    %     - [ru] {\bullet}*+!U{E_{7}}
         },
        \]
then the local delta invariants $\delta_{p}(S)$ of $S$ at $p \in S$ are as follows.
\begin{table}[H]
{\renewcommand\arraystretch{1.5}
 \begin{tabular}{|c||c|c|c|c|} \hline
        $p \in S$                  &  $F_{1}\setminus F_2$   & $F_2$ & $E_{1}\setminus F_{2}$,  $E_{2}\setminus F_{2}$& 
    $S \setminus  \bigcup_{i,j} (E_{i}\cup F_j)$ \\ \hline 
  $\delta_{p}(S)$ & $\frac{3}{4}$ &  $\frac{3}{5}$   &  $\frac{4}{5}$   &  $1$ \\ \hline
  \end{tabular}
  }
\end{table}
%%%%%%%%%%%%%%%

\noindent
$(5)$ If the configuration of negative curves of $S$ is
\[
\xygraph{
   {\circ}* +!D{ }*+!U{F_1}
                    -[r]
    {\circ}*+!U{F_{2}}  - [r]
    {\bullet}*+!U{E}
   %     - [r] \cdots
%       - [r] {\circ}* +!U{F_{2}}
       - [r] {\circ}* +!U{F_{3}}
       %  - [r] {\bullet}*+!U{E_{4}}
              },
        \]
 then the local delta invariants $\delta_{p}(S)$ of $S$ at $p \in S$ are as follows.
\begin{table}[H]
{\renewcommand\arraystretch{1.5}
 \begin{tabular}{|c||c|c|c|c|c|} \hline
             $p \in S$             &  $F_{1}\setminus F_2$   & $F_2 \setminus E$ & $E$ & 
  $F_{3}\setminus E $ &  $S \setminus  (E \cup \bigcup_j F_j)$ \\ \hline 
  $\delta_{p}(S)$ & $\frac{3}{4}$ &  $\frac{3}{5}$   &  $\frac{1}{2}$   & $\frac{3}{4}$  & $1$ \\ \hline
  \end{tabular}
  }
\end{table}
%%%%%%%%%%%%%%%%

\noindent
$(6)$ If the configuration of negative curves of $S$ is
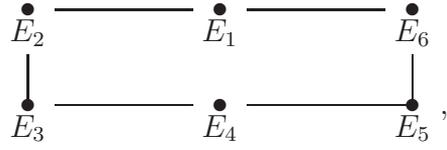
\[
\xygraph{
   {\bullet}* +!D{ }*+!U{E_{3}}
      (- [u] 
    {\bullet}* +!U{E_{2}} - [rr]
    {\bullet}*+!U{E_{1}}  -[rr]
       {\bullet}*+!U{E_6}  
   %     - [r] \cdots
   %     - [r] \bullet 
   %     - [r] \bullet 
       - [d] {\bullet}*+!U{E_{5}}
    ,
                -[rr]
    {\bullet}*+!U{E_{4}} - [rr]
%    {\bullet}*+!U{E_{4}}
   %     - [r] \cdots
  %      - [r] \bullet 
  %      - [r] \bullet 
     %    - [ru] \bullet
         },
        \]
then the local delta invariants $\delta_{p}(S)$ of $S$ at $p \in S$ are as follows.
\begin{table}[H]
{\renewcommand\arraystretch{1.5}
 \begin{tabular}{|c||c|c|} \hline
       $p \in S$                   & $E_{i}$ $(i=1,\cdots,6)$  &  $S \setminus \bigcup_{i}E_i$  \\ \hline 
$\delta_{p}(S)$ &  $1$ & $\frac{6}{5}$  \\ \hline
  \end{tabular}
  }
\end{table}
\end{thm}
%%%%%%%%%%%%%%% degree 7
We state the local delta invariants of weak del Pezzo surfaces with the anti-canonical degree $7$.
It is known that there exist $2$ types of weak del Pezzo surfaces of the anti-canonical degree $7$
in terms of the configuration of negative curves (\cite{CT}, \cite[\S 8.4]{Dol}).
\begin{thm}
Let $S$ be a weak del Pezzo surface with the anti-canonical degree $7$.
The symbols $(E_i, \bullet)$ and $(F , \circ)$ denote $(-1)$-curve and $(-2)$-curve, respectively.
The local delta invariants at $p \in S$ are as follows.

\noindent
$(1)$ If the configuration of negative curves of $S$ is
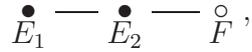
\[
\xygraph{
   {\bullet}* +!D{ }*+!U{E_1}
                    -[r]
    {\bullet}*+!U{E_{2}}  - [r]
    {\circ}*+!U{F}
   %     - [r] \cdots
%       - [r] {\circ}* +!U{F_{2}}
     %  - [r] {\circ}* +!U{F_{3}}
       %  - [r] {\bullet}*+!U{E_{4}}
              },
        \]
        then the local delta invariants $\delta_{p}(S)$ of $S$ at $p \in S$ are as follows.
\begin{table}[H]
{\renewcommand\arraystretch{1.5}
 \begin{tabular}{|c||c|c|c|c|} \hline
           $p \in S$               &  $E_{1}\setminus E_2$   & $E_2$ & $F \setminus E_{2}$ & 
    $S \setminus  (E_1 \cup E_2 \cup F )$ \\ \hline 
  $\delta_{p}(S)$ & $\frac{21}{25}$ &  $\frac{21}{31}$   &  $\frac{7}{9}$   &  $\frac{21}{23}$ \\ \hline
  \end{tabular}
  }
\end{table}
%%%%%%%%%%%%%%

\noindent
$(2)$ If the configuration of negative curves of $S$ is
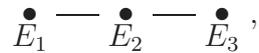
\[
\xygraph{
   {\bullet}* +!D{ }*+!U{E_1}
                    -[r]
    {\bullet}*+!U{E_{2}}  - [r]
    {\bullet}*+!U{E_3}
   %     - [r] \cdots
%       - [r] {\circ}* +!U{F_{2}}
     %  - [r] {\circ}* +!U{F_{3}}
       %  - [r] {\bullet}*+!U{E_{4}}
              },
        \]
then the local delta invariants $\delta_{p}(S)$ of $S$ at $p \in S$ are as follows.
\begin{table}[H]
{\renewcommand\arraystretch{1.5}
 \begin{tabular}{|c||c|c|c|} \hline
          $p \in S$                &  $E_{1}\setminus E_2$   & $E_2$ &  $S \setminus  (E_1 \cup E_2 \cup F )$ \\ \hline 
  $\delta_{p}(S)$ & $\frac{21}{23}$ &  $\frac{21}{25}$   &  $\frac{21}{22}$    \\ \hline
  \end{tabular}
  }
\end{table}
\end{thm}
%%%%%%%%%%%%%%%%degree 8
We state the local delta invariants for weak del Pezzo surfaces of the anti-canonical degree $8$.
Denote by $\pi:\Sigma_n \to \Bb{P}^1$  the $n$-th Hirzebruch surface.
Let $C_0$ be the section of $\pi$ with $C_0^2=-n$ and $\Gamma$ the fiber of $\pi$.
It is known that a weak del Pezzo surface of the anti-canonical degree $8$
is either $\Sigma_0$, $\Sigma_1$ or $\Sigma_2$ (\cite{CT},\cite[\S 8.4]{Dol}).
\begin{thm}
Let $S$ be a weak del Pezzo surface of the anti-canonical degree $8$.

\noindent
$(1)$ If $S = \Sigma_2$, then for any point $p \in S$, it holds that
\begin{equation}
\nonumber
  \delta_{p}(S)=\frac{3}{4}.
  \end{equation}
  
  \noindent
 $(2)$ If $S = \Sigma_1$, then for any point $p \in S$, then it holds that
\begin{equation}
\nonumber
  \delta_{p}(S)=
  \begin{cases}
    \frac{6}{7} & \text{if $p \in C_0$,} \\
   % \frac{21}{31}                 & \text{if $ p \in F_2 \setminus E $,} \\
     %\frac{21}{31}                 & \text{if $p \in E_2 $,} \\
    % \frac{21}{25} & \text{if $p \in  E_2 $} \\
     \frac{12}{13} & \text{if $p \in S \setminus C_0 . $}
  \end{cases}
\end{equation}

  \noindent
$(3)$ If $S = \Sigma_0=\Bb{P}^1 \times \Bb{P}^1$, then for any point $p \in S$, it holds that
\begin{equation}
\nonumber
  \delta_{p}(S)=1.
\end{equation}

%\begin{table}[H]
%{\renewcommand\arraystretch{1.5}
% \begin{tabular}{|c|c|} \hline
  %                       & $S$  \\ \hline 
%$\delta_{p}(S)$ &  $\frac{3}{4}$   \\ \hline
 % \end{tabular}
  %}
%\end{table}
%%%%%%%%%%%%%%%%
%\begin{table}[H]
%{\renewcommand\arraystretch{1.5}
 %\begin{tabular}{|c|c|c|} \hline
   %                      & $C_0$   &  $S \setminus C_0$  \\ \hline 
%$\delta_{p}(S)$ &  $\frac{6}{7}$ & $\frac{12}{13}$  \\ \hline
 % \end{tabular}
 % }
%\end{table}
%\begin{table}[H]
%{\renewcommand\arraystretch{1.5}
% \begin{tabular}{|c|c|} \hline
  %                       & $S$  \\ \hline 
%$\delta_{p}(S)$ &  $1$   \\ \hline
 % \end{tabular}
  %}
%\end{table}
\end{thm}
In the proof, Abban-Zhuang's method \cite[Theorem~3.2]{AZ} and its formula by intersection numbers \cite[Theorem1.106]{FAND} are crucial.
The key point of the proof is the following process:
For each point $p \in S$,
we choose a suitable plt blowing up $\Ti{S} \to S$ that extracts the divisor attaining $\delta_{p}(S)$.
How to choose such a plt blowing up is crucial.
All that remains to give the Zariski decomposition of divisors on $\Ti{S}$,
$\delta_{p}(S)$ is determined.

\textbf{Acknowledgememt}

The author is deeply grateful to Professor Kento Fujita for his valuable advice and support.
The research is supported by JSPS KAKENHI No. 20J20055.

\tableofcontents

\section*{Notation}
In this paper, we tacitly use the following notations.
\begin{itemize}
\item For a real vector $\mathbf{a} \in \Bb{R}^k$, we denote by $\mathbf{a}^{T}$ its transpose.
\item The symbol $\sim$ means the linearly equivalence between Cartier divisors.
\item We denote by $H$ a general hyperplane of $\Bb{P}^2$.
\item We denote by $\Ov{pq}$ the line on $\Bb{P}^2$ passing through two distinct points $p,q \in \Bb{P}^2$.
\item We denote by $\Rm{Bl}_{\{q_1, \cdots , q_k\}}\Bb{P}^2$ the surface obtained by the composition of the blowing-ups at $k$ distinct points $ q_1, \cdots , q_k \in \Bb{P}^2$.
\item Let $ \sigma: Y \to X$ be a birational morphism between projective varieties. For a Cartier divisor $D$ on $X$,
we denote by $\sigma_{\ast}^{-1}D$ the proper transform of it.
\end{itemize}

%\begin{align*}
 %H & : \text{a general hyperplane of }\Bb{P}^2 \\
 %\Ov{pq} & : \text{the line on }\Bb{P}^2 \text{ passing through }p,q \in \Bb{P}^2\\
 %\sigma_{\ast}^{-1}D & : \text{the proper transform of a Cartier divisor on }X \text{ by a birational morphism } \sigma: Y \to X
 %\end{align*}

\section{Preliminaries}
%%%%%%%%%%%%%%%%%%%%%%%%%%%%%%%%%%%%%%%%%%%%%%%%%%%%%
\begin{comment}
The delta invariant is introduced by \cite{FO} to serve as a criterion for the K-stability of Fano varieties.
We recall the important properties of delta invariant.

Let $X$ be a $n$-dimensional weak log Fano variety with at most log terminal singularities and $p \in X$ a closed point. 
Take a projective birational morphism $\sigma: Y \to X$ with smooth variety Y and a prime divisor $E$ on $Y$.
Let
\[
A_{X}(E):=1+ \Rm{ord}_{E}(K_{Y}-\sigma^{\ast}K_{X}),
\]
and we let
\[
S(E):=\frac{1}{(-K_{X})^n}\int_{0}^{\tau}\Rm{vol}(\sigma^{\ast}(-K_{X})-uE)du,
\]
where $\tau$ is the pseudo effective threshold of $E$ with respect to $-K_{X}$.
%If $S$ is a weak del Pezzo surface, then we write $S(E)=S_{-K_{S}}(E)$.
By \cite{BJ}, we may define the local delta invariant of $(S,L)$ at $p \in S$ as follows:
\[
\delta_{p}(X):=\Rm{inf}\left\{\frac{A_{X}(E)}{S(E)}\mid \text{$E$ is a prime divisor over $X$ with $p \in C_{X}(E)$} \right\}
\]
\end{comment}
%%%%%%%%%%%%%%%%%%%%%%%%%%%%%%%%%%%%%%%%%%%%%%%%%%%%%%%
In what follows, we state only for the case of weak del Pezzo surfaces.
By \cite[Theorem~1.106]{FAND}, we may define $S(W_{\bullet,\bullet}^{E},q)$ as follows:
%We need the following theorem to calculate the local delta invariant of weak del Pezzo surfaces.
\begin{define}[\cite{FAND},Theorem~1.106]
\label{S-inv}
Let $S$ be a weak del Pezzo surface,
$p \in S$ a closed point, 
$E$ a prime divisor over $S$ with $p \in C_{S}(E)$ and 
$\sigma:\Ti{S} \to S$ the plt blowing-up of $E$.
Assume $\Ti{S}$ is a Mori dream surface.
If $\Ti{P}(u)+\Ti{N}(u)$ is the Zariski Decomposition of $\sigma^{\ast}(-K_{S})-uE$,
then, for $q \in E$
\[
S(W_{\bullet,\bullet}^{E},q):= \frac{2}{(-K_{S})^2} \int_{0}^{\tau}(\Ti{P}(u)\cdot E)\cdot \Rm{ord}_{q}\Ti{N}(u)|_{E}\; du
+ \frac{1}{(-K_{S})^2} \int_{0}^{\tau} (\Ti{P}(u)\cdot E)^2 du,
\]
where $\tau = \tau(-K_{S})$ is the pseudo effective threshold of $E$ with respect to $-K_{S}$.
\end{define}

\begin{thm}[\cite{AZ},Theorem~3.2]
Let $S$ be a weak del Pezzo surface,
%$L$ a big $\Bb{Q}$ divisor,
$p \in S$ a closed point and $E$ a prime divisor over $S$ with $p \in C_{S}(E)$.
If $\sigma:\Ti{S} \to S$ is the plt blowing-up of $E$, then it holds that
\[
\delta_{p}(S) \geq \Rm{min}\left\{\frac{A_{S}(E)}{S(E)},\Rm{inf}\left\{\frac{A_{E,\Delta_{E}}(q)}{S(W_{\bullet,\bullet}^{E},q)}\mid q \in E,\sigma(q)=p\right\} \right\},
\]
where $\Delta_{E}$ is the $\Bb{Q}$-divisor such that $(K_{E}+\Delta_{E})=(K_{\Ti{S}}+E)|_{E}$ and 
$A_{E,\Delta_{E}}(q)$ is the log discrepancy of $q$.
\end{thm}

\begin{cor}
\label{AZ}
Let the notations as above.
\begin{itemize}
\item[(1)]
If $E$ is a smooth prime divisor passing through $p \in S$, then it holds 
\[
\delta_{p}(S) \geq \Rm{min}\left\{\frac{1}{S(E)},\frac{1}{S(W_{\bullet,\bullet}^{E},p)} \right\}.
\]
\item[(2)]
%Assume $S$ is smooth.
If $E$ is the exceptional curve by the ordinary blowing up $\sigma:\Ti{S} \to S$ at a point $p$,
then it holds 
\[
\delta_{p}(S) \geq \Rm{min}\left\{\frac{2}{S(E)},\Rm{inf}\left\{\frac{1}{S(W_{\bullet,\bullet}^{E},q)}\mid q \in E\right\} \right\}.
\]
\end{itemize}
\end{cor}

%\begin{rmk}
%For weak del Pezzo surface $S$ and a divisor $E$ over $S$, we write 
%\[
%S(E)=S_{-K_{S}}(E), \quad 
%\delta_{p}(S) = \delta_{p}(S,-K_{S}).
%\]
%\end{rmk}

\section{The case of the anti-canonical degree $5$}

It is known that there exist $7$ types of weak del Pezzo surfaces of the anti-canonical degree $5$
in terms of the configuration of negative curves (\cite{CT}).

\begin{prop}
Let $S$ be the anti-canonical degree $5$ weak del Pezzo surface such
that the dual graph of negative curves is 
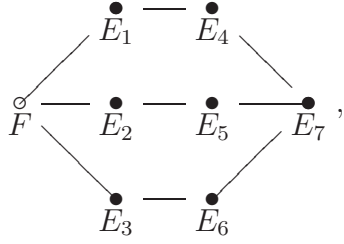
\[
\xygraph{
   {\circ}* +!D{ }*+!U{F}
      (- [ru] \bullet 
    {\bullet}* +!U{E_{1}} - [r]
    {\bullet}*+!U{E_{4}}  
   %     - [r] \cdots
   %     - [r] \bullet 
   %     - [r] \bullet 
       - [rd] \bullet
       ,
                -[r]
    {\bullet}*+!U{E_{2}}  - [r]
    {\bullet}*+!U{E_{5}}
   %     - [r] \cdots
   %     - [r] \bullet 
    %    - [r] \bullet 
         - [r] 
         ,
                -[rd]
    {\bullet}*+!U{E_{3}} - [r]
    {\bullet}*+!U{E_{6}}
   %     - [r] \cdots
  %      - [r] \bullet 
  %      - [r] \bullet 
         - [ru] {\bullet}*+!U{E_{7}}
         },
        \]

where $E_{i}$ $(i=1, \cdots, 7)$ is a $(-1)$-curve and $F$ is a $(-2)$-curve.
Then, for a point $p \in S$, it holds that
\begin{equation}
\nonumber
  \delta_{p}(S)=
  \begin{cases}
    \frac{15}{17} & \text{if $p \in F$,} \\
    1                 & \text{if $p \in E_{i}\setminus F$ for $i=1,2,3$,} \\
    \frac{15}{13}       & \text{if $p \in E_{i+3}\setminus E_{i}$ for $i=1,2,3$,}\\
   \frac{15}{13} & \text{if $p \in E_{7}$,}\\
   \frac{4}{3} & \text{if $p \in S \setminus (F \cup \bigcup_{i=1}^{7} E_{i})$.}
  \end{cases}
\end{equation}
\end{prop}

\begin{proof}
We recall the construction of $S$.
Take non-colinear three points $q_{0}, q_{1}, q_{3} \in \Bb{P}^{2}$ and $q_{2} \in \Ov{q_{1}q_{3}}\setminus\{q_1,q_3\}$.
Then  $S$ is obtained by $\rho: S = \Rm{Bl}_{\{q_{1}, q_{2}, q_{3}, q_{4}\}}\Bb{P}^{2} \to \Bb{P}^2$.
%We denote by $\rho: S \to \Bb{P}^2 $ this birational morphism.
Moreover, we have 
$F=\rho^{-1}_{\ast}\Ov{q_1 q_3}$, $E_1 = \rho^{-1}(q_1)$, $E_2 = \rho^{-1}(q_2)$,
$E_3 = \rho^{-1}(q_3)$, $E_4 = \rho^{-1}_{\ast}(\Ov{q_0 q_1})$,
$E_5 = \rho^{-1}_{\ast}(\Ov{q_0 q_2})$, $E_6 = \rho^{-1}_{\ast}(\Ov{q_0 q_3})$ and 
$E_7 = \rho^{-1}(q_0)$.
We denote a divisor $D=\sum_{i=1}^{7}a_{i} E_{i} + bF \in \Rm{Div}(S)$ ($a_i , b \in \Z$) by 
$D=(a_1,a_2,a_3,a_4,a_5,a_6,a_7,b)$.
The intersection matrix of $\{E_1, E_2, E_3, E_4, E_5, E_6, E_7, F \}$ is 
\[A:=
\left(
\begin{array}{ccccccc|c}
  -1 & 0 & 0  & 1  &  0  &  0  &  0  & 1 \\
  0 & -1 & 0  & 0  &  1  &  0  &  0  & 1 \\
  0 & 0 & -1  & 0  &  0 &  1  &  0  & 1 \\
  1 & 0 & 0  & -1  &  0  &  0  &  1  & 0 \\
  0 & 1 & 0  & 0  &  -1  &  0  &  1  & 0 \\
  0 & 0 & 1  & 0  &  0  &  -1  &  1  & 0 \\
  0 & 0 & 0  & 1 &  1  &  1  &  -1  & 0 \\ \hline
  1 & 1 & 1  & 0  &  0  &  0  &  0  & -2 
\end{array}
\right).
\]
We note that $-K_{S}\sim (0,0,0,1,1,1,2,0)$.

\noindent
(1)\;The case $p \in F$.

\noindent
We calculate $S(F)$ and $S(W_{\bullet,\bullet}^{F},p)$
in order to apply Corollary~\ref{AZ} for prime divisor $F$.
Take $u \in \R_{\geq 0}$. 
Let $P(u)+N(u)$ be the Zariski decomposition of $-K_{S}-uF$,
where $P(u)$ is the positive part and $N(u)$ is the negative part.
%%%%%%%%%%%%%%%%%%%%%%%%%%%%%%%%%%%%%%%%%%%%%%
\begin{comment}
Then we have
\begin{align*}
\nonumber
 &P(u)=
  \begin{cases}
   -K_{S}-uF & (u\in [0,1]),\\
                    & \\
    (3-u)F+(2-u)(E_{1}+E_{2}+E_{3})+(1-u)(E_{4}+E_{5}+E_{6})+2E_{7}      &( u\in [1,2]), \\
  \end{cases}\\
    &N(u)=
  \begin{cases}
    0 & (u\in [0,1]) \\
                     &  \\
    (u-1)E_{4}+E_{5}+E_{6})     & (u\in [1,2])\\
  \end{cases}
\end{align*}
\end{comment}
%
If $u\in [0,1]$, then we have
\begin{align*}
&P(u)=(0,0,0,1,1,1,2,-u),\\
& N(u)=0.
\end{align*}
If $u\in [1,2]$, then we have
\begin{align*}
&P(u)=(1-u,1-u,1-u,1,1,1,2,-u),\\
& N(u)=(u-1,u-1,u-1,0,0,0,0,0).
\end{align*}
We note that $-K_{S}-uF$ is not pseudo effective for $u>2$.
Therefore, 
if $u\in [0,1]$, then we have
\begin{align*}
P(u)F
=&(0,0,0,1,1,1,2,-u)A(0,0,0,0,0,0,0,1)^{T}  \\
=&2u, \\
P(u)^{2}
=&(0,0,0,1,1,1,2,-u)A(0,0,0,1,1,1,2,-u)^{T} \\
=&5-2u^2,
\end{align*}
and
if $u\in [1,2]$, then we have
\begin{align*}
P(u)F
=&(1-u,1-u,1-u,1,1,1,2,-u)A(0,0,0,0,0,0,0,1)^{T}\\
=&3-u, \\
P(u)^{2}
=&(1-u,1-u,1-u,1,1,1,2,-u)A(1-u,1-u,1-u,1,1,1,2,-u)^{T}\\
=&(4-u)(2-u).
\end{align*}
Hence we get 
\begin{align*}
S(F)=\frac{1}{5}\int_{0}^{1}5-2u^2 du
+ \frac{1}{5}\int_{1}^{2} (4-u)(2-u) du
=\frac{17}{15}
\end{align*}
by the definition of $S(F)$ and 
\begin{align*}
S(W_{\bullet,\bullet}^{F},p)
&=%\frac{2}{5}\int_{0}^{1}2u \Rm{ord}_{p}(N(u)|_{F}) du
\frac{1}{5}\int_{0}^{1}4u^2 du
+ \frac{2}{5}\int_{1}^{2}(3-u) \Rm{ord}_{p}(N(u)|_{F})  du
+\frac{1}{5}\int_{1}^{2}(3-u)^2 du\\
 &= 
 \begin{cases}
  \frac{11}{15} & \text{if $p \in E_{i}\cap F $ for $i=1,2,3$},\\
                    & \\
   \frac{7}{15}   & \text{if $p \in F \setminus \bigcup_{i=1}^{3} E_{i}$}, \\
  \end{cases}
\end{align*}
by Definition~\ref{S-inv}.
Hence we have 
\begin{align*}
\delta_{p}(S) \geq  \Rm{min}\left\{ \frac{1}{S(F)}, \frac{1}{S(W_{\bullet,\bullet}^{F},p)} \right\}
= \frac{15}{17} 
\end{align*}
from Corollary~\ref{AZ}.
On the other hand, we have 
\begin{align*}
\frac{A_{S}(F)}{S(F)}=\frac{15}{17} \geq \delta_{p}(S)   
\end{align*}
by the definition of the local delta invariant. Thus, we have $\delta_{p}(S)= 15/17$ in this case.

%%%%%%%%%%%%%%%%%%%%% (2)
\noindent
(2)\;The case $p \in E_{i} \setminus (F \cup E_{i+3}) $ for $i=1,2,3$.

\noindent
%To apply \Red{Theorem~Abban-Zhang} for prime divisor $F$,
%We show only the case $E_1$.
%\Red{$E_{1}$についてのみ示す,$E_2,E_3$も同様}\\
We calculate $S(E_{1})$ and $S(W_{\bullet,\bullet}^{E_{1}},p)$.
Take $u \in \R_{\geq 0}$. 
Let $P(u)+N(u)=-K_{S}-uE_{1}$ be the Zariski decomposition,
where $P(u)$ is the positive part and $N(u)$ is the negative part.
If $u\in [0,1]$, then we have
\begin{align*}
&P(u)=\left(  -u  , 0   , 0    ,  1   , 1    ,1     ,2     , -\frac{u}{2}\right), \\
&N(u)=\left(   0 ,  0  , 0    ,  0   , 0    , 0    , 0    , \frac{u}{2}\right).
\end{align*}
If $u\in [1,2]$, then we have
\begin{align*}
&P(u)=\left( -u   , 0   , 0    , 2-u    ,1     ,1     ,2     ,-\frac{u}{2}\right),\\
& N(u)=\left(  0  ,  0  ,  0   , u-1    ,0     ,0     ,0     ,\frac{u}{2}\right).
\end{align*}
We note that $-K_{S}-uE_{1}$ is not pseudo effective for $u>2$.
Therefore, 
if $u\in [0,1]$, then we have
\begin{align*}
P(u)E_{1}=\frac{u+2}{2}, \quad
P(u)^{2}=5-2u-\frac{u^2}{2},
\end{align*}
and
if $u\in [1,2]$, then we have
\begin{align*}
P(u)E_{1}=\frac{4-u}{2}, \quad
P(u)^{2}=\frac{1}{2}(6-u)(2-u).
\end{align*}
Hence we get 
\begin{align*}
S(E_1)=\frac{1}{5}\int_{0}^{1}\left(5-2u-\frac{u^2}{2}\right) du
+ \frac{1}{5}\int_{1}^{2} \frac{1}{2}(6-u)(2-u)du
=1
\end{align*}
by the definition of $S(E_1)$ and 
\begin{align*}
S(W_{\bullet,\bullet}^{E_{1}},p)
%&=\frac{2}{5}\int_{0}^{1} \frac{u+2}{2}\Rm{ord}_{p}(N(u)|_{F}) du
&=\frac{1}{5}\int_{0}^{1}\left(\frac{u+2}{2}\right)^2 du
%
%+ \frac{2}{5}\int_{1}^{2}\frac{4-u}{2}\Rm{ord}_{p}(N(u)|_{E_{1}})  du
+\frac{1}{5}\int_{1}^{2}\left(\frac{4-u}{2}\right)^2 du = \frac{19}{30}
\end{align*}
by Definition~\ref{S-inv}.
Hence we have 
\begin{align*}
1 \geq \delta_{p}(S) \geq  \Rm{min}\left\{ \frac{1}{S(E_1)}, \frac{1}{S(W_{\bullet,\bullet}^{E_1},p)} \right\}
= 1
\end{align*}
from Corollary~\ref{AZ}.
Thus, we have $\delta_{p}(S)= 1$ in this case.
We can show $\delta_{p}(S)= 1$ for $p \in E_{i} \setminus (F \cup E_{i+3}) (i=2,3)$ by the same calculation. 
%%%%%%%%%%%%%%%%%%%%%%%%%%%%%%%%%%%%%%%%%%%%%%%%(3)

\noindent
(3)\;The case $p \in E_{i}\setminus E_7$ for $i=4,5,6$.

\noindent
%We show only for $E_{4}$.
%\Red{$E_{4}$についてのみ示す. $E_5,E_6$も同様}\\
%To apply \Red{Theorem~Abban-Zhang} for prime divisor $F$,
We calculate $S(E_{4})$ and $S(W_{\bullet,\bullet}^{E_{4}},p)$.
Take $u \in \R_{\geq 0}$. 
Let $P(u)+N(u)$ be the Zariski decomposition of $-K_{S}-uE_{4}$,
where $P(u)$ is the positive part and $N(u)$ is the negative part.
If $u\in [0,1]$, then we have
\begin{align*}
&P(u)=\left(0,0 ,0,1-u,1,1,2,0\right), \\
%(1-u)E_4 + E_5 + E_6 +2E_7\sim -K_{S}-uE_4, \\
&N(u)=0.
\end{align*}
If $u\in [1,2]$, then we have
\begin{align*}
&P(u)=\left(   2-2u , 0, 0, 1-u, 1    , 1    , 3-u ,  1-u\right),\\
%(1-u)F+(2-2u)E_1+ (1-u)E_4 + E_5 + E_6 + (3-u)E_7 ,\\
& N(u)=\left(  2(u-1),0,0,0, 0,  0   , u-1,   u-1\right).
%(u-1)(F+2E_1+E_7).
\end{align*}
We note that $-K_{S}-uE_{4}$ is not pseudo effective for $u>2$.
Therefore, 
if $u\in [0,1]$, then we have
\begin{align*}
P(u)E_{4}=1+u, \quad
P(u)^{2}=5-2u-u^2,
\end{align*}
and
if $u\in [1,2]$, then we have
\begin{align*}
P(u)E_{4}=4-2u, \quad
P(u)^{2}=2(2-u)^2.
\end{align*}
Hence we get 
\begin{align*}
S(E_4)=\frac{1}{5}\int_{0}^{1}(5-2u-u^2) du
+ \frac{1}{5}\int_{1}^{2} 2(2-u)^2du
=\frac{13}{15}
\end{align*}
by the definition of $S(E_4)$ and 
\begin{align*}
S(W_{\bullet,\bullet}^{E_{4}},p)
&=%\frac{2}{5}\int_{0}^{1}2u \Rm{ord}_{p}(N(u)|_{F}) du
\frac{1}{5}\int_{0}^{1}(1+u)^2 du
+ \frac{2}{5}\int_{1}^{2}(4-2u) \Rm{ord}_{p}(N(u)|_{E_{4}})  du
+\frac{1}{5}\int_{1}^{2}(4-2u)^2 du\\
 &= 
 \begin{cases}
  1 & \text{if $p \in (E_{1} \cap E_{4})$},\\
                    & \\
   \frac{11}{15}   &  \text{if $p \in E_{4}\setminus(E_1 \cup E_{7}) $}, \\
  \end{cases}
\end{align*}
by Definition~\ref{S-inv}.
Hence we have 
\begin{align*}
\delta_{p}(S) \geq  \Rm{min}\left\{ \frac{1}{S(E_4)}, \frac{1}{S(W_{\bullet,\bullet}^{E_4},p)} \right\}
= \begin{cases}
  1 & \text{if $p \in (E_{1} \cap E_{4})$},\\
                    & \\
   \frac{15}{13}   &  \text{if $p \in E_{4}\setminus(E_1 \cup E_{7}) $}, \\
  \end{cases}
\end{align*}
from Corollary~\ref{AZ}.
 Hence we have $\delta_{p}(S)= 15/13$ for $p \in E_{4}\setminus(E_1 \cup E_{7})$.
 If $\{p\}=E_{1} \cap E_{4}$, we have $1=S(E_1) \geq \delta_{p}(S)$ by the calculation in (2).
 Thus, we have 
 \begin{align*}
 \delta_{p}(S)= 
 \begin{cases}
  1 & \text{if $p \in (E_{1} \cap E_{4})$},\\
                    & \\
   \frac{15}{13}   &  \text{if $p \in E_{4}\setminus(E_1 \cup E_{7}) $}. \\
  \end{cases}
 \end{align*}
 We can show 
 \begin{align*}
 \delta_{p}(S)= 
 \begin{cases}
  1 & \text{if $p \in (E_{i-3} \cap E_{i})$},\\
                    & \\
   \frac{15}{13}   &  \text{if $p \in E_{i}\setminus(E_{i-3} \cup E_{7}) $}, \\
  \end{cases}
 \end{align*}
for $i=5,6$ by the same calculation.

%%%%%%%%%%%%%%%%%%%%%%%%%%%%%%%%%%%%%%%%%%%%(4)
\noindent
(4)\;The case $p \in E_{7}$.

\noindent
%To apply \Red{Theorem~Abban-Zhang} for prime divisor $F$,
We calculate $S(E_{7})$ and $S(W_{\bullet,\bullet}^{E_{7}},p)$.
Take $u \in \R_{\geq 0}$. 
Let $P(u)+N(u)$ be the Zariski decomposition $-K_{S}-uE_{7}$,
where $P(u)$ is the positive part and $N(u)$ is the negative part.
If $u\in [0,1]$, then we have
\begin{align*}
&P(u)=\left( 0,0,0,1,1,1,2-u,0    \right),\\
%E_4 + E_5 + E_6 +(2-u)E_7\sim -K_{S}-uE_7, \\
&N(u)=0.
\end{align*}
If $u\in [1,2]$, then we have
\begin{align*}
&P(u)=(2-u)\left( 0, 0,0,1,1,1 ,1 ,0 \right), \\
%(E_{4}+E_{5}+E_{6}+E_{7}),\\
& N(u)=(u-1)\left( 0,0 ,0 ,1 ,1 ,1 ,0 , 0\right).
%(E_{4}+E_{5}+E_{6}).
\end{align*}
We note that $-K_{S}-uE_{7}$ is not pseudo effective for $u>2$.
Therefore, 
if $u\in [0,1]$, then we have
\begin{align*}
P(u)E_{7}=1+u, \quad
P(u)^{2}=5-2u-u^2,
\end{align*}
and
if $u\in [1,2]$, then we have
\begin{align*}
P(u)E_{7}=4-2u, \quad
P(u)^{2}=2(2-u)^2.
\end{align*}
Hence we get 
\begin{align*}
S(E_7)=\frac{1}{5}\int_{0}^{1}(5-2u-u^2) du
+ \frac{1}{5}\int_{1}^{2} 2(2-u)^2du
=\frac{13}{15}
\end{align*}
by the definition of $S(E_7)$ and 
\begin{align*}
S(W_{\bullet,\bullet}^{E_{7}},p)
&=%\frac{2}{5}\int_{0}^{1}2u \Rm{ord}_{p}(N(u)|_{F}) du
\frac{1}{5}\int_{0}^{1}(1+u)^2 du
+ \frac{2}{5}\int_{1}^{2}(4-2u) \Rm{ord}_{p}(N(u)|_{E_{7}})  du
+\frac{1}{5}\int_{1}^{2}(4-2u)^2 du\\
 &= 
 \begin{cases}
  \frac{13}{15} & \text{if $p \in E_{i}\cap E_{7} $ for $i=4,5,6$},\\
                    & \\
   \frac{11}{15}   & \text{if $p \in E_{7} \setminus \bigcup_{i=4}^{6} E_{i}$}, \\
  \end{cases}
\end{align*}
by Definition~\ref{S-inv}.
Hence we have 
\begin{align*}
\frac{15}{13}  \geq \delta_{p}(S) \geq  \Rm{min}\left\{ \frac{1}{S(E_7)}, \frac{1}{S(W_{\bullet,\bullet}^{E_7},p)} \right\}
= \frac{15}{13} 
\end{align*}
from Corollary~\ref{AZ}.
 Thus, we have $\delta_{p}(S)= 15/13$ in this case.
%%%%%%%%%%%%%%%%%%%%% (5)

\noindent
(5)\;The case $p \in S \setminus \left(F \cup \bigcup_{i=1}^{7} E_{i} \right)$.

\noindent
%To apply \Red{Theorem~Abban-Zhang} for prime divisor $F$,
Consider a blowing up $\sigma:\Ti{S} \to S$ at $p$.
Let $E$ be the exceptional curve of $p$,
$\Ti{F}$ and $\Ti{E}_{i}$ be the proper transform of $F$ and $E_{i}$, respectively.
Put $G_{i}:=(\rho\sigma)_{\ast}^{-1}\Ov{\rho(p)q_i}$ for $i=0,1,2,3$.
Then we have $\sigma^{\ast}(-K_{S})-uE \sim G_0 + G_2 + \Ti{F} + \Ti{E}_2 + (2-u)E$.
We calculate $S(E)$ and $S(W_{\bullet,\bullet}^{E},p)$.
Take $u \in \R_{\geq 0}$. 
Let $\Ti{P}(u)+\Ti{N}(u)$ be the Zariski decomposition of $\sigma^{\ast}(-K_{S})-uE$,
where $\Ti{P}(u)$ is the positive part and $\Ti{N}(u)$ is the negative part.
If $u\in [0,2]$, then we have
\begin{align*}
&\Ti{P}(u)=\Ti{E}_2 + \Ti{F} +  G_0 + G_2 +   (2-u)E,\\
&\Ti{N}(u)=0.
\end{align*}
If $u\in [2,5/2]$, then we have
\begin{align*}
&\Ti{P}(u)=\Ti{E}_2 +  (3-u) \Ti{F} + (5-2u)G_0 + (2-u)G_1+ (3-u)G_2 + (2-u)G_3 + (2-u)E,\\
&\Ti{N}(u)=(u-2)F + (2u-4)G_0 + (u-2)(G_1 + G_2 + G_3).
\end{align*}

We note that $\sigma^{\ast}(-K_{S})-uE$ is not pseudo effective for $u>5/2$.
If $u\in [0,2]$, then we have
\begin{align*}
\Ti{P}(u)E=u, \quad
\Ti{P}(u)^{2}=5-u^2,
\end{align*}
If $u\in [2,5/2]$, then we have
\begin{align*}
\Ti{P}(u)E=2(5-2u), \quad
\Ti{P}(u)^{2}=(5-2u)^2.
\end{align*}

Hence we get 
\begin{align*}
S(E)=\frac{1}{5}\int_{0}^{2} 5 -u^2 du + \frac{1}{5}\int_{2}^{\frac{5}{2}} (5-2u)^2 du =\frac{3}{2}
\end{align*}
by the definition of $S(E)$ and 
\begin{align*}
S(W_{\bullet,\bullet}^{E},p)
&=%\frac{2}{5}\int_{0}^{1}2u \Rm{ord}_{p}(N(u)|_{F}) du
\frac{1}{5}\int_{0}^{2}u^2 du
+ \frac{2}{5}\int_{2}^{\frac{5}{2}}2(5-2u) \Rm{ord}_{p}(\Ti{N}(u)|_{E})  du
+\frac{1}{5}\int_{2}^{\frac{5}{2}}4(5-2u)^2 du\\
 &= 
 \begin{cases}
  \frac{11}{15} & \text{if $p \in E \cap G_0 $},\\
                    & \\
   \frac{7}{10}   & \text{if $p \in E\cap G_i$ for $i=1,2,3$}, \\
    & \\
   \frac{2}{3}   & \text{if $p \in E \setminus \bigcup_{i=0}^{3} G_{i}$}, \\
  \end{cases}
\end{align*}
by Definition~\ref{S-inv}.
Hence we have 
\begin{align*}
\frac{4}{3}  \geq \delta_{p}(S) \geq  \Rm{min}\left\{ \frac{2}{S(E)}, \frac{1}{S(W_{\bullet,\bullet}^{E},p)} \right\}
= \frac{4}{3} 
\end{align*}
from Corollary~\ref{AZ}.
 Thus, we have $\delta_{p}(S)= 4/3$ in this case.
\end{proof}

\begin{prop}
Let $S$ be the anti-canonical degree $5$ weak del Pezzo surface such
that the dual graph of negative curves is
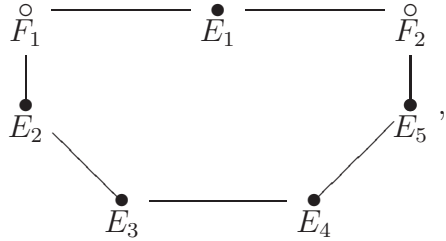
\[
\xygraph{
   {\bullet}* +!D{ }*+!U{E_{2}}
      (- [u] 
    {\circ}* +!U{F_{1}} - [rr]
    {\bullet}*+!U{E_{1}}  -[rr]
       {\circ}*+!U{F_{2}}  
   %     - [r] \cdots
   %     - [r] \bullet 
   %     - [r] \bullet 
       - [d] {\bullet}*+!U{E_{5}}
    ,
                -[rd]
    {\bullet}*+!U{E_{3}} - [rr]
    {\bullet}*+!U{E_{4}}
   %     - [r] \cdots
  %      - [r] \bullet 
  %      - [r] \bullet 
         - [ru] \bullet
         },
        \]
where $E_{i}$ $(i=1, \cdots, 5)$ is a $(-1)$-curve and $F_{j}$ $(j=1,2)$ is a $(-2)$-curve.
Then, for a point $p \in S$, it holds that
\begin{equation}
\nonumber
  \delta_{p}(S)=
  \begin{cases}
    \frac{15}{19} & \text{if $p \in E_1$,} \\
    \frac{15}{17}                & \text{if $p \in F_{1}\setminus E_{1}$ or $p \in F_{2}\setminus E_{1}$,} \\
    1       & \text{if $p \in E_{2}\setminus F_{1}$ or $p \in E_{5}\setminus F_{2}$,}\\
   \frac{15}{13} & \text{if $p \in E_{3} \setminus E_{2}$ or $p \in E_{4} \setminus E_{5}$,}\\
   \frac{4}{3} & \text{if $p \in S \setminus \left( \bigcup_{i,j}( E_{i} \cup F_{j})\right)$ .}
  \end{cases}
\end{equation}
\end{prop}

\begin{proof}
We can assume that we get $S$ from $\Bb{P}^2$ as follows.
\begin{itemize}
\item[(1)]Let $\rho_{1}:S_{1}=\Rm{Bl}_{\{q_{1}, q_{2}, q_{3}\}}\Bb{P}^{2} \to \Bb{P}^2$ be a 
blowing-up at non-colinear points $q_1$, $q_2$, $q_3$.
\item[(2)]Let $q_4$ be a point 
at which $\rho_{1}^{-1}(q_4)$ and $(\rho_{1})_{\ast}^{-1}\Ov{q_1 q_2}$ meet.
Take a blowing-up $\rho_{2}: S_2 \to S_1$ at $q_4$. 
Then $S=S_{2}$. Put $\rho=\rho_{1} \rho_{2}: S \to \Bb{P}^2$.
\end{itemize}
Moreover, we have $E_{1}= \rho_{2}^{-1}(q_4)$, 
$E_2 = \rho^{-1}(q_2)$, $E_{3}=\rho^{-1}_{\ast}(\Ov{q_2 q_3})$, $E_4=\rho^{-1}(q_3)$,
$E_5 = \rho^{-1}_{\ast}(\Ov{q_3 q_1})$, $F_1= \rho^{-1}_{\ast}(\Ov{q_1 q_2})$ and 
$F_2 =  (\rho_{2})^{-1}_{\ast}(\rho_{1}^{-1}(q_1))$.
We denote a divisor $D=\sum_{i=1}^{5}a_{i} E_{i} + \sum_{j=1}^{2}b_{j}F_{j} \in \Rm{Div}(S)$ ($a_i , b \in \Z$) by
$D=(a_1,a_2,a_3,a_4,a_5,b_1, b_2)$.
The intersection matrix of $\{E_1, E_2, E_3, E_4, E_5, F_1, F_2 \}$ is 
\[A:=
\left(
\begin{array}{ccccc|cc}
  -1 & 0 & 0  & 0  &  0  &  1  &  1   \\
  0 & -1 & 1  & 0  &  0  &  1  &  0   \\
  0 & 1 & -1  & 1  &  0 &  0 &  0   \\
  0 & 0 & 1 & -1  &  1 &  0  &  0   \\
  0 & 0 & 0  & 1 &  -1  &  0  &  1   \\ \hline
  1 & 1 & 0  & 0  &  0  &  -2  &  0   \\
  1 & 0 & 0  & 0 &  1  &  0  &  -2    \\ 
\end{array}
\right).
\]
 We note that $-K_{S} \sim \sum_{i=1}^{5}E_{i} + \sum_{i=j}^{2}F_{j}=(1,1,1,1,1,1,1)$.
 
%%%%%%%%%%%%%%%%%%%%% (1)

\noindent
(1)\;The case $p \in E_1$.

\noindent
We calculate $S(E_{1})$ and $S(W_{\bullet,\bullet}^{E_{1}},p)$.
Take $u \in \R_{\geq 0}$. 
Let $P(u)+N(u)$ be the Zariski decomposition of $-K_{S}-uE_{1}$.
%where $P(u)$ is the positive part and $N(u)$ is the negative part.
%
If $u\in [0,2]$, then we have
\begin{align*}
&P(u)= \left( 1-u,1 ,1 ,1 ,1 ,1-\frac{u}{2},1-\frac{u}{2} \right) , \\
&N(u)= \left(0,0 ,0 ,0 ,0 ,\frac{u}{2},\frac{u}{2} \right).
\end{align*}
If $u\in [2,3]$, then we have
\begin{align*}
&P(u)= \left( 1-u,3-u ,1 ,1 ,3-u ,2-u,2-u \right), \\
%(1-u)E_1 + (3-u) E_2 + E_{3} + E_4 +  (3-u)E_5 + (2-u)(F_1 + F_2), \\
& N(u)= \left(0,u-2 ,0 ,0 ,u-2 ,u-1,u-1\right).
%(u-2)(E_2 + E_5) + (u-1)(F_1 + F_2).
\end{align*}
We note that $-K_{S}-uE_{1}$ is not pseudo effective for $u>3$.
If $u\in [0,2]$, then we have 
\begin{align*}
P(u)^2=(5-2u) ,\quad P(u)E_1 = 1.
\end{align*}
If $u\in [2,3]$, then we have 
\begin{align*}
P(u)^2=(3-u)^2,\quad P(u)E_1 = (3-u).
\end{align*}
Therefore,  we get 
\begin{align*}
S(E_1)=\frac{1}{5}\int_{0}^{2}(5-2u) du
+ \frac{1}{5}\int_{2}^{3} (3-u)^2 du
=\frac{19}{15}
\end{align*}
by the definition of $S(E_1)$ and 
\begin{align*}
S(W_{\bullet,\bullet}^{E_{1}},p)
&=\frac{2}{5}\int_{0}^{2} \Rm{ord}_{p}(N(u)|_{E_{1}}) du
+\frac{1}{5}\int_{0}^{2} 1 du\\
&+ \frac{2}{5}\int_{2}^{3}(3-u) \Rm{ord}_{p}(N(u)|_{E_{1}})  du
+\frac{1}{5}\int_{2}^{3}(3-u)^2 du\\
 &= 
 \begin{cases}
  \frac{17}{15} & \text{if $p \in E_{1}\cap F_{j} $ for $j=1,2$},\\
                    & \\
   \frac{7}{15}   & \text{if $p \in E_{1} \setminus \bigcup_{j=1}^{2} F_{j}$}, \\
  \end{cases}
\end{align*}
by Definition~\ref{S-inv}.
Hence we have 
\begin{align*}
\frac{15}{19}  \geq \delta_{p}(S) \geq  \Rm{min}\left\{ \frac{1}{S(E_1)}, \frac{1}{S(W_{\bullet,\bullet}^{E_1},p)} \right\}
= \frac{15}{19} 
\end{align*}
from Corollary~\ref{AZ}.
Thus, we have $\delta_{p}(S)= 15/19$ in this case.

%%%%%%%%%%%%%%%%%%%%% (2)

\noindent
(2)\;The case $p \in F_1\setminus E_1$.

\noindent
We calculate $S(F_{1})$ and $S(W_{\bullet,\bullet}^{F_{1}},p)$.
Take $u \in \R_{\geq 0}$. 
Let $P(u)+N(u)$ be the Zariski decomposition of $-K_{S}-uF_{1}$.
%where $P(u)$ is the positive part and $N(u)$ is the negative part.
%
If $u\in [0,1]$, then we have
\begin{align*}
&P(u)=  \left(1,1 ,1 ,1 ,1 ,1-u,1\right),\\
%\sum_{i=1}^5 E_{i} + (1-u)F_1 + F_2 , \\
&N(u)=0.
\end{align*}
If $u\in [1,2]$, then we have
\begin{align*}
&P(u)=  \left(3-2u,2-u ,1 ,1 ,1 ,1-u, 2-u\right),\\
%(3-2u)E_1 + (2-u) E_2 + \sum_{i=3}^{5} E_i + (1-u)F_1 + (2-u)F_2, \\
& N(u)=  \left(2(u-1),u-1 ,0 ,0 ,0 ,0,u-1\right).
%(u-1)(2E_1+ E_2 + F_2).
\end{align*}
We note that $-K_{S}-uF_{1}$ is not pseudo effective for $u>2$.
If $u\in [0,1]$, then we have 
\begin{align*}
P(u)^2=(5-2u^2) ,\quad P(u)F_1 = 2u.
\end{align*}
If $u\in [1,2]$, then we have 
\begin{align*}
P(u)^2=(2-u)(4-u),\quad P(u)F_1 = (3-u).
\end{align*}
Therefore,  we get 
\begin{align*}
S(F_1)=\frac{1}{5}\int_{0}^{1}(5-2u^2) du
+ \frac{1}{5}\int_{1}^{2} (2-u)(4-u) du
=\frac{17}{15}
\end{align*}
by the definition of $S(F_1)$ and 
\begin{align*}
S(W_{\bullet,\bullet}^{F_{1}},p)
%&=\frac{2}{5}\int_{0}^{2} \Rm{ord}_{p}(N(u)|_{F}) du
&=\frac{1}{5}\int_{0}^{1} 4u^2 du
+ \frac{2}{5}\int_{1}^{2}(3-u) \Rm{ord}_{p}(N(u)|_{F_{1}})  du
+\frac{1}{5}\int_{1}^{2}(3-u)^2 du\\
 &= 
 \begin{cases}
  1 & \text{if $p \in F_{1}\cap E_{2} $},\\
                    & \\
   \frac{11}{15}   & \text{if $p \in F_{1} \setminus (E_1 \cup E_{2}) $}, \\
  \end{cases}
\end{align*}
by Definition~\ref{S-inv}.
Hence we have 
\begin{align*}
\frac{15}{17}  \geq \delta_{p}(S) \geq  \Rm{min}\left\{ \frac{1}{S(F_1)}, \frac{1}{S(W_{\bullet,\bullet}^{F_1},p)} \right\}
= \frac{15}{17} 
\end{align*}
from Corollary~\ref{AZ}.
Thus, we have $\delta_{p}(S)= 15/17$ in this case.

%%%%%%%%%%%%%%%%%%%%% (3)

\noindent
(3)\;The case $p \in E_2 \setminus F_1$.

\noindent
We calculate $S(E_{2})$ and $S(W_{\bullet,\bullet}^{E_{2}},p)$.
Take $u \in \R_{\geq 0}$. 
Let $P(u)+N(u)$ be the Zariski decomposition of $-K_{S}-uE_{2}$.
%where $P(u)$ is the positive part and $N(u)$ is the negative part.
%
If $u\in [0,1]$, then we have
\begin{align*}
&P(u)=\left(1,1-u ,1 ,1 ,1 ,1-\frac{u}{2},1\right),\\
%E_1 + (1-u)E_2 + \sum_{i=3}^5 E_{i} + \frac{2-u}{2}F_1 + F_2 , \\
&N(u)=\left(0,0 ,0 ,0 ,0 ,\frac{u}{2},0\right).
%\frac{u}{2}F_1.
\end{align*}
If $u\in [1,2]$, then we have
\begin{align*}
&P(u)=\left(1,1-u ,2-u ,1 ,1 ,1-\frac{u}{2},1\right),\\
%E_1 + (1-u)E_2 + (2-u)E_{3} +E_4 + E_5  + \frac{2-u}{2}F_1 + F_2, \\
& N(u)=\left(0,0 ,u-1 ,0 ,0 ,\frac{u}{2},0\right).
%\frac{u}{2}F_1 + (u-1)E_3
\end{align*}
We note that $-K_{S}-uE_{2}$ is not pseudo effective for $u>2$.
If $u\in [0,1]$, then we have 
\begin{align*}
P(u)^2=5-2u-\frac{u^2}{2} ,\quad P(u)E_2 = \frac{u+2}{2}.
\end{align*}
If $u\in [1,2]$, then we have 
\begin{align*}
P(u)^2=\frac{1}{2}(6-u)(2-u),\quad P(u)E_2 = \frac{4-u}{2}.
\end{align*}
Therefore,  we get 
\begin{align*}
S(E_2)=\frac{1}{5}\int_{0}^{1}(5-2u-\frac{u^2}{2}) du
+ \frac{1}{5}\int_{1}^{2} \frac{1}{2}(6-u)(2-u) du
=1
\end{align*}
by the definition of $S(E_2)$ and 
\begin{align*}
S(W_{\bullet,\bullet}^{E_{2}},p)
%&=\frac{2}{5}\int_{0}^{2} \Rm{ord}_{p}(N(u)|_{F}) du
&=\frac{1}{5}\int_{0}^{1} \frac{(2+u)^2}{4} du
+ \frac{2}{5}\int_{1}^{2} \frac{4-u}{2} \Rm{ord}_{p}(N(u)|_{E_{2}})  du
+\frac{1}{5}\int_{1}^{2}\frac{(4-u)^2}{4} du\\
 &= 
 \begin{cases}
  \frac{13}{15} & \text{if $p \in E_{2}\cap E_{3} $},\\
                  & \\
   \frac{19}{30}   & \text{if $p \in E_{2} \setminus (F_1 \cup E_{3}) $}, \\
  \end{cases}
\end{align*}
by Definition~\ref{S-inv}.
Hence we have 
\begin{align*}
1  \geq \delta_{p}(S) \geq  \Rm{min}\left\{ \frac{1}{S(E_2)}, \frac{1}{S(W_{\bullet,\bullet}^{E_2},p)} \right\}
= 1
\end{align*}
from Corollary~\ref{AZ}.
Thus, we have $\delta_{p}(S)= 1$ in this case.

%%%%%%%%%%%%%%%%%%%%% (4)

\noindent
(4)\;The case $p \in E_3 \setminus E_2$.

\noindent
We calculate $S(E_{3})$ and $S(W_{\bullet,\bullet}^{E_{3}},p)$.
Take $u \in \R_{\geq 0}$. 
Let $P(u)+N(u)$ be the Zariski decomposition of $-K_{S}-uE_{3}$.
%where $P(u)$ is the positive part and $N(u)$ is the negative part.
%
If $u\in [0,1]$, then we have
\begin{align*}
&P(u)=\left(1,1 ,1-u ,1 ,1 ,1,1\right), \\
%\sum_{i=1}^5 E_{i} + (1-u)F_1 + F_2 , \\
&N(u)=0.
\end{align*}
If $u\in [1,2]$, then we have
\begin{align*}
&P(u)=\left(1,3-2u,1-u ,2-u ,1 ,2-u,1\right),\\
%E_1 + (3-2u)E_2 + (1-u)E_3 + (2-u)E_4 + E_5 + (2-u)F_1 + F_2, \\
& N(u)=\left(0,2(u-1) ,0 ,u-1 ,0 ,u-1,0\right).
%(u-1)(2E_2 + E_4 + F_1)
\end{align*}
We note that $-K_{S}-uE_{3}$ is not pseudo effective for $u>2$.
If $u\in [0,1]$ we have 
\begin{align*}
P(u)^2=5-2u-u^2 ,\quad P(u)E_3 = 1+u.
\end{align*}
If $u\in [1,2]$ we have 
\begin{align*}
P(u)^2=2(2-u)^2,\quad P(u)E_3 = 4-2u.
\end{align*}
Therefore,  we get 
\begin{align*}
S(E_3)=\frac{1}{5}\int_{0}^{1}(5-2u-u^2) du
+ \frac{1}{5}\int_{1}^{2} 2(2-u)^2 du
=\frac{13}{15}
\end{align*}
by the definition of $S(E_3)$ and 
\begin{align*}
S(W_{\bullet,\bullet}^{E_{3}},p)
%&=\frac{2}{5}\int_{0}^{2} \Rm{ord}_{p}(N(u)|_{F}) du
&=\frac{1}{5}\int_{0}^{1} (1+u)^2 du
+ \frac{2}{5}\int_{1}^{2}(4-2u) \Rm{ord}_{p}(N(u)|_{E_{3}})  du
+\frac{1}{5}\int_{1}^{2}(4-2u)^2 du\\
 &= 
 \begin{cases}
  \frac{7}{15} + \frac{2}{15} +\frac{4}{15} = \frac{13}{15} &  \text{if $p \in E_{3}\cap E_{4} $},\\
                  & \\
   \frac{7}{15} + \frac{4}{15} = \frac{11}{15}   & \text{if $p \in E_{3} \setminus (E_{2} \cup E_{4}) $}, \\
  \end{cases}
\end{align*}
by Definition~\ref{S-inv}.
Hence we have 
\begin{align*}
\frac{15}{13}  \geq \delta_{p}(S) \geq  \Rm{min}\left\{ \frac{1}{S(E_3)}, \frac{1}{S(W_{\bullet,\bullet}^{E_3},p)} \right\}
= \frac{15}{13}
\end{align*}
from Corollary~\ref{AZ}.
Thus, we have $\delta_{p}(S)= \frac{15}{13}$ in this case.

%%%%%%%%%%%%%%%%%%%%%%%%%%%%%%%%%%%%%(5)
\noindent
(5)\;The case $p \in S \setminus \left( \bigcup_{i,j}( E_{i} \cup F_{j})\right)$.

\noindent
%To apply \Red{Theorem~Abban-Zhang} for prime divisor $F$,
Consider a blowing up $\sigma:\Ti{S} \to S$ at $p$.
Let $E$ be the exceptional curve of $p$,
 $\Ti{E}_{i}$ and $\Ti{F}_{j}$ be the proper transform of $E_{i}$ and $F_{j}$, respectively.
 Put $G_{i}:=(\rho\sigma)_{\ast}^{-1}\Ov{\rho(p)q_i}$ for $i=1,2,3$.
Then we have $\sigma^{\ast}(-K_{S})-uE \sim
\Ti{F}_1 + \Ti{E}_2 +  G_2 + G_3+ (2-u)E$.
We calculate $S(E)$ and $S(W_{\bullet,\bullet}^{E},p)$.
Take $u \in \R_{\geq 0}$. 
Let $\Ti{P}(u)+\Ti{N}(u)$ be the Zariski decomposition of $\sigma^{\ast}(-K_{S})-uE$,
where $\Ti{P}(u)$ is the positive part and $\Ti{N}(u)$ is the negative part.
If $u\in [0,2]$, then we have
\begin{align*}
&\Ti{P}(u)=\Ti{F}_1 + \Ti{E}_2 +  G_2 + G_3+ (2-u)E,\\
&\Ti{N}(u)=0.
\end{align*}
If $u\in [2,\frac{5}{2}]$, then we have
\begin{align*}
&\Ti{P}(u)=\Ti{E}_2+ (3-u)\Ti{F}_1 + (2-u)\Ti{F}_2 + 2(2-u)G_1 + (3-u)G_2+ (5-2u)G_3+ (2-u)E,\\
&\Ti{N}(u)=(u-2)\Ti{F}_1+ (u-2)\Ti{F}_2 + 2(u-2)G_1 + (u-2)G_2+ 2(u-2)G_3.
\end{align*}
We note that $\sigma^{\ast}(-K_{S})-uE$ is not pseudo effective for $u>5/2$.
If $u\in [0,2]$, then we have
\begin{align*}
\Ti{P}(u)^2=5-u^2 ,\quad \Ti{P}(u)E = u.
\end{align*}
If $u\in [2,\frac{5}{2}]$, then we have 
\begin{align*}
\Ti{P}(u)^2=(5-2u)^2,\quad \Ti{P}(u)E = 2(5-2u).
\end{align*}
Therefore, we get 
\begin{align*}
S(E)=\frac{1}{5}\int_{0}^{2} 5 -u^2 du + \frac{1}{5}\int_{2}^{\frac{5}{2}} (5-2u)^2 du =\frac{3}{2}
\end{align*}
by the definition of $S(E)$ and 
\begin{align*}
S(W_{\bullet,\bullet}^{E},p)
&=%\frac{2}{5}\int_{0}^{1}2u \Rm{ord}_{p}(N(u)|_{F}) du
\frac{1}{5}\int_{0}^{2}u^2 du
+ \frac{2}{5}\int_{2}^{\frac{5}{2}}2(5-2u) \Rm{ord}_{p}(\Ti{N}(u)|_{E})  du
+\frac{1}{5}\int_{2}^{\frac{5}{2}}4(5-2u)^2 du\\
 &= 
 \begin{cases}
  \frac{11}{15} & \text{if $p \in E \cap  G_1 $},\\
                    & \\
   \frac{7}{10}   & \text{if $p \in E\cap \left( G_2 \cup G_3 \right)$}, \\
    & \\
   \frac{2}{3}   & \text{if $p \in E \setminus\left(G_1 \cup G_2 \cup G_3\right)$}, \\
  \end{cases}
\end{align*}
by Definition~\ref{S-inv}.
Hence we have 
\begin{align*}
\frac{4}{3}  \geq \delta_{p}(S) \geq  \Rm{min}\left\{ \frac{2}{S(E)}, \frac{1}{S(W_{\bullet,\bullet}^{E},p)} \right\}
= \frac{4}{3} 
\end{align*}
from Corollary~\ref{AZ}.
 Thus, we have $\delta_{p}(S)= 4/3$ in this case.
 \end{proof}

\begin{prop}
Let $S$ be the anti-canonical degree $5$ weak del Pezzo surface such
that the dual graph of negative curves is
\[
\xygraph{
   {\bullet}* +!D{ }*+!U{E_1}
                    -[r]
    {\bullet}*+!U{E_{2}}  - [r]
    {\circ}*+!U{F_{1}}
   %     - [r] \cdots
       - [r] {\circ}* +!U{F_{2}}
       - [r] {\bullet}* +!U{E_{3}}
         - [r] {\circ}*+!U{F_{3}}
              },
        \]
where $E_{i}$ $(i=1,2,3)$ is a $(-1)$-curve and $F_j$ $(j=1,2,3)$ is a $(-2)$-curve.
Then, for a point $p \in S$, it holds that
\begin{equation}
\nonumber
  \delta_{p}(S)=
  \begin{cases}
    \frac{15}{13} & \text{if $p \in E_1 \setminus E_2$,} \\
    \frac{15}{17}                 & \text{if $p \in E_2 \setminus F_1$,} \\
    \frac{15}{19}      & \text{if $p \in F_1 \setminus F_2$,}\\
   \frac{5}{7} & \text{if $p \in F_2 \setminus E_3$,}\\
   \frac{15}{23} &  \text{if $p \in E_3$,}\\
  \frac{15}{17} & \text{if $p \in F_3 \setminus E_3$,}\\
   \frac{30}{23} & \text{if $p \in S \setminus \left( \bigcup_{i,j}( E_{i} \cup F_{j})\right)$.}
  \end{cases}
\end{equation}
\end{prop}

\begin{proof}
We can assume that we get $S$ from $\Bb{P}^2$ as follows.
\begin{itemize}
\item[(1)]
Take two distinct points $q_1,q_4 \in \Bb{P}^2$ and a line $l (\neq \Ov{q_1 q_4})$ passing through $q_1$.
Let $\rho_{1}:S_{1}=\Rm{Bl}_{\{q_{1}, q_{4}\}}\Bb{P}^{2} \to \Bb{P}^2$ be a 
blowing-up at points $q_1$, $q_4$,
let $l_1 = (\rho_1)^{-1}_{\ast}l$ and let $q_2$ be a point at which $l_1$ and $\rho_{1}^{-1}(q_1)$ meet.
\item[(2)] Let $\rho_{2}: S_2 \to S_1$ be a blowing-up at $q_2$,
let $l_2 = (\rho_2)_{\ast}^{-1}l_1$ and let 
$q_3$ be a point at which $l_2$ and $\rho_{2}^{-1}(q_2)$ meet.
\item[(3)] Let $\rho_{3}: S_3 \to S_2$ be a blowing-up at $q_3$.
Then $S=S_3$. Put $\rho=\rho_1 \rho_2 \rho_3$.
\end{itemize}
Moreover, we have $E_{1}=\rho^{-1}(q_4)$, 
$E_2 = \rho^{-1}_{\ast}(\Ov{q_1 q_4})$, 
$F_{1}=(\rho_2 \rho_3)^{-1}_{\ast}(\rho_1^{-1}(q_1))$, 
$F_2= (\rho_3)^{-1}_{\ast}(\rho^{-1}(q_2))$,
$E_3 =\rho_{3}^{-1}(q_3)$, 
$F_3 = \rho^{-1}_{\ast}l$.
We denote $D=\sum_{i=1}^{3}a_{i} E_{i} + \sum_{j=1}^{3}b_{j}F_{j} \in \Rm{Div}(S)$ ($a_i , b \in \Z$)
by $D=(a_1,a_2,a_3, b_1, b_2 ,  b_3)$.
The intersection matrix of $\{E_1, E_2, E_3, F_1, F_2 , F_3 \}$ is 
\[A:=
\left(
\begin{array}{ccc|ccc}
  -1 & 1 & 0  & 0  &  0  &  0     \\
  1 & -1 & 0  & 1  &  0  &  0     \\
  0 & 0 & -1  & 0  &  1 &  1    \\ \hline
  0 & 1 & 0 & -2  &  1 &  0     \\
  0 & 0 & 1  & 1 &  -2  &  0     \\ 
  0 & 0 & 1  & 0  &  0  &  -2    \\
\end{array}
\right).
\]
 We note that $-K_{S} \sim 2E_1+ 3E_2 + 2F_1 + F_2 = (2,3,0,2,1,0)$.
%%%%%%%%%%%%%%%%%%%%% (1)

\noindent
(1)\;The case $p \in E_1 \setminus E_2$.

\noindent
We calculate $S(E_{1})$ and $S(W_{\bullet,\bullet}^{E_{1}},p)$.
Take $u \in \R_{\geq 0}$. 
Let $P(u)+N(u)$ be the Zariski decomposition of $-K_{S}-uE_{1}$.
%where $P(u)$ is the positive part and $N(u)$ is the negative part.
%
If $u\in [0,1]$, then we have
\begin{align*}
&P(u)=\left(2-u,3 ,0 ,2 ,1 ,0\right),\\
%(2-u)E_1 + 3E_2 + 2F_1 + F_2, \\
&N(u)=0.
\end{align*}
If $u\in [1,2]$, then we have
\begin{align*}
&P(u)= \left(2-u,3(2-u) ,0 ,2(2-u) ,2-u ,0\right),\\
%(2-u)E_1 + (4-3u)E_2 + (3-2u)F_1 + (2-u)F_2, \\
& N(u)=\left(0,3(u-1) ,0 ,2(u-1) ,u-1 ,0\right).
%3(u-1)E_2+2(u-1)F_1 + (u-1)F_2.
\end{align*}
%(\Red{多分$P(u)= (2-u)E_1 + 3(2-u)E_2 + 2(2-u)F_1 + (2-u)F_2,$である})
We note that $-K_{S}-uE_{1}$ is not pseudo effective for $u>2$.
If $u\in [0,1]$, then we have 
\begin{align*}
P(u)^2=5-2u-u^2 ,\quad P(u)E_1 = 1+u.
\end{align*}
If $u\in [1,2]$, then we have 
\begin{align*}
P(u)^2=2(2-u)^2,\quad P(u)E_1 = 4-2u.
\end{align*}
Therefore,  we get 
\begin{align*}
S(E_1)=\frac{1}{5}\int_{0}^{1}(5-2u-u^2) du
+ \frac{1}{5}\int_{1}^{2} 2(2-u)^2 du
=\frac{13}{15}
\end{align*}
by the definition of $S(E_1)$ and 
\begin{align*}
S(W_{\bullet,\bullet}^{E_{1}},p)
%=\frac{2}{5}\int_{0}^{1} \Rm{ord}_{p}(N(u)|_{F}) du
&=\frac{1}{5}\int_{0}^{1} (1+u)^2 du
+ \frac{2}{5}\int_{1}^{2}(4-2u) \Rm{ord}_{p}(N(u)|_{E_{1}})  du
+\frac{1}{5}\int_{1}^{2}(4-2u)^2 du\\
 &= \frac{7}{15}+\frac{4}{15} = \frac{11}{15}
% \begin{cases}
 % \frac{17}{15} & \text{if $p \in E_{1}\cap F_{j} $ for $j=1,2$},\\
 %                   & \\
 %  \frac{7}{15}   & \text{if $p \in E_{1} \setminus \bigcup_{j=1}^{2} F_{j}$}, \\
%  \end{cases}
\end{align*}
by Definition~\ref{S-inv}.
Hence we have 
\begin{align*}
\frac{15}{13}  \geq \delta_{p}(S) \geq  \Rm{min}\left\{ \frac{1}{S(E_1)}, \frac{1}{S(W_{\bullet,\bullet}^{E_1},p)} \right\}
= \frac{15}{13} 
\end{align*}
from Corollary~\ref{AZ}.
Thus, we have $\delta_{p}(S)= 15/13$ in this case.

%%%%%%%%%%%%%%%%%%%%% (2)

\noindent
(2)\;The case $p \in E_2 \setminus F_1$.

\noindent
We calculate $S(E_{2})$ and $S(W_{\bullet,\bullet}^{E_{2}},p)$.
Take $u \in \R_{\geq 0}$. 
Let $P(u)+N(u)$ be the Zariski decomposition of $-K_{S}-uE_{2}$.
%where $P(u)$ is the positive part and $N(u)$ is the negative part.
%
If $u\in [0,1]$, then we have
\begin{align*}
&P(u)=\left(2,3-u ,0 ,2-\frac{2}{3}u ,1-\frac{u}{3} ,0\right),\\
%2E_1 + (3-u)E_2 + (2-\frac{2}{3}u)F_1 + (1-\frac{1}{3}u) F_2, \\
&N(u)=\left(0,0 ,0 ,\frac{2}{3}u ,\frac{u}{3} ,0\right).
%\frac{2}{3}uF_1 + \frac{1}{3}uF_2.
\end{align*}
If $u\in [1,3]$, then we have
\begin{align*}
&P(u)= \left(3-u,3-u ,0 ,2-\frac{2}{3}u ,1-\frac{u}{3} ,0\right), \\
%(3-u)E_1 + (3-u)E_2 + (2-\frac{2}{3}u)F_1 + (1-\frac{1}{3}u) F_2, \\
& N(u)= \left(u-1,0 ,0 ,\frac{2}{3}u ,\frac{u}{3} ,0\right).
%(u-1)E_1+ \frac{2}{3}uF_1 + \frac{1}{3}uF_2.
\end{align*}
We note that $-K_{S}-uE_{2}$ is not pseudo effective for $u>3$.
If $u\in [0,1]$, then we have 
\begin{align*}
P(u)^2=5-2u-\frac{u^2}{3} ,\quad P(u)E_2 =  \frac{3+u}{3}.
\end{align*}
If $u\in [1,3]$, then we have 
\begin{align*}
P(u)^2=\frac{2}{3}(3-u)^2,\quad P(u)E_2 =  2-\frac{2}{3}u.
\end{align*}
Therefore,  we get 
\begin{align*}
S(E_2)=\frac{1}{5}\int_{0}^{1}(5-2u-\frac{u^2}{3}) du
+ \frac{1}{5}\int_{1}^{3} \frac{2}{3}(3-u)^2 du
=\frac{17}{15}
\end{align*}
by the definition of $S(E_2)$ and 
\begin{align*}
S(W_{\bullet,\bullet}^{E_{2}},p)
&=\frac{2}{5}\int_{0}^{1} \frac{3+u}{3}\Rm{ord}_{p}(N(u)|_{E_2}) du
+\frac{1}{5}\int_{0}^{1} \frac{(3+u)^2}{9}du \\
&+ \frac{2}{5}\int_{1}^{3}\left( 2-\frac{2}{3}u\right)\Rm{ord}_{p}(N(u)|_{E_{2}})  du
+\frac{1}{5}\int_{1}^{3} \left(2-\frac{2}{3}u\right)^2 du\\
 &=
 \begin{cases}
  \frac{13}{15} & \text{if $p \in E_{2}\cap E_{1}$},\\
                    & \\
   \frac{23}{45}   & \text{if $p \in E_{2} \setminus (E_1 \cup F_1) $}, \\
 \end{cases}
\end{align*}
by Definition~\ref{S-inv}.
Hence we have 
\begin{align*}
\frac{15}{17}  \geq \delta_{p}(S) \geq  \Rm{min}\left\{ \frac{1}{S(E_2)}, \frac{1}{S(W_{\bullet,\bullet}^{E_2},p)} \right\}
= \frac{15}{17} 
\end{align*}
from Corollary~\ref{AZ}.
Thus, we have $\delta_{p}(S)= 15/17$ in this case.

%%%%%%%%%%%%%%%%%%%%% (3)

\noindent
(3)\;The case $p \in F_1 \setminus F_2$.

\noindent
We calculate $S(F_{1})$ and $S(W_{\bullet,\bullet}^{F_{1}},p)$.
Take $u \in \R_{\geq 0}$. 
Let $P(u)+N(u)$ be the Zariski decomposition of $-K_{S}-uF_{1}$.
%where $P(u)$ is the positive part and $N(u)$ is the negative part.
%
If $u\in [0,1]$, then we have
\begin{align*}
&P(u)= \left(2,3 ,0 ,2-u ,1-\frac{u}{2} ,0\right),\\
%2E_1 + 3E_2 + (2-u)F_1 + (1-\frac{1}{2}u) F_2, \\
&N(u)=\left(0,0 ,0 ,0 ,\frac{u}{2} ,0\right).
% \frac{1}{2}uF_2.
\end{align*}
If $u\in [1,2]$, then we have
\begin{align*}
&P(u)= \left(2,4-u ,0 ,2-u ,1-\frac{u}{2} ,0\right),\\
%2E_1 + (4-u) E_2 + (2-u)F_1 + (1-\frac{1}{2}u) F_2, \\
& N(u)= \left(0,u-1 ,0 , 0 , \frac{u}{2} ,0\right).
%(u-1)E_2+ \frac{1}{2}uF_2 .
\end{align*}
We note that $-K_{S}-uF_{1}$ is not pseudo effective for $u>2$.
If $u\in [0,1]$, then we have 
\begin{align*}
P(u)^2=5-\frac{3}{2}u^2 ,\quad P(u)F_1 =  \frac{3u}{2}.
\end{align*}
If $u\in [1,2]$, then we have 
\begin{align*}
P(u)^2= \frac{1}{2}(2-u)(6+u),\quad P(u)F_1 =  \frac{2+u}{2}.
\end{align*}
Therefore,  we get 
\begin{align*}
S(F_1)=\frac{1}{5}\int_{0}^{1} \left(5-\frac{3}{2}u^2 \right)du
+ \frac{1}{5}\int_{1}^{2} \frac{1}{2}(2-u)(6+u) du
=\frac{19}{15}
\end{align*}
by the definition of $S(F_1)$ and 
\begin{align*}
S(W_{\bullet,\bullet}^{F_{1}},p)
%&=\frac{2}{5}\int_{0}^{1} \Rm{ord}_{p}(N(u)|_{F_1}) du
&=\frac{1}{5}\int_{0}^{1}  \frac{9}{4}u^2 du 
+ \frac{2}{5}\int_{1}^{2}  \frac{2+u}{2} \Rm{ord}_{p}(N(u)|_{F_1})  du
+\frac{1}{5}\int_{1}^{2} \frac{(2+u)^2}{4} du\\
 &=
 \begin{cases}
  \frac{17}{15} & \text{if $p \in F_{1}\cap E_{2}$},\\
                    & \\
   \frac{23}{30}   & \text{if $p \in F_{1} \setminus (E_2 \cup F_2) $}, \\
 \end{cases}
\end{align*}
by Definition~\ref{S-inv}.
Hence we have 
\begin{align*}
\frac{15}{19}  \geq \delta_{p}(S) \geq  \Rm{min}\left\{ \frac{1}{S(F_1)}, \frac{1}{S(W_{\bullet,\bullet}^{F_1},p)} \right\}
= \frac{15}{19} 
\end{align*}
from Corollary~\ref{AZ}.
Thus, we have $\delta_{p}(S)= 15/19$ in this case.

%%%%%%%%%%%%%%%%%%%%% (4)

\noindent
(4)\;The case $p \in F_2 \setminus E_3$.

\noindent
We calculate $S(F_{2})$ and $S(W_{\bullet,\bullet}^{F_{2}},p)$.
Take $u \in \R_{\geq 0}$. 
Let $P(u)+N(u)$ be the Zariski decomposition of $-K_{S}-uF_{2}$.
%where $P(u)$ is the positive part and $N(u)$ is the negative part.
%
If $u\in [0,1]$, then we have
\begin{align*}
&P(u)= \left(2, 3 ,0 ,2-\frac{u}{2} ,1-u, 0\right),\\
%2E_1 + 3E_2 + \left(2-\frac{1}{2}u\right)F_1 + (1-u) F_2, \\
&N(u)= \left(0, 0 , 0 ,\frac{u}{2} ,0, 0\right).
%\frac{1}{2}u F_1.
\end{align*}
If $u\in [1,2]$, then we have
\begin{align*}
&P(u)=\left(2,3 ,2(1-u) ,2-\frac{u}{2} ,1-u, 1-u \right),\\
% 2E_1 + 3E_2 + \left(2-\frac{1}{2}u\right)F_1 + (1-u) F_2 + 2(1-u)E_3 + (1-u)F_3, \\
& N(u)=\left(0, 0 ,2(u-1) ,\frac{u}{2} ,0, u-1\right).
%\frac{1}{2}u F_1+2(u-1)E_3 + (u-1)F_3 .
\end{align*}
If $u\in [2,3]$, then we have
\begin{align*}
&P(u)= \left(2, 5-u ,2(1-u) ,3-u ,1-u, 1-u\right),\\
%2E_1 + (5-u)E_2 + (3-u)F_1 + (1-u) F_2 + 2(1-u)E_3 + (1-u)F_3, \\
& N(u)= \left(0, u-2 ,2(u-1) ,u-1 ,0, u-1\right).
%(u-2)E_2 +  (u-1) F_1+2(u-1)E_3 + (u-1)F_3 .
\end{align*}
We note that $-K_{S}-uF_{2}$ is not pseudo effective for $u>3$.
If $u\in [0,1]$, then we have
\begin{align*}
P(u)^2=\frac{1}{2}(10-3u^2) ,\quad P(u)F_2 =  \frac{3u}{2}.
\end{align*}
If $u\in [1,2]$, then we have 
\begin{align*}
P(u)^2= \frac{1}{2}(u^2 -8u +14),\quad P(u)F_2 = 2-\frac{u}{2}.
\end{align*}
If $u\in [2,3]$, then we have 
\begin{align*}
P(u)^2= (3-u)^2,\quad P(u)F_2 =  3-u.
\end{align*}
Therefore,  we get 
\begin{align*}
S(F_2)=\frac{1}{5}\int_{0}^{1} \frac{1}{2}(10-3u^2) du
+ \frac{1}{5}\int_{1}^{2} \frac{1}{2}(u^2 -8u +14) du
+ \frac{1}{5}\int_{2}^{3} (3-u)^2 du
=\frac{7}{5}
\end{align*}
by the definition of $S(F_2)$ and 
\begin{align*}
S(W_{\bullet,\bullet}^{F_{2}},p)
&=\frac{2}{5}\int_{0}^{1}\frac{3u}{2} \Rm{ord}_{p}(N(u)|_{F_2}) du
+\frac{1}{5}\int_{0}^{1}  \frac{9}{4}u^2 du \\
&+ \frac{2}{5}\int_{1}^{2}  \left(2-\frac{u}{2}\right) \Rm{ord}_{p}(N(u)|_{F_2})  du
+\frac{1}{5}\int_{1}^{2}  \left(2-\frac{u}{2}\right)^2 du\\
&+ \frac{2}{5}\int_{2}^{3}  (3-u)\Rm{ord}_{p}(N(u)|_{F_2})  du
+\frac{1}{5}\int_{2}^{3} (3-u)^2 du\\
 &=
 \begin{cases}
  \frac{23}{30} & \text{if $p \in F_{2}\cap F_{1}$},\\
                    & \\
   \frac{8}{15}   & \text{if $p \in F_{2} \setminus (E_3 \cup F_1) $}, \\
 \end{cases}
\end{align*}
by Definition~\ref{S-inv}.
Hence we have 
\begin{align*}
\frac{5}{7}  \geq \delta_{p}(S) \geq  \Rm{min}\left\{ \frac{1}{S(F_2)}, \frac{1}{S(W_{\bullet,\bullet}^{F_2},p)} \right\}
= \frac{5}{7} 
\end{align*}
from Corollary~\ref{AZ}.
Thus, we have $\delta_{p}(S)= 5/7$ in this case.

%%%%%%%%%%%%%%%%%%%%% (5)

\noindent
(5)\;The case $p \in E_3 $.

\noindent
We calculate $S(E_{3})$ and $S(W_{\bullet,\bullet}^{E_{3}},p)$.
Take $u \in \R_{\geq 0}$. 
Let $P(u)+N(u)$ be the Zariski decomposition of $-K_{S}-uE_{3}$.
%where $P(u)$ is the positive part and $N(u)$ is the negative part.
%
If $u\in [0,3]$, then we have
\begin{align*}
&P(u)= \left(2, 3 ,-u ,2-\frac{u}{3} ,1-\frac{2u}{3}, -\frac{u}{2}\right),\\
 %2E_1 + 3E_2 + \left(2-\frac{1}{3}u \right)F_1 + \left(1-\frac{2}{3}u \right)F_2 -uE_3 -\frac{u}{2} F_3\\
&N(u)= \left(0, 0 ,0 ,\frac{u}{3} ,\frac{2u}{3}, \frac{u}{2}\right).
%\frac{1}{3}u F_1 + \frac{2}{3}u F_2 + \frac{u}{2}F_3.
\end{align*}
If $u\in [3,4]$, then we have
\begin{align*}
&P(u)=\left(2, 6-u ,-u ,4-u ,2-u, -\frac{u}{2}\right),\\
 %2E_1+ (6-u)E_2 -uF_1 + (2-u)F_2 + -uE_3 -\frac{u}{2}F_3, \\
& N(u)=  \left(0, u-3 ,0 ,u-2 ,u-1 ,\frac{u}{2}\right).
%(u-3)E_2 + (u-2)F_1 + (u-1)F_2 + \frac{u}{2}F_3  .
\end{align*}
We note that $-K_{S}-uE_{3}$ is not pseudo effective for $u>4$.
If $u\in [0,3]$, then we have
\begin{align*}
P(u)^2=5-2u+\frac{u^2}{6} ,\quad P(u)E_3 =  \frac{6-u}{6}.
\end{align*}
If $u\in [3,4]$, then we have
\begin{align*}
P(u)^2=  \frac{1}{2}(4-u)^2,\quad P(u)E_3 =  \frac{4-u}{2}.
\end{align*}
Therefore,  we get 
\begin{align*}
S(E_3)=\frac{1}{5}\int_{0}^{3} \left(5-2u+\frac{u^2}{6}\right) du
+ \frac{1}{5}\int_{3}^{4} \frac{1}{2}(4-u)^2 du
=\frac{23}{15}
\end{align*}
by the definition of $S(E_3)$ and 
\begin{align*}
S(W_{\bullet,\bullet}^{E_{3}},p)
&=\frac{2}{5}\int_{0}^{3}\frac{6-u}{6} \Rm{ord}_{p}(N(u)|_{E_3}) du
+\frac{1}{5}\int_{0}^{3}  \left(\frac{6-u}{6}\right)^2 du \\
&+ \frac{2}{5}\int_{3}^{4}  \frac{4-u}{2} \Rm{ord}_{p}(N(u)|_{E_3})  du
+\frac{1}{5}\int_{3}^{4}   \left(\frac{4-u}{2}\right)^2 du\\
 &=
 \begin{cases}
  \frac{7}{5} & \text{if $p \in E_{3}\cap F_{2}$},\\
                    & \\
   \frac{17}{15}   &   \text{if $p \in E_{3}\cap F_{3}$},\\
                    & \\
    \frac{11}{30} & \text{if $p \in E_{3} \setminus (F_2 \cup F_3) $},
 \end{cases}
\end{align*}
by Definition~\ref{S-inv}.
Hence we have 
\begin{align*}
\frac{15}{23}  \geq \delta_{p}(S) \geq  \Rm{min}\left\{ \frac{1}{S(E_3)}, \frac{1}{S(W_{\bullet,\bullet}^{E_3},p)} \right\}
= \frac{15}{23} 
\end{align*}
from Corollary~\ref{AZ}.
Thus, we have $\delta_{p}(S)= 15/23$ in this case.

%%%%%%%%%%%%%%%%%%%%% (6)

\noindent
(6)\;The case $p \in F_3 \setminus E_3$.

\noindent
We calculate $S(F_{3})$ and $S(W_{\bullet,\bullet}^{F_{3}},p)$.
Take $u \in \R_{\geq 0}$. 
Let $P(u)+N(u)$ be the Zariski decomposition of $-K_{S}-uF_{3}$.
%where $P(u)$ is the positive part and $N(u)$ is the negative part.
%
If $u\in [0,1]$, then we have
\begin{align*}
&P(u)= \left(2, 3 ,0 ,2 ,1 , -u\right),\\
%2E_1 + 3E_2 + 2F_1 + F_2 -u F_3\\
&N(u)=0.
\end{align*}
If $u\in [1,2]$, then we have
\begin{align*}
&P(u)= \left(2, 3 ,3-3u ,3-u ,3-2u, -u\right),\\
%2E_1 + 3E_2 + (3-u)F_1 + (3-2u)F_2 + (3-3u)E_3 -uF_3, \\
& N(u)= \left(0, 0 ,3(u-1) ,u-1, 2(u-1), 0\right).
%(u-1)(F_1 + 2F_2 + 3E_3 )  .
\end{align*}
We note that $-K_{S}-uF_{3}$ is not pseudo effective for $u>2$.
If $u\in [0,1]$, then we have 
\begin{align*}
P(u)^2=5-2u^2 ,\quad P(u)F_3 =  2u.
\end{align*}
If $u\in [1,2]$, then we have 
\begin{align*}
P(u)^2=  (4-u)(2-u),\quad P(u)F_3 =  3-u.
\end{align*}
Therefore,  we get 
\begin{align*}
S(F_3)=\frac{1}{5}\int_{0}^{1} \left(5-2u^2\right) du
+ \frac{1}{5}\int_{1}^{2} (4-u)(2-u) du
=\frac{17}{15}
\end{align*}
by the definition of $S(F_3)$ and 
\begin{align*}
S(W_{\bullet,\bullet}^{F_{3}},p)
%&=\frac{2}{5}\int_{0}^{3}\frac{6-u}{u} \Rm{ord}_{p}(N(u)|_{F_3}) du
&=\frac{1}{5}\int_{0}^{1} 4u^2 du 
+ \frac{2}{5}\int_{1}^{2}  (3-u) \Rm{ord}_{p}(N(u)|_{F_3})  du
+\frac{1}{5}\int_{1}^{2}   \left(3-u\right)^2 du
 =\frac{17}{15}
 %\begin{cases}
 % \frac{7}{5} & \text{if $p \in E_{3}\cap F_{2}$},\\
  %                  & \\
  % \frac{17}{15}   &   \text{if $p \in E_{3}\cap F_{3}$}.
%                  & \\
%    \frac{11}{30} & \text{if $p \in E_{3} \setminus (F_2 \cup F_3) $},
 %\end{cases}
\end{align*}
by Theorem~\ref{S-inv}.
Hence we have 
\begin{align*}
\frac{15}{17}  \geq \delta_{p}(S) \geq  \Rm{min}\left\{ \frac{1}{S(F_3)}, \frac{1}{S(W_{\bullet,\bullet}^{F_3},p)} \right\}
= \frac{15}{17} 
\end{align*}
from Corollary~\ref{AZ}.
Thus, we have $\delta_{p}(S)= 15/17$ in this case.

%%%%%%%%%%%%%%%%%%%%%%%%%%%%%%%%%%%%%(7)
\noindent
(7)\;The case $p \in S \setminus \left( \bigcup_{i,j}( E_{i} \cup F_{j})\right)$.

\noindent
%To apply \Red{Theorem~Abban-Zhang} for prime divisor $F$,
Consider a blowing up $\sigma:\Ti{S} \to S$ at $p$.
Let $E$ be the exceptional curve of $p$,
 $\Ti{E}_{i}$ and $\Ti{F}_{j}$ be the proper transform of $E_{i}$ and $F_{j}$, respectively.
Take two $(-1)$-curves $G_1 := (\rho\sigma)_{\ast}^{-1}(\Ov{\rho\sigma(p)q_4})$ 
and $G_2 := (\rho\sigma)_{\ast}^{-1}(\Ov{\rho\sigma(p)q_1})$ on $\Ti{S}$.
Since $\Ov{\rho\sigma(p)q_4}+\Ov{\rho\sigma(p)q_1}+l \in |-K_{\Bb{P}^2}|$,
%where $l \in |H|$ is the line such that $F_3 = \rho_{-1}^{\ast}l$.
we have 
\[
\sigma^{\ast}(-K_{S})-uE \sim
\Ti{E}_3 + \Ti{F}_1 + \Ti{F}_2 +  \Ti{F}_3  + G_1 + G_2 + (2-u)E.
\]
We calculate $S(E)$ and $S(W_{\bullet,\bullet}^{E},p)$.
Take $u \in \R_{\geq 0}$. 
Let $\Ti{P}(u)+\Ti{N}(u)$ be the Zariski decomposition of $\sigma^{\ast}(-K_{S})-uE$,
where $\Ti{P}(u)$ is the positive part and $\Ti{N}(u)$ is the negative part.
If $u\in [0,2]$, then we have
\begin{align*}
&\Ti{P}(u)=\Ti{E}_3  + \Ti{F}_1 + \Ti{F}_2 +  \Ti{F}_3 + G_1 + G_2 + (2-u)E,\\
&\Ti{N}(u)=0.
\end{align*}
If $u\in [2,\frac{5}{2}]$, then we have
\begin{align*}
&\Ti{P}(u)=\Ti{E}_3  + (5-2u)\Ti{F}_1 + (3-u)(\Ti{F}_2 +  \Ti{F}_3) + (5-2u)G_1 + (7-3u)G_2 + (2-u)E,\\
&\Ti{N}(u)=(u-2)(2F_1+ F_2 + F_3 + 2G_1 + 3G_2).
\end{align*}
We note that $\sigma^{\ast}(-K_{S})-uE$ is not pseudo effective for $u>5/2$.
If $u\in [0,2]$, then we have
\begin{align*}
\Ti{P}(u)^2=5-u^2 ,\quad \Ti{P}(u)E = u.
\end{align*}
If $u\in [2,\frac{5}{2}]$, then we have
\begin{align*}
\Ti{P}(u)^2=(5-2u)^2,\quad \Ti{P}(u)E = 2(5-2u).
\end{align*}
Therefore, we get 
\begin{align*}
S(E)=\frac{1}{5}\int_{0}^{2} 5 -u^2 du + \frac{1}{5}\int_{2}^{\frac{5}{2}} (5-2u)^2 du =\frac{3}{2}
\end{align*}
by the definition of $S(E)$ and 
\begin{align*}
S(W_{\bullet,\bullet}^{E},p)
&=%\frac{2}{5}\int_{0}^{1}2u \Rm{ord}_{p}(N(u)|_{F}) du
\frac{1}{5}\int_{0}^{2}u^2 du
+ \frac{2}{5}\int_{2}^{\frac{5}{2}}2(5-2u) \Rm{ord}_{p}(\Ti{N}(u)|_{E})  du
+\frac{1}{5}\int_{2}^{\frac{5}{2}}4(5-2u)^2 du\\
 &= 
 \begin{cases}
  \frac{11}{15} & \text{if $p \in E \cap  G_1 $},\\
                    & \\
   \frac{23}{30}   & \text{if $p \in E\cap G_2 $}, \\
    & \\
   \frac{2}{3}   & \text{if $p \in E \setminus\left(G_1 \cup G_2 \right)$}, \\
  \end{cases}
\end{align*}
by Definition~\ref{S-inv}.
Hence we have 
\begin{align*}
\frac{4}{3}  \geq \delta_{p}(S) \geq  \Rm{min}\left\{ \frac{2}{S(E)}, \frac{1}{S(W_{\bullet,\bullet}^{E},p)} \right\}
= \frac{30}{23} 
\end{align*}
from Corollary~\ref{AZ}.

We also calculate $S(G_2)$.
Take $u \in \R_{\geq 0}$. 
Let $\Ti{P}(u)+\Ti{N}(u)$ be the Zariski decomposition of $\sigma^{\ast}(-K_{S})-uG_2$,
where $\Ti{P}(u)$ is the positive part and $\Ti{N}(u)$ is the negative part.
If $u\in [0,\frac{3}{2}]$, then we have
\begin{align*}
&\Ti{P}(u)=\Ti{E}_3  +\left(1-\frac{2}{3}u \right) \Ti{F}_1 + \left(1-\frac{1}{3}u \right)\Ti{F}_2 +  \Ti{F}_3 + G_1 + (1-u)G_2 + (2-u)E,\\
&\Ti{N}(u)=\frac{2}{3}u  \Ti{F}_1 +  \frac{1}{3}u \Ti{F}_2+ uE.
\end{align*}
If $u\in [\frac{3}{2},2]$, then we have
\begin{align*}
&\Ti{P}(u)=(3-2u)\Ti{E}_2+\Ti{E}_3  + (3-2u)\Ti{F}_1 + (2-u)\Ti{F}_2 +  \Ti{F}_3 + G_1 + (1-u)G_2 + (2-u)E,\\
&\Ti{N}(u)=(2u-3)\Ti{E}_2+ 2(u-1)\Ti{F}_1 + (u-1)\Ti{F}_2 + uE .
\end{align*}
We note that $\sigma^{\ast}(-K_{S})-uG_2$ is not pseudo effective for $u>2$.
If $u\in [0,\frac{3}{2}]$, then we have
\begin{align*}
\Ti{P}(u)^2=5-4u+\frac{2u^2}{3}.
\end{align*}
If $u\in [\frac{3}{2},2]$, then we have
\begin{align*}
\Ti{P}(u)^2=2(2-u)^2.
\end{align*}
Therefore, we get 
\begin{align*}
S(G_2)=\frac{1}{5}\int_{0}^{\frac{3}{2}}5-4u+\frac{2u^2}{3} du + \frac{1}{5}\int_{\frac{3}{2}}^{2} 2(2-u)^2 du =\frac{23}{30}
\end{align*}
by the definition of $S(G_2)$.
Hence we have $\frac{30}{23}  \geq \delta_{p}(S)$.
Therefore, we get
\begin{align*}
 \delta_{p}(S) =   \frac{30}{23}.
\end{align*}

\end{proof}

\begin{prop}
Let $S$ be the anti-canonical degree $5$ weak del Pezzo surface such
that the dual graph of negative curves is
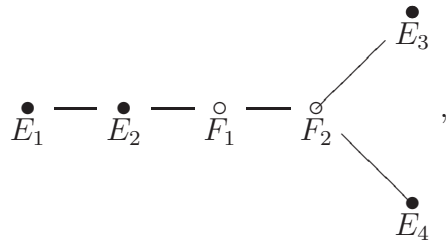
\[
\xygraph{
    {\bullet}* +!U{E_{1}} - [r]
    {\bullet}*+!U{E_{2}}  -[r]
    {\circ}* +!U{F_{1}} - [r]
   {\circ}* +!D{ }*+!U{F_{2}}
      (- [ru] \bullet 
    {\bullet}* +!U{E_{3}} 
   % {\bullet}*+!U{E_{4}}  
   %     - [r] \cdots
   %     - [r] \bullet 
   %     - [r] \bullet 
     %  - [rd] \bullet
              ,
                -[rd]
    {\bullet}*+!U{E_{4}} 
   % {\bullet}*+!U{E_{6}}
   %     - [r] \cdots
  %      - [r] \bullet 
  %      - [r] \bullet 
    %     - [ru] {\bullet}*+!U{E_{7}}
         },
        \]

where $E_{i}$ $(i=1, \cdots, 4)$ is a $(-1)$-curve and $F_j$ $(j=1, 2)$ is a $(-2)$-curve.
Then, for a point $p \in S$, it holds that
\begin{equation}
\nonumber
  \delta_{p}(S)=
  \begin{cases}
    \frac{15}{13} & \text{if $p \in E_1 \setminus E_2$,} \\
    \frac{15}{17}                 & \text{if $p \in E_{2}\setminus F_1$,} \\
    \frac{15}{19}       & \text{if $p \in F_{1}\setminus F_2$,} \\
   \frac{5}{7} & \text{if $p \in F_{2}$,}\\
   \frac{30}{31}       & \text{if $p \in E_{i}\setminus F_2 $ for $i=3,4$,}\\
   \frac{30}{23} & \text{if $p \in S \setminus \left( \bigcup_{i,j}( E_{i} \cup F_{j})\right)$.}
  \end{cases}
\end{equation}
\end{prop}

\begin{proof}
We can assume that we get $S$ from $\Bb{P}^2$ as follows.
\begin{itemize}
\item[(1)]
Take three distinct co-linear points $q_1, q_3, q_4 \in \Bb{P}^2$ and a line $l (\neq \Ov{q_1q_3})$ passing through $q_1$.
Let $\rho_{1}:S_{1}=\Rm{Bl}_{\{q_{1}, q_3,  q_{4}\}}\Bb{P}^{2} \to \Bb{P}^2$ be a 
blowing-up at points $q_1$, $q_3$, $q_4$,
and let $q_2 \in S_1$ be a point at which of $(\rho_{1})^{-1}_{\ast}l$ and $\rho_{1}^{-1}(q_1)$ meet.
% and let $l_1$ be the proper transform of $l$ by $\rho_1$.
\item[(2)] Let $\rho_{2}: S_2 \to S_1$ be a blowing-up at $q_2$.
Then $S=S_2$. Put $\rho=\rho_1 \rho_2 $.
\end{itemize}
Moreover, we have $E_{1}=\rho_{\ast}^{-1}l$, 
$E_2 = \rho_{2}^{-1}(q_2)$, 
$F_{1}=(\rho_2)^{-1}_{\ast}(\rho_1^{-1}(q_1))$, 
$F_2= (\rho)^{-1}_{\ast}(\Ov{q_1 q_3})$,
$E_3 =\rho^{-1}(q_3)$, 
$E_4 =\rho^{-1}(q_4)$.
We denote $D=\sum_{i=1}^{4}a_{i} E_{i} + \sum_{j=1}^{2}b_{j}F_{j} \in \Rm{Div}(S)$ ($a_i , b \in \Z$)
by $D=(a_1,a_2,a_3,a_4, b_1, b_2)$.
The intersection matrix of $\{E_1, E_2, E_3, E_4, F_1, F_2  \}$ is 
\[A:=
\left(
\begin{array}{cccc|cc}
  -1 & 1 & 0  & 0  &  0  &  0     \\
  1 & -1 & 0  & 0  &  1  &  0     \\
  0 & 0 & -1  & 0  &  0 &   1    \\ 
  0 & 0 & 0 & -1  &  0 &    1     \\ \hline
  0 & 1 & 0  & 0 &  -2  &   1     \\ 
  0 & 0 & 1  & 1  &  1  &  -2    \\
\end{array}
\right).
\]
 We note that 
 \[
 -K_{S} \sim 2E_1+ 3E_2 + 2F_1 + F_2 = \left(2, 3 ,0 ,0 , 2, 1\right).
 %\sim E_2 + 2F_1 + 3F_2 + 2E_3 + 2E_4.
\]
%%%%%%%%%%%%%%%%%%%%% (1)

\noindent
(1)\;The case $p \in E_1 \setminus E_2$.

\noindent
We calculate $S(E_{1})$ and $S(W_{\bullet,\bullet}^{E_{1}},p)$.
Take $u \in \R_{\geq 0}$. 
Let $P(u)+N(u)$ be the Zariski decomposition of $-K_{S}-uE_{1}$.
%where $P(u)$ is the positive part and $N(u)$ is the negative part.
%
If $u\in [0,1]$, then we have
\begin{align*}
&P(u)= \left(2-u, 3 ,0 ,0 , 2, 1\right),\\
%(2-u)E_1 + 3E_2 + 2F_1 + F_2, \\
&N(u)=0.
\end{align*}
If $u\in [1,2]$, then we have
\begin{align*}
&P(u)= \left(2-u, 3(2-u) ,0 ,0 , 2(2-u), 2-u\right),\\
%(2-u)E_1 + 3(2-u)E_2 + 2(2-u)F_1 + (2-u)F_2, \\
& N(u)=\left(0, 3(u-1) ,0 ,0 , 2(u-1), u-1\right).
%3(u-1)E_2+2(u-1)F_1 + (u-1)F_2.
\end{align*}
We note that $-K_{S}-uE_{1}$ is not pseudo effective for $u>2$.
If $u\in [0,1]$, then we have 
\begin{align*}
P(u)^2=5-2u-u^2 ,\quad P(u)E_1 = 1+u.
\end{align*}
If $u\in [1,2]$, then we have 
\begin{align*}
P(u)^2=2(2-u)^2,\quad P(u)E_1 = 4-2u.
\end{align*}
Therefore,  we get 
\begin{align*}
S(E_1)=\frac{1}{5}\int_{0}^{1}(5-2u-u^2) du
+ \frac{1}{5}\int_{1}^{2} 2(2-u)^2 du
=\frac{13}{15}
\end{align*}
by the definition of $S(E_1)$ and 
\begin{align*}
S(W_{\bullet,\bullet}^{E_{1}},p)
%=\frac{2}{5}\int_{0}^{1} \Rm{ord}_{p}(N(u)|_{F}) du
&=\frac{1}{5}\int_{0}^{1} (1+u)^2 du
+ \frac{2}{5}\int_{1}^{2}(4-2u) \Rm{ord}_{p}(N(u)|_{E_{1}})  du
+\frac{1}{5}\int_{1}^{2}(4-2u)^2 du\\
 &= \frac{7}{15}+\frac{4}{15} = \frac{11}{15}
% \begin{cases}
 % \frac{17}{15} & \text{if $p \in E_{1}\cap F_{j} $ for $j=1,2$},\\
 %                   & \\
 %  \frac{7}{15}   & \text{if $p \in E_{1} \setminus \bigcup_{j=1}^{2} F_{j}$}, \\
%  \end{cases}
\end{align*}
by Definition~\ref{S-inv}.
Hence we have 
\begin{align*}
\frac{15}{13}  \geq \delta_{p}(S) \geq  \Rm{min}\left\{ \frac{1}{S(E_1)}, \frac{1}{S(W_{\bullet,\bullet}^{E_1},p)} \right\}
= \frac{15}{13} 
\end{align*}
from Corollary~\ref{AZ}.
Thus, we have $\delta_{p}(S)= 15/13$ in this case.

%%%%%%%%%%%%%%%%%%%%% (2)

\noindent
(2)\;The case $p \in E_2 \setminus F_1$.

\noindent
We calculate $S(E_{2})$ and $S(W_{\bullet,\bullet}^{E_{2}},p)$.
Take $u \in \R_{\geq 0}$. 
Let $P(u)+N(u)$ be the Zariski decomposition of $-K_{S}-uE_{2}$.
%where $P(u)$ is the positive part and $N(u)$ is the negative part.
%
If $u\in [0,1]$, then we have
\begin{align*}
&P(u)= \left(2, 3-u ,0 ,0 , 2-\frac{2}{3}u, 1-\frac{u}{3}\right),\\
%2E_1 + (3-u)E_2 + \left(2-\frac{2}{3}u \right)F_1 + \left(1-\frac{1}{3}u \right)F_2, \\
&N(u)= \left(0, 0 ,0 ,0 , \frac{2}{3}u , \frac{1}{3}u \right).
%\frac{2}{3}u F_1 + \frac{1}{3}u F_2.
\end{align*}
If $u\in [1,3]$, then we have
\begin{align*}
&P(u)= \left(3-u, 3-u ,0 ,0 , 2-\frac{2}{3}u, 1-\frac{u}{3}\right),\\
%(3-u)E_1 + (3-u)E_2 + \left(2-\frac{2}{3}u \right)F_1 + \left(1-\frac{1}{3}u \right)F_2, \\
& N(u)=\left(u-1, 3 ,0 ,0 , \frac{2}{3}u, \frac{u}{3}\right).
%(u-1)E_1 + \frac{2}{3}u F_1 + \frac{1}{3}u F_2.
\end{align*}
We note that $-K_{S}-uE_{2}$ is not pseudo effective for $u>3$.
If $u\in [0,1]$, then we have 
\begin{align*}
P(u)^2=5-2u-\frac{u^2}{3} ,\quad P(u)E_2 =  \frac{3+u}{3}.
\end{align*}
If $u\in [1,3]$, then we have 
\begin{align*}
P(u)^2=\frac{2}{3}(3-u)^2,\quad P(u)E_2 =  2-\frac{2}{3}u.
\end{align*}
Therefore,  we get 
\begin{align*}
S(E_2)=\frac{1}{5}\int_{0}^{1}(5-2u-\frac{1}{3}u^2) du
+ \frac{1}{5}\int_{1}^{3} \frac{2}{3}(3-u)^2 du
=\frac{17}{15}
\end{align*}
by the definition of $S(E_2)$ and 
\begin{align*}
S(W_{\bullet,\bullet}^{E_{2}},p)
&=\frac{2}{5}\int_{0}^{1} \frac{1}{3}(3+u)\Rm{ord}_{p}(N(u)|_{E_{2}}) du 
+\frac{1}{5}\int_{0}^{1} \frac{1}{9}(3+u)^2 du \\
&+ \frac{2}{5}\int_{1}^{3}\frac{2}{3}(3-u)\Rm{ord}_{p}(N(u)|_{E_{2}})  du
+\frac{1}{5}\int_{1}^{3}\frac{4}{9}(3-u)^2 du\\
& =
\begin{cases}
  \frac{13}{15} & \text{if $p \in E_{1}\cap E_2 $},\\
                    & \\
   \frac{23}{45}   & \text{if $p \in E_{1} \setminus (E_2 \cup F_1)$}, \\
  \end{cases}
\end{align*}
by Definition~\ref{S-inv}.
Hence we have 
\begin{align*}
\frac{15}{17}  \geq \delta_{p}(S) \geq  \Rm{min}\left\{ \frac{1}{S(E_2)}, \frac{1}{S(W_{\bullet,\bullet}^{E_2},p)} \right\}
= \frac{15}{17} 
\end{align*}
from Corollary~\ref{AZ}.
Thus, we have $\delta_{p}(S)= 15/17$ in this case.

%%%%%%%%%%%%%%%%%%%%%%%%(3)
\noindent
(3)\;The case $p \in F_1 \setminus F_2 $.

\noindent
We calculate $S(F_{1})$ and $S(W_{\bullet,\bullet}^{F_{1}},p)$.
Take $u \in \R_{\geq 0}$. 
Let $P(u)+N(u)$ be the Zariski decomposition of $-K_{S}-uF_{1}$.
%where $P(u)$ is the positive part and $N(u)$ is the negative part.
%
If $u\in [0,1]$, then we have
\begin{align*}
&P(u)=\left(2, 3 ,0 ,0 , 2-u, 1-\frac{u}{2}\right),\\
%2E_1 + 3E_2 + \left(2-u \right)F_1 + \left(1-\frac{1}{2}u \right)F_2, \\
&N(u)= \left(0, 0 ,0 ,0 , 0, \frac{u}{2}\right).
%\frac{1}{2}u F_2.
\end{align*}
If $u\in [1,2]$, then we have
\begin{align*}
&P(u)= \left(2, 4-u ,0 ,0 , 2-u, 1-\frac{u}{2}\right),\\
%2E_1 + (4-u)E_2 + \left(2-u \right)F_1 + \left(1-\frac{1}{2}u \right)F_2, \\
& N(u)= \left(0, u-1 ,0 ,0 , 0, \frac{u}{2}\right).
%(u-1)E_2  + \frac{1}{2}u F_2.
\end{align*}
We note that $-K_{S}-uF_{1}$ is not pseudo effective for $u>2$.
If $u\in [0,1]$, then we have 
\begin{align*}
P(u)^2=5-\frac{3}{2}u^2 ,\quad P(u)F_1 =  \frac{3u}{2}.
\end{align*}
If $u\in [1,2]$, then we have 
\begin{align*}
P(u)^2= \frac{1}{2}(2-u)(6+u),\quad P(u)F_1 =  \frac{2+u}{2}.
\end{align*}
Therefore,  we get 
\begin{align*}
S(F_1)=\frac{1}{5}\int_{0}^{1}5 - \frac{3}{2}u^2 du
+ \frac{1}{5}\int_{1}^{2} 6 -2u - \frac{u^2}{2} du
=\frac{19}{15}
\end{align*}
by the definition of $S(F_1)$ and 
\begin{align*}
S(W_{\bullet,\bullet}^{F_{1}},p)
&=\frac{2}{5}\int_{0}^{1} \frac{3u}{2}\Rm{ord}_{p}(N(u)|_{F_{1}}) du 
+\frac{1}{5}\int_{0}^{1} \frac{9}{4}u^2 du \\
&+ \frac{2}{5}\int_{1}^{2}\left(1 + \frac{u}{2}\right)\Rm{ord}_{p}(N(u)|_{F_{1}})  du
+\frac{1}{5}\int_{1}^{2}\left(1 + \frac{u}{2}\right)^2 du\\
& =
\begin{cases}
  \frac{17}{15} & \text{if $p \in (F_{1}\cap E_2) $},\\
                    & \\
   \frac{23}{30}   & \text{if $p \in F_{1} \setminus (E_2 \cup F_2) $}, \\
  \end{cases}
\end{align*}
by Definition~\ref{S-inv}.
Hence we have 
\begin{align*}
\frac{15}{19}  \geq \delta_{p}(S) \geq  \Rm{min}\left\{ \frac{1}{S(F_1)}, \frac{1}{S(W_{\bullet,\bullet}^{F_1},p)} \right\}
= \frac{15}{19} 
\end{align*}
from Corollary~\ref{AZ}.
Thus, we have $\delta_{p}(S)= 15/19$ in this case.

%%%%%%%%%%%%%%%%%%%%%%%%(4)
\noindent
(4)\;The case $p \in F_2 $.

\noindent
We calculate $S(F_{2})$ and $S(W_{\bullet,\bullet}^{F_{2}},p)$.
Take $u \in \R_{\geq 0}$. 
Let $P(u)+N(u)$ be the Zariski decomposition of $-K_{S}-uF_{2}$.
%where $P(u)$ is the positive part and $N(u)$ is the negative part.
%
If $u\in [0,1]$, then we have
\begin{align*}
&P(u)= \left(2, 3 ,0 ,0 , 2-\frac{u}{2}, 1-u\right),\\
%2E_1 + 3E_2 + \left(2-\frac{1}{2}u \right)F_1 + \left(1-u \right)F_2, \\
&N(u)=\left(0, 0 ,0 ,0 , \frac{u}{2}, 0\right).
% \frac{1}{2}u F_1.
\end{align*}
If $u\in [1,2]$, then we have
\begin{align*}
&P(u)=\left(2, 3 ,1-u ,1-u , 2-\frac{u}{2}, 1-u\right),\\
%2E_1 + 3E_2 + (1-u)E_3 + (1-u) E_4 + \left(2-\frac{1}{2}u \right)F_1 + \left(1-u \right)F_2 , \\
&N(u)=  \left(0, 0 ,u-1, u-1 , \frac{u}{2}, 0 \right).
%(u-1)E_3 + (u-1) E_4 + \frac{1}{2}u F_1 .
\end{align*}
If $u\in [2,3]$, then we have
\begin{align*}
&P(u)= \left(2, 5-u ,1-u ,1-u , 3-u, 1-u\right),\\
%2E_1 + (5-u)E_2 + (1-u)E_3 + (1-u) E_4 + \left(3- u \right)F_1 + \left(1-u \right)F_2 , \\
&N(u)= \left(0, u-2 ,u-1 ,u-1 , u-1,  0 \right).
%(u-2)E_2 +  (u-1)E_3 + (u-1) E_4 + (u-1) F_1 .
\end{align*}
We note that $-K_{S}-uF_{2}$ is not pseudo effective for $u>3$.
If $u\in [0,1]$, then we have 
\begin{align*}
P(u)^2=\frac{1}{2}(10-3u^2) ,\quad P(u)F_2 =  \frac{3u}{2}.
\end{align*}
If $u\in [1,2]$, then we have 
\begin{align*}
P(u)^2= \frac{1}{2}(u^2 -8u +14),\quad P(u)F_2 = 2-\frac{u}{2}.
\end{align*}
If $u\in [2,3]$, then we have 
\begin{align*}
P(u)^2= (3-u)^2,\quad P(u)F_2 =  3-u.
\end{align*}
Therefore,  we get 
\begin{align*}
S(F_2)=\frac{1}{5}\int_{0}^{1} 5 - \frac{3}{2}u^2 du
+ \frac{1}{5}\int_{1}^{2} 7 -4u + \frac{u^2}{2} du
+ \frac{1}{5}\int_{2}^{3} (3-u)^2 du
=\frac{7}{5}
\end{align*}
by the definition of $S(F_2)$ and 
\begin{align*}
S(W_{\bullet,\bullet}^{F_{2}},p)
&=\frac{2}{5}\int_{0}^{1} \frac{3u}{2}\Rm{ord}_{p}(N(u)|_{F_{2}}) du 
+\frac{1}{5}\int_{0}^{1} \frac{9}{4}u^2 du \\
&+ \frac{2}{5}\int_{1}^{2}\left(2 - \frac{u}{2}\right)\Rm{ord}_{p}(N(u)|_{F_{2}})  du
+\frac{1}{5}\int_{1}^{2}\left(2 - \frac{u}{2}\right)^2 du\\
&+ \frac{2}{5}\int_{2}^{3}\left(3-u\right)\Rm{ord}_{p}(N(u)|_{F_{2}})  du
+\frac{1}{5}\int_{2}^{3}\left(3-u\right)^2 du\\
& =
\begin{cases}
  \frac{19}{15} & \text{if $p \in (F_{1}\cap F_2) $ },\\
                    & \\
   \frac{31}{30}   & \text{if $p \in F_{2} \cap E_i (i =3, 4) $}, \\
    & \\
   \frac{8}{15} &\text{if $p \in F_{2} \setminus (F_1 \cup E_3 \cup E_4) $}.\\
  \end{cases}
\end{align*}
by Definition~\ref{S-inv}.
Hence we have 
\begin{align*}
\frac{5}{7}  \geq \delta_{p}(S) \geq  \Rm{min}\left\{ \frac{1}{S(F_2)}, \frac{1}{S(W_{\bullet,\bullet}^{F_2},p)} \right\}
= \frac{5}{7} 
\end{align*}
from Corollary~\ref{AZ}.
Thus, we have $\delta_{p}(S)= 5/7$ in this case.

%%%%%%%%%%%%%%%%%%%%%%%%(5)
\noindent
(5)\;The case $p \in E_3 \setminus F_2$.

\noindent
We calculate $S(E_{3})$ and $S(W_{\bullet,\bullet}^{E_{3}},p)$.
Take $u \in \R_{\geq 0}$. 
Let $P(u)+N(u)$ be the Zariski decomposition of $-K_{S}-uE_{3}$.
%where $P(u)$ is the positive part and $N(u)$ is the negative part.
%
If $u\in \left[0,\frac{3}{2}\right]$, then we have
\begin{align*}
&P(u)= \left(2, 3 ,-u ,0 , 2-\frac{u}{3}, 1-\frac{2}{3}u\right),\\
%2E_1 + 3E_2 + \left(2-\frac{1}{3}u \right)F_1 + \left(1-\frac{2}{3}u \right)F_2 - u E_3, \\
&N(u)= \left(0, 0 ,0 ,0 , \frac{u}{3}, \frac{2}{3}u\right).
% \frac{u}{3} F_1 + \frac{2u}{3}F_2.
\end{align*}
If $u\in [\frac{3}{2},2]$, then we have
\begin{align*}
&P(u)=\left(2, 3 ,-u , 3-2u, 3-u, 3-2u\right),\\
%2E_1 + 3E_2 + \left(3-u \right)F_1 + \left(3-2u \right)F_2 -u E_{3} + (3-2u)E_4 , \\
&N(u)= \left(0, 0 ,0 , 2u-3, u-1 , 2(u-1) \right).
% (u-1) F_1 +2(u-1) F_2 + (2u-3)E_4.
\end{align*}
%If $u\in [2,3]$, it holds
%\begin{align*}
%&P(u)=2E_1 + (5-u)E_2 + (1-u)E_3 + (1-u) E_4 + \left(3- u \right)F_1 + \left(1-u \right)F_2 , \\
%&N(u)= (u-2)E_2 +  (u-1)E_3 + (u-1) E_4 + (u-1) F_1 .
%\end{align*}
We note that $-K_{S}-uE_{3}$ is not pseudo effective for $u>2$.
If $u\in [0,\frac{3}{2}]$, then we have 
\begin{align*}
P(u)^2=5-2u-\frac{u^2}{3} ,\quad P(u)E_3 =  1+ \frac{u}{3}.
\end{align*}
If $u\in [\frac{3}{2},2]$, then we have 
\begin{align*}
P(u)^2=   8-6u + u^2 ,\quad P(u)E_3 =  3-u.
\end{align*}
Therefore,  we get 
\begin{align*}
S(E_3)=\frac{1}{5}\int_{0}^{\frac{3}{2}}\left( 5 - 2u - \frac{u^2}{3} \right)du
+ \frac{1}{5}\int_{\frac{3}{2}}^{2}\left( 8-6u + u^2 \right)du
%+ \frac{1}{5}\int_{2}^{3} (3-u)^2 du
=\frac{31}{30}
\end{align*}
by the definition of $S(E_3)$ and 
\begin{align*}
S(W_{\bullet,\bullet}^{E_{3}},p)
&=\frac{2}{5}\int_{0}^{\frac{3}{2}}\left(1+ \frac{u}{3}\right)\Rm{ord}_{p}(N(u)|_{E_3}) du 
+\frac{1}{5}\int_{0}^{\frac{3}{2}}\left(1+ \frac{u}{3}\right) ^2 du \\
&+ \frac{2}{5}\int_{\frac{3}{2}}^{2}\left(3-u\right)\Rm{ord}_{p}(N(u)|_{E_3})  du
+\frac{1}{5}\int_{\frac{3}{2}}^{2}(3-u)^2 du\\
%&+ \frac{2}{5}\int_{2}^{3}\left(3-u\right)\Rm{ord}_{p}(N(u)|_{F_{2}})  du
%+\frac{1}{5}\int_{2}^{3}\left(3-u\right)^2 du\\
& = \frac{19}{30}
%\begin{cases}
 % \frac{19}{15} & \text{if $p \in (F_{1}\cap F_2) $ },\\
  %                  & \\
 %  \frac{31}{30}   & \text{if $p \in F_{2} \cap E_i (i =3, 4) $}, \\
 %   & \\
 %  \frac{8}{15} &\text{if $p \in F_{2} \setminus (F_1 \cup E_3 \cup E_4) $}.\\
%  \end{cases}
\end{align*}
by Definition~\ref{S-inv}.
Hence we have 
\begin{align*}
\frac{30}{31}  \geq \delta_{p}(S) \geq  \Rm{min}\left\{ \frac{1}{S(E_3)}, \frac{1}{S(W_{\bullet,\bullet}^{E_3},p)} \right\}
= \frac{30}{31} 
\end{align*}
from Corollary~\ref{AZ}.
Thus, we have $\delta_{p}(S)= 30/31$ in this case.
%We can have $\delta_{p}(S)= 30/31$ for $p \in E_4 \setminus F_2$ 

%%%%%%%%%%%%%%%%%%%%%%%%(6)
\noindent
(6)\;The case $p \in S \setminus \left( \bigcup_{i,j}( E_{i} \cup F_{j})\right)$.

\noindent
Let $L := \rho^{-1}_{\ast}\Ov{\rho(p)q_1}$.
We note that $L \in |\rho^{\ast}H-E_{2} -F_{1}|$ and $L \sim E_1 + E_2$.
Hence we have 
$-K_{S}-uL \sim (2-u)E_{1} + (3-u)E_{2} + 2 F_1 + F_2=\left(2-u, 3-u ,0 , 0, 2, 1\right)$.
We calculate $S(L)$ and $S(W_{\bullet,\bullet}^{L},p)$.
Take $u \in \R_{\geq 0}$. 
Let $P(u)+N(u)$ be the Zariski decomposition of $\sigma^{\ast}(-K_{S})-uL$,
where $P(u)$ is the positive part and $N(u)$ is the negative part.
If $u\in \left[0,\frac{3}{2}\right]$, then we have
\begin{align*}
&P(u)= \left(2-u, 3-u ,0 , 0, 2-\frac{2}{3}u, 1-\frac{u}{3}\right),\\
%(2-u)E_1 + (3-u)E_2 + \left(2-\frac{2}{3}u \right)F_1 + \left(1-\frac{1}{3}u \right)F_2 , \\
&N(u)=\left(0, 0 ,0 , 0, \frac{2}{3}u, \frac{u}{3}\right).
% \frac{2u}{3} F_1 + \frac{u}{3}F_2.
\end{align*}
If $u\in \left[\frac{3}{2},2\right]$, then we have
\begin{align*}
&P(u)=\left(2-u, 3(2-u) ,0 , 0, 2(2-u), 2-u \right),\\
%(2-u)E_1 + 3(2-u)E_2 + 2\left(2-u \right)F_1 + \left(2-u \right)F_2 , \\
&N(u)=\left(0, 2u-3 ,0 , 0, 2(u-1), u-1 \right).
%(2u-3)E_2 + 2(u-1)F_1 + (u-1)F_2.
\end{align*}
We note that $\sigma^{\ast}(-K_{S})-uL$ is not pseudo effective for $u>2$.
If $u\in [0,\frac{3}{2}]$, then we have 
\begin{align*}
P(u)^2=5 -4u + \frac{2u^{2}}{3}  ,\quad P(u)L = 2-\frac{2}{3}u.
\end{align*}
If $u\in [\frac{3}{2},2]$, then we have 
\begin{align*}
P(u)^2=2(2-u)^2 ,\quad P(u)L = 2(2-u).
\end{align*}
Therefore, we get 
\begin{align*}
S(L)=\frac{1}{5}\int_{0}^{\frac{3}{2}} \left(5 -4u + \frac{2u^{2}}{3} \right) du + \frac{1}{5}\int_{\frac{3}{2}}^{2} 2(2-u)^2 du   =\frac{23}{30}
\end{align*}
by the definition of $S(L)$ and 
\begin{align*}
S(W_{\bullet,\bullet}^{L},p)
&=\frac{2}{5}\int_{0}^{\frac{3}{2}}\left(2-\frac{2}{3}u \right) \Rm{ord}_{p}(N(u)|_{L}) du + 
\frac{1}{5}\int_{0}^{\frac{3}{2}}\left(2-\frac{2}{3}u \right)^2 du  \\
&+\frac{2}{5}\int_{\frac{3}{2}}^{2}2(2-u) \Rm{ord}_{p}(N(u)|_{L}) du + 
\frac{1}{5}\int_{\frac{3}{2}}^{2} 4(2-u)^2 du \\
&= \frac{22}{30}
\end{align*}
by Definition~\ref{S-inv}.
Hence we have 
\begin{align*}
\frac{30}{23}  \geq \delta_{p}(S) \geq  \Rm{min}\left\{ \frac{2}{S(L)}, \frac{1}{S(W_{\bullet,\bullet}^{L},p)} \right\}
= \frac{30}{23} 
\end{align*}
from Corollary~\ref{AZ}.
 Thus, we have $\delta_{p}(S)= 30/23$ in this case.
\end{proof}

\begin{prop}
Let $S$ be the anti-canonical degree $5$ weak del Pezzo surface such
that the dual graph of negative curves is
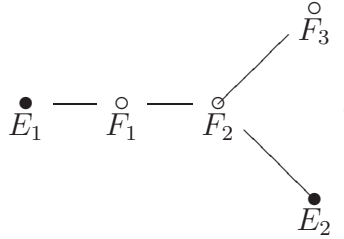
\[
\xygraph{
    {\bullet}* +!U{E_{1}} - [r]
    {\circ}*+!U{F_{1}}  -[r]
 %   {\circ}* +!U{F_{1}} - [r]
   {\circ}* +!D{ }*+!U{F_{2}}
      (- [ru] 
    {\circ}* +!U{F_{3}} 
   % {\bullet}*+!U{E_{4}}  
   %     - [r] \cdots
   %     - [r] \bullet 
   %     - [r] \bullet 
     %  - [rd] \bullet
              ,
                -[rd]
    {\bullet}*+!U{E_{2}} 
   % {\bullet}*+!U{E_{6}}
   %     - [r] \cdots
  %      - [r] \bullet 
  %      - [r] \bullet 
    %     - [ru] {\bullet}*+!U{E_{7}}
         },
        \]
where $E_{i}$ $(i=1,2)$ is a $(-1)$-curve and $F_j$ $(j=1,2,3)$ is a $(-2)$-curve.
Then, for a point $p \in S$, it holds that
\begin{equation}
\nonumber
  \delta_{p}(S)=
  \begin{cases}
    \frac{15}{16} & \text{if $p \in E_1 \setminus F_1$,} \\
    \frac{30}{43}                 & \text{if $p \in F_1 \setminus F_2$,} \\
    \frac{5}{9}      & \text{if $p \in F_2 $,}\\
   \frac{15}{19} & \text{if $p \in F_3 \setminus F_2$,}\\
   \frac{10}{13} &  \text{if $p \in E_2 \setminus F_2$,}\\
  \frac{5}{4} & \text{if $p \in S \setminus \left( \bigcup_{i,j}( E_{i} \cup F_{j})\right)$.}\\
 %  \frac{15}{11} & \text{if $p \in S \setminus \left( \bigcup_{i,j}( E_{i} \cup F_{j})\right)$}
  \end{cases}
\end{equation}
\end{prop}

\begin{proof}
We can assume that we get $S$ from $\Bb{P}^2$ as follows.
\begin{itemize}
\item[(1)]
Take two distinct points $q_1, q_4 \in \Bb{P}^2$.
% and a line $l \in |H|$ passing through $q_1$.
Let $\rho_{1}:S_{1}=\Rm{Bl}_{\{q_{1}, q_{4}\}}\Bb{P}^{2} \to \Bb{P}^2$ be the
composition of blowing-ups at points $q_1$, $q_4$ 
and let $q_2 \in S_1$ be the point at which $(\rho_1)^{-1}_{\ast}(\Ov{q_1 q_4})$ and $\rho_{1}^{-1}(q_1)$ meet.
% and let $l_1$ be the proper transform of $l$ by $\rho_1$.
\item[(2)] Let $\rho_{2}: S_2 \to S_1$ be a blowing-up at $q_2$.
Take a point 
\[
q_3 \in \rho_{2}^{-1}(q_2) \setminus \left((\rho_1 \rho_2)^{-1}_{\ast}(\Ov{q_1 q_4}) \cup
(\rho_2)^{-1}_{\ast}\left(\rho_{1}^{-1}(q_1)\right) \right).
\]
\item[(3)] Let $\rho_{3}: S_3 \to S_2$ be a blowing-up at $q_3$.
Then $S=S_3$. Put $\rho=\rho_1 \rho_2 \rho_3$.
\end{itemize}
Moreover, we have $E_{1}=(\rho_2 \rho_3)^{-1}_{\ast}(\rho_1^{-1}(q_4))$, 
$F_1 = \rho^{-1}_{\ast}(\Ov{q_1 q_4})$, 
%$F_{1}=(\rho_2)_{-1}^{\ast}(\rho_1^{-1}(q_1))$, 
$F_2= (\rho_3)^{-1}_{\ast}(\rho_2^{-1}(q_2))$,
$F_3 =(\rho_2 \rho_3)^{-1}_{\ast}(\rho_1^{-1}(q_1))$, 
$E_2 =\rho_{3}^{-1}(q_3)$.
We denote $D=\sum_{i=1}^{2}a_{i} E_{i} + \sum_{j=1}^{3}b_{j}F_{j} \in \Rm{Div}(S)$ ($a_i , b_j \in \Z$)
by $D=(a_1,a_2, b_1, b_2, b_3)$.
The intersection matrix of $\{E_1, E_2, F_1, F_2 , F_3  \}$ is 
\[A:=
\left(
\begin{array}{cc|ccc}
  -1 & 0 &  1  &  0  &  0     \\
  0 & -1 &  0  &  1  &  0     \\ \hline
  1 & 0 &  -2  &  1 &   0     \\ 
  0 & 1 &  1 &  -2  &   1     \\ 
  0 & 0 &  0  &  1  &  -2    \\
\end{array}
\right).
\]
 We note that 
 $
 -K_{S} \sim 2E_1+ 3E_2 + 3F_1 + 4F_2 + 2F_3 = (2,3,3,4,2).
$
%%%%%%%%%%%%%%%%%%%%% (1)

\noindent
(1)\;The case $p \in E_1 \setminus F_1$.

\noindent
We calculate $S(E_{1})$ and $S(W_{\bullet,\bullet}^{E_{1}},p)$.
Take $u \in \R_{\geq 0}$. 
Let $P(u)+N(u)$ be the Zariski decomposition of $-K_{S}-uE_{1}$.
%where $P(u)$ is the positive part and $N(u)$ is the negative part.
%
If $u\in [0,2]$, then we have
\begin{align*}
&P(u)=\left(2-u, 3 ,3-\frac{3}{4}u , 4-\frac{u}{2}, 2-\frac{u}{4}\right),\\
%(2-u)E_1+ 3E_2 + \left(3-\frac{3}{4}u\right)F_1 + \left(4-\frac{u}{2}\right)F_2 + \left(2-\frac{u}{4}\right)F_3. \\
&N(u)=\left(0, 0 ,\frac{3}{4}u , \frac{u}{2}, \frac{u}{4}\right).
%\frac{3}{4}uF_1 + \frac{u}{2}F_2 +  \frac{u}{4}F_3.
\end{align*}
%If $u\in [1,2]$, it holds
%\begin{align*}
%&P(u)= (2-u)E_1 + 3(2-u)E_2 + 2(2-u)F_1 + (2-u)F_2, \\
%& N(u)=3(u-1)E_2+2(u-1)F_1 + (u-1)F_2.
%\end{align*}
We note that $-K_{S}-uE_{1}$ is not pseudo effective for $u>2$.
If $u\in [0,2]$, then we have
\begin{align*}
P(u)^2=5-2u-\frac{u^2}{4} ,\quad P(u)E_1 = 1+\frac{u}{4}.
\end{align*}
Therefore,  we get 
\begin{align*}
S(E_1)=\frac{1}{5}\int_{0}^{2}\left(5-2u-\frac{u^2}{4}\right) du
%+ \frac{1}{5}\int_{1}^{2} 2(2-u)^2 du
=\frac{16}{15}
\end{align*}
by the definition of $S(E_1)$ and 
\begin{align*}
S(W_{\bullet,\bullet}^{E_{1}},p)
%=\frac{2}{5}\int_{0}^{1} \Rm{ord}_{p}(N(u)|_{F}) du
&=%\frac{1}{5}\int_{0}^{1} (1+u)^2 du
+ \frac{2}{5}\int_{0}^{2}\left(1+\frac{u}{4}\right) \Rm{ord}_{p}(N(u)|_{E_{1}})  du
+\frac{1}{5}\int_{0}^{2}\left(1+\frac{u}{4}\right)^2 du\\
 &=\frac{4}{5}
% \begin{cases}
 % \frac{17}{15} & \text{if $p \in E_{1}\cap F_{j} $ for $j=1,2$},\\
 %                   & \\
 %  \frac{7}{15}   & \text{if $p \in E_{1} \setminus \bigcup_{j=1}^{2} F_{j}$}, \\
%  \end{cases}
\end{align*}
by Definition~\ref{S-inv}.
Hence we have 
\begin{align*}
\frac{15}{16}  \geq \delta_{p}(S) \geq  \Rm{min}\left\{ \frac{1}{S(E_1)}, \frac{1}{S(W_{\bullet,\bullet}^{E_1},p)} \right\}
= \frac{15}{16} 
\end{align*}
from Corollary~\ref{AZ}.
Thus, we have $\delta_{p}(S)= 15/16$ in this case.

%%%%%%%%%%%%%%%%%%%%% (2)

\noindent
(2)\;The case $p \in F_1 \setminus F_2$.

\noindent
We calculate $S(F_{1})$ and $S(W_{\bullet,\bullet}^{F_{1}},p)$.
Take $u \in \R_{\geq 0}$. 
Let $P(u)+N(u)$ be the Zariski decomposition of $-K_{S}-uF_{1}$.
%where $P(u)$ is the positive part and $N(u)$ is the negative part.
%
If $u\in [0,1]$, then we have
\begin{align*}
&P(u)=\left(2, 3 ,3-u , 4-\frac{2}{3}u, 2-\frac{u}{3}\right),\\
%2E_1+ 3E_2 + \left(3-u\right)F_1 + \left(4-\frac{2}{3}u\right)F_2 + \left(2-\frac{u}{3}\right)F_3. \\
&N(u)=\left(0, 0 ,0 , \frac{2}{3}u, \frac{u}{3}\right).
% \frac{2}{3}uF_2 +  \frac{u}{3}F_3.
\end{align*}
If $u\in \left[1,\frac{3}{2}\right]$, then we have
\begin{align*}
&P(u)=\left(3-u, 3 ,3-u , 4-\frac{2}{3}u, 2-\frac{u}{3}\right),\\
%(3-u)E_1+ 3E_2 + \left(3-u\right)F_1 + \left(4-\frac{2}{3}u\right)F_2 + \left(2-\frac{u}{3}\right)F_3. \\
&N(u)=\left(u-1, 0 ,0 , \frac{2}{3}u, \frac{u}{3}\right).
% (u-1)E_1 + \frac{2}{3}uF_2 +  \frac{u}{3}F_3.
\end{align*}
If $u\in \left[\frac{3}{2},3\right]$, then we have
\begin{align*}
&P(u)=\left(3-u, 2(3-u) ,3-u , 2(3-u), 3-u \right),\\
%(3-u)E_1+ 2(3-u)E_2 + \left(3-u\right)F_1 + 2\left(3-u\right)F_2 + \left(3-u\right)F_3. \\
&N(u)= \left(u-1, 2u-3 ,0 , 2(u-1), u-1 \right).
%(u-1)E_1 + (2u-3)E_2+ 2(u-1)F_2 +  (u-1)F_3.
\end{align*}
We note that $-K_{S}-uF_{1}$ is not pseudo effective for $u>3$.
If $u\in [0,1]$, then we have 
\begin{align*}
P(u)^2=5-\frac{4}{3}u^2 ,\quad P(u)F_1 =  \frac{4u}{3}.
\end{align*}
If $u\in [1,\frac{3}{2}]$, then we have 
\begin{align*}
P(u)^2= 6 -2u -\frac{u^2}{3},\quad P(u)F_1 =  1+ \frac{u}{3}.
\end{align*}
If $u\in [\frac{3}{2},3]$, then we have 
\begin{align*}
P(u)^2= (3-u)^2,\quad P(u)F_1 =  3-u.
\end{align*}
Therefore,  we get 
\begin{align*}
S(F_1)=\frac{1}{5}\int_{0}^{1}\left(5-\frac{4u^2}{3}\right) du
+ \frac{1}{5}\int_{1}^{\frac{3}{2}} \left(6 -2u -\frac{u^2}{3}\right) du
+ \frac{1}{5}\int_{\frac{3}{2}}^{3} (3-u)^2 du
=\frac{43}{30}
\end{align*}
by the definition of $S(F_1)$ and 
\begin{align*}
S(W_{\bullet,\bullet}^{F_{1}},p)
&=\frac{2}{5}\int_{0}^{1} \frac{4}{3}u\cdot\Rm{ord}_{p}(N(u)|_{F_{1}}) du 
+\frac{1}{5}\int_{0}^{1} \frac{16}{9}u^2 du \\
&+ \frac{2}{5}\int_{1}^{\frac{3}{2}}\left(1+\frac{u}{3}\right)\Rm{ord}_{p}(N(u)|_{F_{1}})  du
+\frac{1}{5}\int_{1}^{\frac{3}{2}}\left(1+\frac{u}{3}\right)^2 du\\
&+ \frac{2}{5}\int_{\frac{3}{2}}^{3}\left(3-u\right)\Rm{ord}_{p}(N(u)|_{F_{1}})  du
+\frac{1}{5}\int_{\frac{3}{2}}^{3}\left(3-u\right)^2 du\\
& =
\begin{cases}
  \frac{9}{5} & \text{if $p \in (F_{1}\cap F_2)$},\\
                    & \\
%   \frac{31}{30}   & \text{if $p \in F_{2} \cap E_i (i =3, 4) $}, \\
 %   & \\
   \frac{49}{90} &\text{if $p \in F_{1} \setminus (E_1\cup F_2)  $}.\\
  \end{cases}
\end{align*}
by Definition~\ref{S-inv}.
Hence we have 
\begin{align*}
\frac{30}{43}  \geq \delta_{p}(S) \geq  \Rm{min}\left\{ \frac{1}{S(F_1)}, \frac{1}{S(W_{\bullet,\bullet}^{F_1},p)} \right\}
= \frac{30}{43} 
\end{align*}
by Corollary~\ref{AZ}.
Thus, we have $\delta_{p}(S)= 30/43$ in this case.

%%%%%%%%%%%%%%%%%%%%% (3)

\noindent
(3)\;The case $p \in F_2$.

\noindent
We calculate $S(F_{2})$ and $S(W_{\bullet,\bullet}^{F_{2}},p)$.
Take $u \in \R_{\geq 0}$. 
Let $P(u)+N(u)$ be the Zariski decomposition of $-K_{S}-uF_{2}$.
%where $P(u)$ is the positive part and $N(u)$ is the negative part.
%
If $u\in [0,1]$, then we have
\begin{align*}
&P(u)=\left(2,3,3-\frac{u}{2},(4-u),2-\frac{u}{2} \right), \\
&N(u)= \left(0,0,\frac{u}{2},0,\frac{u}{2} \right).
\end{align*}
If $u\in \left[1,2\right]$, then we have
\begin{align*}
&P(u)=\left(2,4-u,3-\frac{u}{2},(4-u),2-\frac{u}{2} \right), \\
&N(u)= \left(0,u-1,\frac{u}{2},0,\frac{u}{2} \right).
\end{align*}
If $u\in \left[2,3\right]$, then we have
\begin{align*}
&P(u)=\left(4-u,4-u,4-u,4-u,2-\frac{u}{2} \right), \\
&N(u)= \left(u-2,u-1,u-1,0,\frac{u}{2} \right).
\end{align*}
We note that $-K_{S}-uF_{2}$ is not pseudo effective for $u>3$.
If $u\in [0,1]$, then we have 
\begin{align*}
P(u)^2=5-u^2 ,\quad P(u)F_2 = u.
\end{align*}
If $u\in [1,2]$, then we have 
\begin{align*}
P(u)^2= 6-2u ,\quad P(u)F_2 = 1.
\end{align*}
If $u\in [2,3]$, then we have 
\begin{align*}
P(u)^2= 8 -4u +\frac{u^2}{2},\quad P(u)F_2 =  2-\frac{u}{2}.
\end{align*}
Therefore,  we get 
\begin{align*}
S(F_2)=\frac{1}{5}\int_{0}^{1}\left(5-u^2 \right) du
+ \frac{1}{5}\int_{1}^{2} \left(6 -2u\right) du
+ \frac{1}{5}\int_{2}^{3} \left(8 -4u +\frac{u^2}{2} \right)du
=\frac{9}{5}
\end{align*}
by the definition of $S(F_2)$ and 
\begin{align*}
S(W_{\bullet,\bullet}^{F_{2}},p)
&=\frac{2}{5}\int_{0}^{1} u\cdot\Rm{ord}_{p}(N(u)|_{F_{2}}) du 
+\frac{1}{5}\int_{0}^{1} u^2 du \\
&+ \frac{2}{5}\int_{1}^{2}1\cdot \Rm{ord}_{p}(N(u)|_{F_{2}})  du
+\frac{1}{5}\int_{1}^{2} 1 du\\
&+ \frac{2}{5}\int_{2}^{3}\left(2-\frac{u}{2}\right)\Rm{ord}_{p}(N(u)|_{F_{2}})  du
+\frac{1}{5}\int_{2}^{3}\left(2-\frac{u}{2}\right)^2 du\\
& =
\begin{cases}
  \frac{13}{10} & \text{if $p \in (F_{2}\cap F_3)$},\\
                    & \\
   \frac{43}{30}   & \text{if $p \in (F_{2} \cap F_1)$}, \\
   & \\
    \frac{19}{15}   & \text{if $p \in (F_{2} \cap E_2)$}, \\
   & \\
   \frac{2}{5} &\text{if $p \in F_{2} \setminus ( E_2 \cup F_1\cup F_3)  $}.\\
  \end{cases}
\end{align*}
by Definition~\ref{S-inv}.
Hence we have 
\begin{align*}
\frac{5}{9}  \geq \delta_{p}(S) \geq  \Rm{min}\left\{ \frac{1}{S(F_2)}, \frac{1}{S(W_{\bullet,\bullet}^{F_2},p)} \right\}
= \frac{5}{9} 
\end{align*}
by Corollary~\ref{AZ}.
Thus, we have $\delta_{p}(S)= 5/9$ in this case.

%%%%%%%%%%%%%%%%%%%%% (4)

\noindent
(4)\;The case $p \in E_2 \setminus F_2$.

\noindent
We calculate $S(E_{2})$ and $S(W_{\bullet,\bullet}^{E_{2}},p)$.
Take $u \in \R_{\geq 0}$. 
Let $P(u)+N(u)$ be the Zariski decomposition of $-K_{S}-uE_{2}$.
%where $P(u)$ is the positive part and $N(u)$ is the negative part.
%
If $u\in [0,2]$, then we have
\begin{align*}
&P(u)=\left(2,3-u,3-\frac{u}{2},(4-u),2-\frac{u}{2} \right), \\
&N(u)= \left(0,0,\frac{u}{2},u,\frac{u}{2} \right).
\end{align*}
If $u\in \left[2,3\right]$, then we have
\begin{align*}
&P(u)=\left(2(3-u),3-u,2(3-u),2(3-u),3-u \right), \\
&N(u)= \left(2u-4,0,2u-3,2(u-1),u-1\right).
\end{align*}
%If $u\in \left[2,3\right]$, it holds
%\begin{align*}
%&P(u)=\left(4-u,4-u,4-u,4-u,2-\frac{u}{2} \right) \\
%&N(u)= \left(u-2,u-1,u-1,0,\frac{u}{2} \right)
%\end{align*}

We note that $-K_{S}-uE_{2}$ is not pseudo effective for $u>3$.
If $u\in [0,2]$, then we have 
\begin{align*}
P(u)^2=5-2u ,\quad P(u)E_2 =  1.
\end{align*}
If $u\in [2,3]$, then we have 
\begin{align*}
P(u)^2=(3-u)^2,\quad P(u)E_2 =  3-u.
\end{align*}
Therefore,  we get 
\begin{align*}
S(E_2)=\frac{1}{5}\int_{0}^{2}\left(5-2u \right) du
+ \frac{1}{5}\int_{2}^{3} \left(3-u\right)^2 du
%+ \frac{1}{5}\int_{2}^{3} \left(8 -4u +\frac{u^2}{2} \right)du
=\frac{19}{15}
\end{align*}
by the definition of $S(E_2)$ and 
\begin{align*}
S(W_{\bullet,\bullet}^{E_{2}},p)
&=\frac{2}{5}\int_{0}^{2} 1\cdot \Rm{ord}_{p}(N(u)|_{E_{2}}) du 
+\frac{1}{5}\int_{0}^{2}  1 du \\
&+ \frac{2}{5}\int_{2}^{3}(3-u) \Rm{ord}_{p}(N(u)|_{E_{2}})  du
+\frac{1}{5}\int_{2}^{3} (3-u)^2 du\\
%&+ \frac{2}{5}\int_{2}^{3}\left(2-\frac{u}{2}\right)\Rm{ord}_{p}(N(u)|_{F_{1}})  du
%+\frac{1}{5}\int_{2}^{3}\left(2-\frac{u}{2}\right)^2 du\\
& = \frac{7}{15}
%\begin{cases}
 % \frac{13}{10} & \text{if $p \in (F_{2}\cap F_3)$},\\
  %                  & \\
  % \frac{43}{30}   & \text{if $p \in (F_{2} \cap F_1)$}, \\
 %  & \\
  %  \frac{19}{15}   & \text{if $p \in (F_{2} \cap E_2)$}, \\
 %  & \\
 %  \frac{2}{5} &\text{if $p \in F_{2} \setminus ( E_2 \cup F_1\cup F_3)  $}.\\
 % \end{cases}
\end{align*}
by Definition~\ref{S-inv}.
Hence we have 
\begin{align*}
\frac{15}{19}  \geq \delta_{p}(S) \geq  \Rm{min}\left\{ \frac{1}{S(E_2)}, \frac{1}{S(W_{\bullet,\bullet}^{E_2},p)} \right\}
= \frac{15}{19} 
\end{align*}
by Corollary~\ref{AZ}.
Thus, we have $\delta_{p}(S)= 15/19$ in this case.

%%%%%%%%%%%%%%%%%%%%% (5)

\noindent
(5)\;The case $p \in F_3 \setminus F_2$.

\noindent
We calculate $S(F_{3})$ and $S(W_{\bullet,\bullet}^{F_{3}},p)$.
Take $u \in \R_{\geq 0}$. 
Let $P(u)+N(u)$ be the Zariski decomposition of $-K_{S}-uF_{3}$.
%where $P(u)$ is the positive part and $N(u)$ is the negative part.
%
If $u\in \left[0,\frac{3}{2}\right]$, then we have
\begin{align*}
&P(u)=\left(2,3,3-\frac{u}{3},(4-\frac{2}{3}u),2-u \right), \\
&N(u)= \left(0,0,\frac{u}{3},\frac{2}{3}u,0 \right).
\end{align*}
If $u\in \left[\frac{3}{2},2\right]$, then we have
\begin{align*}
&P(u)=\left(2,2(3-u),4-u,2(3-u),2-u \right), \\
&N(u)= \left(0,2u-3,u-1,2(u-1),0\right).
\end{align*}
%If $u\in \left[2,3\right]$, it holds
%\begin{align*}
%&P(u)=\left(4-u,4-u,4-u,4-u,2-\frac{u}{2} \right) \\
%&N(u)= \left(u-2,u-1,u-1,0,\frac{u}{2} \right)
%\end{align*}
We note that $-K_{S}-uF_{3}$ is not pseudo effective for $u>2$.
If $u\in [0,\frac{3}{2}]$, then we have 
\begin{align*}
P(u)^2=5-\frac{4}{3}u^2 ,\quad P(u)F_3 = \frac{4u}{3}.
\end{align*}
If $u\in [\frac{3}{2},2]$, then we have 
\begin{align*}
P(u)^2= 4\left(2-u\right) ,\quad P(u)F_3 = 2.
\end{align*}
Therefore, we get 
\begin{align*}
S(F_3)=\frac{1}{5}\int_{0}^{\frac{3}{2}}\left(5-\frac{4}{3}u^2 \right) du
+ \frac{1}{5}\int_{\frac{3}{2}}^{2} 4\left(2-u\right) du
%+ \frac{1}{5}\int_{2}^{3} \left(8 -4u +\frac{u^2}{2} \right)du
=\frac{13}{10}
\end{align*}
by the definition of $S(F_3)$ and 
\begin{align*}
S(W_{\bullet,\bullet}^{F_{3}},p)
&=\frac{2}{5}\int_{0}^{\frac{3}{2}} \frac{4u}{3} \Rm{ord}_{p}(N(u)|_{F_{3}}) du 
+\frac{1}{5}\int_{0}^{\frac{3}{2}}  \frac{16}{9}u^2 du \\
&+ \frac{2}{5}\int_{\frac{3}{2}}^{2}2 \cdot \Rm{ord}_{p}(N(u)|_{F_{3}})  du
+\frac{1}{5}\int_{\frac{3}{2}}^{2} 4 du\\
%&+ \frac{2}{5}\int_{2}^{3}\left(2-\frac{u}{2}\right)\Rm{ord}_{p}(N(u)|_{F_{1}})  du
%+\frac{1}{5}\int_{2}^{3}\left(2-\frac{u}{2}\right)^2 du\\
& = \frac{4}{5}
%\begin{cases}
 % \frac{13}{10} & \text{if $p \in (F_{2}\cap F_3)$},\\
  %                  & \\
  % \frac{43}{30}   & \text{if $p \in (F_{2} \cap F_1)$}, \\
 %  & \\
  %  \frac{19}{15}   & \text{if $p \in (F_{2} \cap E_2)$}, \\
 %  & \\
 %  \frac{2}{5} &\text{if $p \in F_{2} \setminus ( E_2 \cup F_1\cup F_3)  $}.\\
 % \end{cases}
\end{align*}
by Definition~\ref{S-inv}.
Hence we have 
\begin{align*}
\frac{10}{13}  \geq \delta_{p}(S) \geq  \Rm{min}\left\{ \frac{1}{S(F_3)}, \frac{1}{S(W_{\bullet,\bullet}^{F_3},p)} \right\}
= \frac{10}{13} 
\end{align*}
by Corollary~\ref{AZ}.
Thus, we have $\delta_{p}(S)= 10/13$ in this case.

%%%%%%%%%%%%%%%%%%%%%%%%(6)
\noindent
(6)\;The case $p \in S \setminus \left( \bigcup_{i,j}( E_{i} \cup F_{j})\right)$.

\noindent
Let $L:= \rho_{\ast}^{-1} \Ov{\rho(p)q_1}$.
We note that $L \in |\rho^{\ast}H-E_{2} -F_{2} - F_{3}|$ and $L \sim E_1 + E_2 + F_1 + F_2$.
Hence we have 
$-K_{S}-uL \sim (2-u)E_{1} + (3-u)E_{2} + (3-u) F_1 + (4-u) F_2 + 2F_3$.
We calculate $S(L)$ and $S(W_{\bullet,\bullet}^{L},p)$.
Take $u \in \R_{\geq 0}$. 
Let $P(u)+N(u)$ be the Zariski decomposition of $\sigma^{\ast}(-K_{S})-uL$,
where $P(u)$ is the positive part and $N(u)$ is the negative part.
If $u\in \left[0,2\right]$, then we have
\begin{align*}
&P(u)=\left(2-u,3-u,3- \frac{5}{4}u,4 - \frac{3}{2}u, 2 -\frac{3}{4}u \right), \\
&N(u)= \left(0,0,\frac{u}{4},\frac{u}{2},\frac{3}{4}u\right).
\end{align*}
%f $u\in \left[\frac{3}{2},2\right]$, it holds
%\begin{align*}
%&P(u)=\left(2(3-u),3-u,2(3-u),2(3-u),3-u \right) \\
%&N(u)= \left(2u-4,0,2u-3,2(u-1),u-1\right)
%\end{align*}
We note that $\sigma^{\ast}(-K_{S})-uL$ is not pseudo effective for $u>2$.
If $u\in [0,2]$, then we have 
\begin{align*}
P(u)^2=5 -4u + \frac{3u^{2}}{4}  ,\quad P(u)L = 2-\frac{3}{4}u.
\end{align*}
Therefore, we get 
\begin{align*}
S(L)=\frac{1}{5}\int_{0}^{2} \left(5 -4u + \frac{3u^{2}}{4} \right) du %+ \frac{1}{5}\int_{\frac{3}{2}}^{2} 2(2-u)^2 du   
=\frac{4}{5}
\end{align*}
by the definition of $S(L)$ and 
\begin{align*}
S(W_{\bullet,\bullet}^{L},p)
%&=\frac{2}{5}\int_{0}^{\frac{3}{2}}\left(2-\frac{2}{3}u \right) \Rm{ord}_{p}(N(u)|_{F}) du + 
=\frac{1}{5}\int_{0}^{2}\left(2-\frac{3}{4}u \right)^2 du  
%&+\frac{2}{5}\int_{\frac{3}{2}}^{2}2(2-u) \Rm{ord}_{p}(N(u)|_{F}) du + 
%\frac{1}{5}\int_{\frac{3}{2}}^{2} 4(2-u)^2 du \\
= \frac{7}{10}
\end{align*}
by Definition~\ref{S-inv}.
Hence we have 
\begin{align*}
\frac{5}{4}  \geq \delta_{p}(S) \geq  \Rm{min}\left\{ \frac{2}{S(L)}, \frac{1}{S(W_{\bullet,\bullet}^{L},p)} \right\}
= \frac{5}{4} 
\end{align*}
from Corollary~\ref{AZ}.
 Thus, we have $\delta_{p}(S)= 5/4$ in this case.
 \end{proof}

\begin{prop}
Let $S$ be the anti-canonical degree $5$ weak del Pezzo surface such
that the dual graph of negative curves is
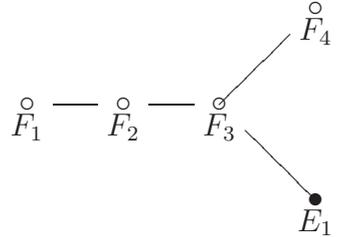
\[
\xygraph{
    {\circ}* +!U{F_{1}} - [r]
    {\circ}*+!U{F_{2}}  -[r]
 %   {\circ}* +!U{F_{1}} - [r]
   {\circ}* +!D{ }*+!U{F_{3}}
      (- [ru] 
    {\circ}* +!U{F_{4}} 
   % {\bullet}*+!U{E_{4}}  
   %     - [r] \cdots
   %     - [r] \bullet 
   %     - [r] \bullet 
     %  - [rd] \bullet
              ,
                -[rd]
    {\bullet}*+!U{E_{1}} 
   % {\bullet}*+!U{E_{6}}
   %     - [r] \cdots
  %      - [r] \bullet 
  %      - [r] \bullet 
    %     - [ru] {\bullet}*+!U{E_{7}}
         },
        \]
where $E_{1}$ is a $(-1)$-curve and $F_j$ $(j=1,2,3,4)$ is a $(-2)$-curve.
Then, for a point $p \in S$, it holds that
\begin{equation}
\nonumber
  \delta_{p}(S)=
  \begin{cases}
    \frac{3}{4} & \text{if $p \in F_1 \setminus F_2$,} \\
    \frac{6}{11}                 & \text{if $p \in F_2 \setminus F_3$,} \\
    \frac{3}{7}      & \text{if $p \in F_3 $,}\\
   \frac{9}{13} & \text{if $p \in F_4 \setminus F_3$,}\\
   \frac{3}{5} &  \text{if $p \in E_1 \setminus F_3$,}\\
  \frac{6}{5} & \text{if $p \in S \setminus \left( E_{1}  \bigcup_{j} F_{j}\right)$.}\\
 %  \frac{15}{11} & \text{if $p \in S \setminus \left( \bigcup_{i,j}( E_{i} \cup F_{j})\right)$}
  \end{cases}
\end{equation}
\end{prop}

\begin{proof}
We can assume that we get $S$ from $\Bb{P}^2$ as follows.
\begin{itemize}
\item[(1)]
Take a point $q_1 \in \Bb{P}^2$ and a line $l$ passing through $q_1$.
Let $\rho_{1}:S_{1}=\Rm{Bl}_{\{q_{1}\}}\Bb{P}^{2} \to \Bb{P}^2$ be the
blowing-up at point $q_1$,
and let $q_2 \in S_1$ be the point at which $(\rho_1)^{-1}_{\ast}l$ and $\rho_{1}^{-1}(q_1)$ meet.
% and let $l_1$ be the proper transform of $l$ by $\rho_1$.
\item[(2)] Let $\rho_{2}: S_2 \to S_1$ be a blowing-up at $q_2$ and 
let $q_3 \in S_2$ be the point at which $(\rho_1 \rho_2)^{-1}_{\ast}l$ and $\rho_{2}^{-1}(q_2)$ meet.
\item[(3)] Let $\rho_{3}: S_3 \to S_2$ be a blowing-up at $q_3$.
Take a point 
\[
q_4 \in \rho_{3}^{-1}(q_3) \setminus \left((\rho_1 \rho_2 \rho _3)^{-1}_{\ast}l \cup
(\rho_3)^{-1}_{\ast}\left(\rho_{2}^{-1}(q_2)\right) \right).
\]
\item[(4)] Let $\rho_{4}: S_4 \to S_3$ be the blowing-up at $q_4$.
Then $S=S_4$. Put $\rho=\rho_1 \rho_2 \rho_3 \rho_4$.
\end{itemize}
Moreover, we have $E_{1}=\rho_{4}^{-1}(q_4)$,
%(\rho_2 \rho_3)_{-1}^{\ast}(\rho_1^{-1}(q_4))$, 
%$F_1 = \rho_{-1}^{\ast}(\Ov{q_1 q_4})$, 
$F_{1}=(\rho_2 \rho_3 \rho_4)^{-1}_{\ast}(\rho_1^{-1}(q_1))$, 
$F_2= (\rho_3 \rho_4)^{-1}_{\ast}(\rho_2^{-1}(q_2))$,
$F_3 =(\rho_4)^{-1}_{\ast}(\rho_3^{-1}(q_3))$, 
$F_4 =\rho^{-1}_{\ast}l$.
We denote $D=a_{1} E_{1} + \sum_{j=1}^{4}b_{j}F_{j} \in \Rm{Div}(S)$ ($a_i , b \in \Z$)
by $D=(a_1, b_1, b_2, b_3, b_4)$.
The intersection matrix of $\{E_1, F_1, F_2 , F_3, F_4  \}$ is 
\[A:=
\left(
\begin{array}{c|cccc}
  -1 & 0 &  0  &  1  &  0     \\ \hline
  0 & -2 &  1  &  0  &  0     \\ 
  0 & 1 &  -2  &  1 &   0     \\ 
  1 & 0 &  1 &  -2  &   1     \\ 
  0 & 0 &  0  &  1  &  -2    \\
\end{array}
\right).
\]
 We note that 
 $
 -K_{S} \sim 5E_1 + 2F_1 + 4F_2 + 6F_3 + 3F_4 = (5,2,4,6,3).
$
%%%%%%%%%%%%%%%%%%%%% (1)

\noindent
(1)\;The case $p \in F_1 \setminus F_2$.

\noindent
We calculate $S(F_{1})$ and $S(W_{\bullet,\bullet}^{F_{1}},p)$.
Take $u \in \R_{\geq 0}$. 
Let $P(u)+N(u)$ be the Zariski decomposition of $-K_{S}-uF_{1}$.
%where $P(u)$ is the positive part and $N(u)$ is the negative part.
%
If $u\in [0,2]$, then we have
\begin{align*}
&P(u)=\left(5, 2-u, 4-\frac{3}{4}u , 6-\frac{u}{2}, 3-\frac{1}{4}u\right),\\
%(2-u)E_1+ 3E_2 + \left(3-\frac{3}{4}u\right)F_1 + \left(4-\frac{u}{2}\right)F_2 + \left(2-\frac{u}{4}\right)F_3. \\
&N(u)=\left(0, 0 ,\frac{3}{4}u , \frac{u}{2}, \frac{1}{4}u\right).
%\frac{3}{4}uF_1 + \frac{u}{2}F_2 +  \frac{u}{4}F_3.
\end{align*}
%If $u\in [1,2]$, it holds
%\begin{align*}
%&P(u)= (2-u)E_1 + 3(2-u)E_2 + 2(2-u)F_1 + (2-u)F_2, \\
%& N(u)=3(u-1)E_2+2(u-1)F_1 + (u-1)F_2.
%\end{align*}
We note that $-K_{S}-uF_{1}$ is not pseudo effective for $u>2$.
If $u\in [0,2]$, then we have 
\begin{align*}
P(u)^2=\frac{5}{4}(4-u^2) ,\quad P(u)F_1 =  \frac{5u}{4}.
\end{align*}
Therefore, we get 
\begin{align*}
S(F_1)=\frac{1}{5}\int_{0}^{2}\frac{5}{4}(4-u^2) du
%+ \frac{1}{5}\int_{1}^{2} 2(2-u)^2 du
=\frac{4}{3}
\end{align*}
by the definition of $S(F_1)$ and 
\begin{align*}
S(W_{\bullet,\bullet}^{F_{1}},p)
%=\frac{2}{5}\int_{0}^{1} \Rm{ord}_{p}(N(u)|_{F}) du
&=%\frac{1}{5}\int_{0}^{1} (1+u)^2 du
+ \frac{2}{5}\int_{0}^{2}\frac{5}{4}u \cdot \Rm{ord}_{p}(N(u)|_{F_{1}})  du
+\frac{1}{5}\int_{0}^{2}\frac{25}{16}u^2 du\\
 &=\frac{11}{6}
% \begin{cases}
 % \frac{17}{15} & \text{if $p \in E_{1}\cap F_{j} $ for $j=1,2$},\\
 %                   & \\
 %  \frac{7}{15}   & \text{if $p \in E_{1} \setminus \bigcup_{j=1}^{2} F_{j}$}, \\
%  \end{cases}
\end{align*}
by Definition~\ref{S-inv}.
Hence we have 
\begin{align*}
\frac{3}{4}  \geq \delta_{p}(S) \geq  \Rm{min}\left\{ \frac{1}{S(F_1)}, \frac{1}{S(W_{\bullet,\bullet}^{F_1},p)} \right\}
= \frac{3}{4} 
\end{align*}
from Corollary~\ref{AZ}.
Thus, we have $\delta_{p}(S)= 3/4$ in this case.

%%%%%%%%%%%%%%%%%%%%% (2)

\noindent
(2)\;The case $p \in F_2 \setminus F_3$.

\noindent
We calculate $S(F_{2})$ and $S(W_{\bullet,\bullet}^{F_{2}},p)$.
Take $u \in \R_{\geq 0}$. 
Let $P(u)+N(u)$ be the Zariski decomposition of $-K_{S}-uF_{2}$.
%where $P(u)$ is the positive part and $N(u)$ is the negative part.
%
If $u\in [0,\frac{3}{2}]$, then we have
\begin{align*}
&P(u)=\left(5, 2-\frac{u}{2} ,4-u , 6-\frac{2}{3}u, 3-\frac{u}{3}\right),\\
%2E_1+ 3E_2 + \left(3-u\right)F_1 + \left(4-\frac{2}{3}u\right)F_2 + \left(2-\frac{u}{3}\right)F_3. \\
&N(u)=\left(0, \frac{u}{2} , 0 , \frac{2}{3}u, \frac{u}{3}\right).
% \frac{2}{3}uF_2 +  \frac{u}{3}F_3.
\end{align*}
If $u\in \left[\frac{3}{2}, 4\right]$, then we have
\begin{align*}
&P(u)=\left(2(4-u), 2-\frac{u}{2} ,4-u , 6-2(u-1), 3-(u-1)\right),\\
%(3-u)E_1+ 3E_2 + \left(3-u\right)F_1 + \left(4-\frac{2}{3}u\right)F_2 + \left(2-\frac{u}{3}\right)F_3. \\
&N(u)=\left(2u-3, \frac{u}{2} ,0 , 2(u-1), u-1\right).
% (u-1)E_1 + \frac{2}{3}uF_2 +  \frac{u}{3}F_3.
\end{align*}
%If $u\in \left[\frac{3}{2},3\right]$, it holds
%\begin{align*}
%&P(u)=\left(3-u, 2(3-u) ,3-u , 2(3-u), 3-u \right)\\
%(3-u)E_1+ 2(3-u)E_2 + \left(3-u\right)F_1 + 2\left(3-u\right)F_2 + \left(3-u\right)F_3. \\
%&N(u)= \left(u-1, 2u-3 ,0 , 2(u-1), u-1 \right)
%(u-1)E_1 + (2u-3)E_2+ 2(u-1)F_2 +  (u-1)F_3.
%\end{align*}
We note that $-K_{S}-uF_{2}$ is not pseudo effective for $u>4$.
If $u\in [0,\frac{3}{2}]$, then we have 
\begin{align*}
P(u)^2=\frac{5}{6}(6-u^2) ,\quad P(u)F_2 = \frac{5}{6}u .
\end{align*}
If $u\in [\frac{3}{2},4]$, then we have 
\begin{align*}
P(u)^2=  \frac{1}{2}\left(4-u\right)^2,\quad P(u)F_2 = \frac{4-u}{2}.
\end{align*}
Therefore,  we get 
\begin{align*}
S(F_2)=\frac{1}{5}\int_{0}^{\frac{3}{2}}\frac{5}{6}(6-u^2) du
+ \frac{1}{5}\int_{\frac{3}{2}}^{4} \frac{1}{2}\left(4-u\right)^2 du
%+ \frac{1}{5}\int_{\frac{3}{2}}^{3} (3-u)^2 du
=\frac{11}{6}
\end{align*}
by the definition of $S(F_2)$ and 
\begin{align*}
S(W_{\bullet,\bullet}^{F_{2}},p)
&=\frac{2}{5}\int_{0}^{\frac{3}{2}} \frac{5}{6}u\cdot\Rm{ord}_{p}(N(u)|_{F_{2}}) du 
+\frac{1}{5}\int_{0}^{\frac{3}{2}} \frac{25}{36}u^2 du \\
&+ \frac{2}{5}\int_{\frac{3}{2}}^{4}\left(\frac{4-u}{2}\right)\Rm{ord}_{p}(N(u)|_{F_{2}})  du
+\frac{1}{5}\int_{\frac{3}{2}}^{4}\left(\frac{4-u}{2}\right)^2 du\\
%&+ \frac{2}{5}\int_{2}^{3}\left(3-u\right)\Rm{ord}_{p}(N(u)|_{F_{1}})  du
%+\frac{1}{5}\int_{2}^{3}\left(3-u\right)^2 du\\
& =
\begin{cases}
  \frac{4}{3} & \text{if $p \in (F_{2}\cap F_1)$},\\
                    & \\
%   \frac{31}{30}   & \text{if $p \in F_{2} \cap E_i (i =3, 4) $}, \\
 %   & \\
   \frac{5}{12} &\text{if $p \in F_{1} \setminus (F_1\cup F_3)  $}.\\
  \end{cases}
\end{align*}
by Definition~\ref{S-inv}.
Hence we have 
\begin{align*}
\frac{6}{11}  \geq \delta_{p}(S) \geq  \Rm{min}\left\{ \frac{1}{S(F_2)}, \frac{1}{S(W_{\bullet,\bullet}^{F_2},p)} \right\}
= \frac{6}{11} 
\end{align*}
by Corollary~\ref{AZ}.
Thus, we have $\delta_{p}(S)= 6/11$ in this case.

%%%%%%%%%%%%%%%%%%%%% (3)

\noindent
(3)\;The case $p \in F_3$.

\noindent
We calculate $S(F_{3})$ and $S(W_{\bullet,\bullet}^{F_{3}},p)$.
Take $u \in \R_{\geq 0}$. 
Let $P(u)+N(u)$ be the Zariski decomposition of $-K_{S}-uF_{3}$.
%where $P(u)$ is the positive part and $N(u)$ is the negative part.
%
If $u\in [0,1]$, then we have
\begin{align*}
&P(u)=\left(5,2-\frac{1}{3}u,4-\frac{2}{3}u,(6-u),3-\frac{u}{2} \right), \\
&N(u)= \left(0,\frac{u}{3},\frac{2}{3}u,0,\frac{u}{2} \right).
\end{align*}
If $u\in \left[1,6\right]$, then we have
\begin{align*}
&P(u)=\left(6-u,2-\frac{1}{3}u,4-\frac{2}{3}u,(6-u),3-\frac{u}{2} \right), \\
&N(u)= \left(u-1,\frac{1}{3}u,\frac{2}{3}u,0,\frac{u}{2} \right).
\end{align*}
%If $u\in \left[2,3\right]$, it holds
%\begin{align*}
%&P(u)=\left(4-u,4-u,4-u,4-u,2-\frac{u}{2} \right) \\
%&N(u)= \left(u-2,u-1,u-1,0,\frac{u}{2} \right)
%\end{align*}

We note that $-K_{S}-uF_{3}$ is not pseudo effective for $u>6$.
If $u\in [0,1]$, then we have 
\begin{align*}
P(u)^2=5-\frac{5}{6}u^2 ,\quad P(u)F_3 = \frac{5}{6}u.
\end{align*}
If $u\in [1,6]$, then we have 
\begin{align*}
P(u)^2= \frac{(6-u)^2}{6} ,\quad P(u)F_3 = \frac{6-u}{6}.
\end{align*}
Therefore,  we get 
\begin{align*}
S(F_3)=\frac{1}{5}\int_{0}^{1}\left(5-\frac{5}{6}u^2 \right) du
+ \frac{1}{5}\int_{1}^{6} \frac{1}{6}\left(6-u \right)^2 du
%+ \frac{1}{5}\int_{2}^{3} \left(8 -4u +\frac{u^2}{2} \right)du
=\frac{7}{3}
\end{align*}
by the definition of $S(F_3)$ and 
\begin{align*}
S(W_{\bullet,\bullet}^{F_{3}},p)
&=\frac{2}{5}\int_{0}^{1} \frac{5}{6}u \cdot\Rm{ord}_{p}(N(u)|_{F_{3}}) du 
+\frac{1}{5}\int_{0}^{1} \frac{25}{36}u^2 du \\
&+ \frac{2}{5}\int_{1}^{6}\frac{6-u}{6}\cdot \Rm{ord}_{p}(N(u)|_{F_{3}})  du
+\frac{1}{5}\int_{1}^{6} \frac{(6-u)^2}{36} du\\
%&+ \frac{2}{5}\int_{2}^{3}\left(2-\frac{u}{2}\right)\Rm{ord}_{p}(N(u)|_{F_{3}})  du
%+\frac{1}{5}\int_{2}^{3}\left(2-\frac{u}{2}\right)^2 du\\
& =
\begin{cases}
  \frac{11}{6} & \text{if $p \in (F_{3}\cap F_2)$},\\
                    & \\
   \frac{13}{9}   & \text{if $p \in (F_{3} \cap F_4)$}, \\
   & \\
    \frac{5}{3}   & \text{if $p \in (F_{3} \cap E_1)$}, \\
   & \\
   \frac{5}{18} &\text{if $p \in F_{3} \setminus ( F_2 \cup F_4\cup E_1)  $}.\\
  \end{cases}
\end{align*}
by Definition~\ref{S-inv}.
Hence we have 
\begin{align*}
\frac{3}{7}  \geq \delta_{p}(S) \geq  \Rm{min}\left\{ \frac{1}{S(F_3)}, \frac{1}{S(W_{\bullet,\bullet}^{F_3},p)} \right\}
= \frac{3}{7} 
\end{align*}
by Corollary~\ref{AZ}.
Thus, we have $\delta_{p}(S)= 3/7$ in this case.

%%%%%%%%%%%%%%%%%%%%% (4)

\noindent
(4)\;The case $p \in F_4 \setminus F_3$.

\noindent
We calculate $S(F_{4})$ and $S(W_{\bullet,\bullet}^{F_{4}},p)$.
Take $u \in \R_{\geq 0}$. 
Let $P(u)+N(u)$ be the Zariski decomposition of $-K_{S}-uF_{4}$.
%where $P(u)$ is the positive part and $N(u)$ is the negative part.
%
If $u\in [0,\frac{4}{3}]$, then we have
\begin{align*}
&P(u)=\left(5,2-\frac{u}{4},4-\frac{u}{2},6-\frac{3}{4}u,3-u \right), \\
&N(u)= \left(0,\frac{u}{4},\frac{u}{2},\frac{3}{4}u,0 \right).
\end{align*}
If $u\in \left[\frac{4}{3},3\right]$, then we have
\begin{align*}
&P(u)=\left(3(3-u),3-u,2(3-u),3(3-u),3-u \right), \\
&N(u)= \left(3u-4,u-1,2(u-1),3(u-1), 0\right).
\end{align*}
%If $u\in \left[2,3\right]$, it holds
%\begin{align*}
%&P(u)=\left(4-u,4-u,4-u,4-u,2-\frac{u}{2} \right) \\
%&N(u)= \left(u-2,u-1,u-1,0,\frac{u}{2} \right)
%\end{align*}
We note that $-K_{S}-uF_{4}$ is not pseudo effective for $u>3$.
If $u\in [0,\frac{4}{3}]$, then we have 
\begin{align*}
P(u)^2=\frac{5}{4}(4-u^2) ,\quad P(u)F_4 = \frac{5}{4}u.
\end{align*}
If $u\in [\frac{4}{3},3]$, then we have 
\begin{align*}
P(u)^2= (3-u)^2 ,\quad P(u)F_4 = 3-u.
\end{align*}
Therefore,  we get 
\begin{align*}
S(F_4)=\frac{1}{5}\int_{0}^{\frac{4}{3}}\frac{5}{4}(4-u^2)  du
+ \frac{1}{5}\int_{\frac{4}{3}}^{3} \left(3-u\right)^2 du
%+ \frac{1}{5}\int_{2}^{3} \left(8 -4u +\frac{u^2}{2} \right)du
=\frac{13}{9}
\end{align*}
by the definition of $S(F_4)$ and 
\begin{align*}
S(W_{\bullet,\bullet}^{F_4},p)
&=\frac{2}{5}\int_{0}^{\frac{4}{3}} \frac{5}{4}u \cdot \Rm{ord}_{p}(N(u)|_{F_4}) du 
+\frac{1}{5}\int_{0}^{\frac{4}{3}}  \frac{25}{16}u^2 du \\
&+ \frac{2}{5}\int_{\frac{4}{3}}^{3}(3-u) \Rm{ord}_{p}(N(u)|_{F_4})  du
+\frac{1}{5}\int_{\frac{4}{3}}^{3} (3-u)^2 du\\
%&+ \frac{2}{5}\int_{2}^{3}\left(2-\frac{u}{2}\right)\Rm{ord}_{p}(N(u)|_{F_{1}})  du
%+\frac{1}{5}\int_{2}^{3}\left(2-\frac{u}{2}\right)^2 du\\
& = \frac{5}{9}
%\begin{cases}
 % \frac{13}{10} & \text{if $p \in (F_{2}\cap F_3)$},\\
  %                  & \\
  % \frac{43}{30}   & \text{if $p \in (F_{2} \cap F_1)$}, \\
 %  & \\
  %  \frac{19}{15}   & \text{if $p \in (F_{2} \cap E_2)$}, \\
 %  & \\
 %  \frac{2}{5} &\text{if $p \in F_{2} \setminus ( E_2 \cup F_1\cup F_3)  $}.\\
 % \end{cases}
\end{align*}
by Definition~\ref{S-inv}.
Hence we have 
\begin{align*}
\frac{9}{13}  \geq \delta_{p}(S) \geq  \Rm{min}\left\{ \frac{1}{S(F_4)}, \frac{1}{S(W_{\bullet,\bullet}^{F_4},p)} \right\}
= \frac{9}{13} 
\end{align*}
by Corollary~\ref{AZ}.
Thus, we have $\delta_{p}(S)= 9/13$ in this case.

%%%%%%%%%%%%%%%%%%%%% (5)

\noindent
(5)\;The case $p \in E_1 \setminus F_3$.

\noindent
We calculate $S(E_1)$ and $S(W_{\bullet,\bullet}^{E_1},p)$.
Take $u \in \R_{\geq 0}$. 
Let $P(u)+N(u)$ be the Zariski decomposition of $-K_{S}-uE_1$.
%where $P(u)$ is the positive part and $N(u)$ is the negative part.
%
If $u\in \left[0,5\right]$, then we have
\begin{align*}
&P(u)=\left(5-u,2-\frac{2}{5}u,4-\frac{4}{5}u,6-\frac{6}{5}u,3-\frac{3}{5}u \right), \\
&N(u)= \left(0,\frac{2}{5}u,\frac{4}{5}u, \frac{6}{5}u, \frac{3}{5}u \right).
\end{align*}
%If $u\in \left[\frac{3}{2},2\right]$, it holds
%\begin{align*}
%&P(u)=\left(2,2(3-u),4-u,2(3-u),2-u \right) \\
%&N(u)= \left(0,2u-3,u-1,2(u-1),0\right)
%\end{align*}
%If $u\in \left[2,3\right]$, it holds
%\begin{align*}
%&P(u)=\left(4-u,4-u,4-u,4-u,2-\frac{u}{2} \right) \\
%&N(u)= \left(u-2,u-1,u-1,0,\frac{u}{2} \right)
%\end{align*}
We note that $-K_{S}-uE_{1}$ is not pseudo effective for $u>5$.
If $u\in [0,5]$, then we have 
\begin{align*}
P(u)^2=\frac{(5-u)^2}{5} ,\quad P(u)E_1 = \frac{5-u}{5} .
\end{align*}
Therefore, we get 
\begin{align*}
S(E_1)=\frac{1}{5}\int_{0}^{5}\frac{1}{5}\left(5-u \right)^2 du
%+ \frac{1}{5}\int_{\frac{3}{2}}^{2} 4\left(2-u\right)^2 du
%+ \frac{1}{5}\int_{2}^{3} \left(8 -4u +\frac{u^2}{2} \right)du
=\frac{5}{3}
\end{align*}
by the definition of $S(E_1)$ and 
\begin{align*}
S(W_{\bullet,\bullet}^{E_1},p)
&=\frac{2}{5}\int_{0}^{5} \frac{5-u}{5} \Rm{ord}_{p}(N(u)|_{E_1}) du 
+\frac{1}{5}\int_{0}^{5}  \frac{(5-u)^2}{25}u^2 du \\
%&+ \frac{2}{5}\int_{\frac{3}{2}}^{2}2 \cdot \Rm{ord}_{p}(N(u)|_{E_1})  du
%+\frac{1}{5}\int_{\frac{3}{2}}^{2} 4 du\\
%&+ \frac{2}{5}\int_{2}^{3}\left(2-\frac{u}{2}\right)\Rm{ord}_{p}(N(u)|_{F_{1}})  du
%+\frac{1}{5}\int_{2}^{3}\left(2-\frac{u}{2}\right)^2 du\\
& = \frac{1}{3}
%\begin{cases}
 % \frac{13}{10} & \text{if $p \in (F_{2}\cap F_3)$},\\
  %                  & \\
  % \frac{43}{30}   & \text{if $p \in (F_{2} \cap F_1)$}, \\
 %  & \\
  %  \frac{19}{15}   & \text{if $p \in (F_{2} \cap E_2)$}, \\
 %  & \\
 %  \frac{2}{5} &\text{if $p \in F_{2} \setminus ( E_2 \cup F_1\cup F_3)  $}.\\
 % \end{cases}
\end{align*}
by Definition~\ref{S-inv}.
Hence we have 
\begin{align*}
\frac{3}{5}  \geq \delta_{p}(S) \geq  \Rm{min}\left\{ \frac{1}{S(E_1)}, \frac{1}{S(W_{\bullet,\bullet}^{E_1},p)} \right\}
= \frac{3}{5} 
\end{align*}
by Corollary~\ref{AZ}.
Thus, we have $\delta_{p}(S)= 3/5$ in this case.

%%%%%%%%%%%%%%%%%%%%%%%%(6)
\noindent
(6)\;The case $p \in S \setminus \left( E_1 \bigcup_{j} F_{j}\right)$.

\noindent
Let $L:= \rho^{-1}_{\ast}\Ov{\rho(p)q_1}$.
We note that $L \sim 2E_1 + F_2 + 2F_3 + F_4$.
Hence we have 
$-K_{S}-uL \sim (5-2u)E_{1} + 2F_1 + (4-u) F_2 + (6-2u)F_3 + (3-u)F_4 = (5-2u, 2, 4-u, 6-2u, 3-u)$.
We calculate $S(L)$ and $S(W_{\bullet,\bullet}^{L},p)$.
Take $u \in \R_{\geq 0}$. 
Let $P(u)+N(u)$ be the Zariski decomposition of $\sigma^{\ast}(-K_{S})-uL$,
where $P(u)$ is the positive part and $N(u)$ is the negative part.
If $u\in \left[0,\frac{5}{2}\right]$, then we have
\begin{align*}
&P(u)=\left(5-2u, 2-\frac{4}{5}u, 4-\frac{8}{5}u, 6-\frac{12}{5}u, 3-\frac{6}{5}u \right), \\
&N(u)= \left(0,\frac{4}{5}u,\frac{3}{5}u,\frac{2}{5}u,\frac{1}{5}u\right).
\end{align*}
%f $u\in \left[\frac{3}{2},2\right]$, it holds
%\begin{align*}
%&P(u)=\left(2(3-u),3-u,2(3-u),2(3-u),3-u \right) \\
%&N(u)= \left(2u-4,0,2u-3,2(u-1),u-1\right)
%\end{align*}
We note that $\sigma^{\ast}(-K_{S})-uL$ is not pseudo effective for $u>5/2$.
If $u\in [0,2]$, then we have 
\begin{align*}
P(u)^2= \frac{1}{5}\left(5-2u\right)^2  ,\quad P(u)L = \frac{5-2u}{5}.
\end{align*}
Therefore, we get 
\begin{align*}
S(L)=\frac{1}{5}\int_{0}^{\frac{5}{2}} \frac{\left(5-2u\right)^2}{5} du %+ \frac{1}{5}\int_{\frac{3}{2}}^{2} 2(2-u)^2 du   
=\frac{5}{6}
\end{align*}
by the definition of $S(L)$ and 
\begin{align*}
S(W_{\bullet,\bullet}^{L},p)
%&=\frac{2}{5}\int_{0}^{\frac{3}{2}}\left(2-\frac{2}{3}u \right) \Rm{ord}_{p}(N(u)|_{F}) du + 
=\frac{1}{5}\int_{0}^{\frac{5}{2}}\frac{(5-2u)^2}{25} du  
%&+\frac{2}{5}\int_{\frac{3}{2}}^{2}2(2-u) \Rm{ord}_{p}(N(u)|_{F}) du + 
%\frac{1}{5}\int_{\frac{3}{2}}^{2} 4(2-u)^2 du \\
= \frac{1}{6}
\end{align*}
by Definition~\ref{S-inv}.
Hence we have 
\begin{align*}
\frac{6}{5}  \geq \delta_{p}(S) \geq  \Rm{min}\left\{ \frac{2}{S(L)}, \frac{1}{S(W_{\bullet,\bullet}^{L},p)} \right\}
= \frac{6}{5} 
\end{align*}
from Corollary~\ref{AZ}.
 Thus, we have $\delta_{p}(S)= 6/5$ in this case.
\end{proof}

\begin{prop}
Let $S$ be the del Pezzo surface with the anti-canonical degree $5$. 
Then, for a point $p \in S$, it holds that
%\begin{table}[htb]
 % \begin{tabular}{|ll|c|c|c|} \hline
 % & & Order & Length of non-trivial orbits &\\ \hline 
 % $\Z_{l}$ &Cyclic group &$ l $  & $1$, $1$ &\\ \hline
 % \end{tabular}
%\end{table}
\begin{equation}
\nonumber
  \delta_{p}(S)=
  \begin{cases}
    \frac{15}{13} & \text{if $p$ lies on a $(-1)$-curve,} \\
    \frac{40}{31}                 & \text{if $p$ does not lies on all $(-1)$-curves.} \\
  \end{cases}
\end{equation}
\end{prop}

\begin{proof}
We recall the construction of $S$.
There exists four points $q_{1}, q_{2}, q_{3},q_{4} \in \Bb{P}^{2}$ 
in general positions such that  $S$ is obtained by $\rho: S = \Rm{Bl}_{\{q_{1}, q_{2}, q_{3}, q_{4}\}}\Bb{P}^{2} \to \Bb{P}^2$.
%We denote by $\rho: S \to \Bb{P}^2 $ this birational morphism.
We have 
$E_1=\rho^{-1}(q_1)$, $E_3=\rho^{-1}(q_2)$, $E_5=\rho^{-1}(q_3)$, $E_7=\rho^{-1}(q_4)$,
$E_2=\rho_{\ast}^{-1}\Ov{q_1 q_2}$, $E_4=\rho_{\ast}^{-1}\Ov{q_2 q_3}$, 
$E_6=\rho_{\ast}^{-1}\Ov{q_3 q_4}$, $E_8=\rho_{\ast}^{-1}\Ov{q_4 q_1}$
$E_9=\rho_{\ast}^{-1}\Ov{q_1 q_3}$ and $E_{10}=\rho_{\ast}^{-1}\Ov{q_4 q_2}$.

%We write a divisor $D=\sum_{i=1}^{10}a_{i} E_{i}  \in \Rm{Div}(S)$ ($a_i , b \in \Z$) as 
%$D=(a_1,a_2,a_3,a_4,a_5,a_6,a_7,a_8, a_9, a_{10})$.
The intersection matrix of $\{E_1, E_2, E_3, E_4, E_5, E_6, E_7, E_8,E_9,E_{10} \}$ is 
\[A:=
\left(
\begin{array}{cccccccccc}
  -1 & 1 & 0  & 0  &  0  &  0  &  0  & 1 & 1   &0\\
  1 & -1 & 1  & 0  &  0  &  1  &  0  & 0 & 0   &  0  \\
  0 & 1 & -1  & 1  &  0  &  0  &  0  & 0 & 0   &1\\
  0 & 0 & 1   & -1 &  1  &  0  &  0  & 1 & 0   & 0 \\
  0 & 1 & 0  & 1  &  -1  &  1  &  0  & 0 & 1   &0 \\
  0 & 1 & 0  & 0  &  1  &  -1  &  1  & 0 & 0   &0\\
  0 & 0 & 0  & 0  &  0  &  1   &  -1 & 1 & 0   &0\\ 
  1 & 0 & 0  & 1  &  0  &  0  &  1   & -1 & 0  &1\\
  1 & 0 & 0  & 0  &  1  &  0  &  0   & 0  & -1 &1\\
  0 & 0 & 1  & 0  &  0  &  0  &  1   & 0  & 1  &-1 \\
\end{array}
\right).
\]

%%%%%%%%%%%%%%%%%%%%% (1)

\noindent
(1)\;The case $p \in E_1$.

\noindent
We denote a divisor $D=a_1E_1 + a_2E_2 + a_8E_8 + a_9E_9  \in \Rm{Div}(S)$ ($a_i , b \in \Z$) by 
$D=(a_1,a_2,a_8, a_9)$.
We give a lower bound of $S(E_{1})$.
Take $u \in \R_{\geq 0}$. 
Let $P(u)+N(u)$ be the Zariski decomposition of $-K_{S}-uE_{1}$.
%where $P(u)$ is the positive part and $N(u)$ is the negative part.
%
If $u\in [0,1]$, it holds
\begin{align*}
&P(u)=\left(2-u, 1,1,1\right),\\
%(2-u)E_1+ 3E_2 + \left(3-\frac{3}{4}u\right)F_1 + \left(4-\frac{u}{2}\right)F_2 + \left(2-\frac{u}{4}\right)F_3. \\
&N(u)=\left(0, 0 ,0,0\right).
%\frac{3}{4}uF_1 + \frac{u}{2}F_2 +  \frac{u}{4}F_3.
\end{align*}
If $u\in [1,2]$, then we have
\begin{align*}
&P(u)=\left(2-u, 2-u,2-u,2-u\right),\\
%(2-u)E_1+ 3E_2 + \left(3-\frac{3}{4}u\right)F_1 + \left(4-\frac{u}{2}\right)F_2 + \left(2-\frac{u}{4}\right)F_3. \\
&N(u)=\left(0, u-1 ,u-1,u-1\right).
%\frac{3}{4}uF_1 + \frac{u}{2}F_2 +  \frac{u}{4}F_3.
\end{align*}
%\Red{We note that $-K_{S}-uE_{1}$ is not pseudo effective for $u>2$.}
%
If $u\in [0,1]$, then we have 
\begin{align*}
P(u)^2=5-2u-u^2 .%\quad P(u)F = 1+u
\end{align*}
If $u\in [1,2]$, then we have 
\begin{align*}
P(u)^2=2(2-u)^2 .%\quad P(u)F = (3-u)
\end{align*}
Therefore,  we get 
\begin{align*}
S(E_1)%=\frac{1}{5}\int_{0}^{2}\left(5-2u-\frac{u^2}{4}\right) du
%+ \frac{1}{5}\int_{1}^{2} 2(2-u)^2 du
\geq \frac{13}{15}
\end{align*}
by the definition of $S(E_1)$.
% and 
%\begin{align*}
%S(W_{\bullet,\bullet}^{E_{1}},p)
%=\frac{2}{5}\int_{0}^{1} \Rm{ord}_{p}(N(u)|_{F}) du
%&=%\frac{1}{5}\int_{0}^{1} (1+u)^2 du
%+ \frac{2}{5}\int_{0}^{2}\left(1+\frac{u}{4}\right) \Rm{ord}_{p}(N(u)|_{E_{1}})  du
%+\frac{1}{5}\int_{0}^{2}\left(1+\frac{u}{4}\right)^2 du\\
% &=\frac{4}{5}
% \begin{cases}
 % \frac{17}{15} & \text{if $p \in E_{1}\cap F_{j} $ for $j=1,2$},\\
 %                   & \\
 %  \frac{7}{15}   & \text{if $p \in E_{1} \setminus \bigcup_{j=1}^{2} F_{j}$}, \\
%  \end{cases}
%\end{align*}
%by \Red{Theorem~$S(W_{\bullet,\bullet}^{F},p)$の明示式}.
Hence we have 
\begin{align*}
\frac{15}{13}  \geq \delta_{p}(S) %\geq  \Rm{min}\left\{ \frac{1}{S(E_1)}, \frac{1}{S(W_{\bullet,\bullet}^{E_1},p)} \right\}
%= \frac{15}{16} 
\end{align*}
%from \Red{Theorem~Abban-Zhang}.
By \cite[\S 2]{FAND}, the (global) delta invariant of del Pezzo surfaces with the anti-canonical degree $5$ is $15/13$.
Hence we have $\delta_{p}(S)= 15/13$ for $p \in E_1$.
Since there exists the automorphisms in $\Rm{Aut}(S)$ which permutes the curves $E_1$, $E_3$, $E_5$ and $E_7$,
we have $\delta_{p}(S)= 15/13$ for $p \in E_i$ for $i=1,3,5,7$.

%%%%%%%%%%%%%%%%%%%%% (2)

\noindent
(2)\;The case $p \in E_2$.

\noindent
We denote a divisor $D= a_2E_2 + a_1E_1 + a_3E_3 + a_6E_6  \in \Rm{Div}(S)$ ($a_i , b \in \Z$) by 
$D=(a_2,a_1,a_3, a_6)$.
We give a lower bound of $S(E_{2})$.
Take $u \in \R_{\geq 0}$. 
Let $P(u)+N(u)$ be the Zariski decomposition of $-K_{S}-uE_{2}$.
%where $P(u)$ is the positive part and $N(u)$ is the negative part.
%
If $u\in [0,1]$, then we have
\begin{align*}
&P(u)=\left(2-u, 1,1,1\right),\\
%(2-u)E_1+ 3E_2 + \left(3-\frac{3}{4}u\right)F_1 + \left(4-\frac{u}{2}\right)F_2 + \left(2-\frac{u}{4}\right)F_3. \\
&N(u)=\left(0, 0 ,0,0\right).
%\frac{3}{4}uF_1 + \frac{u}{2}F_2 +  \frac{u}{4}F_3.
\end{align*}
If $u\in [1,2]$, then we have
\begin{align*}
&P(u)=\left(2-u, 2-u,2-u,2-u\right),\\
%(2-u)E_1+ 3E_2 + \left(3-\frac{3}{4}u\right)F_1 + \left(4-\frac{u}{2}\right)F_2 + \left(2-\frac{u}{4}\right)F_3. \\
&N(u)=\left(0, u-1 ,u-1,u-1\right).
%\frac{3}{4}uF_1 + \frac{u}{2}F_2 +  \frac{u}{4}F_3.
\end{align*}
%\Red{We note that $-K_{S}-uE_{1}$ is not pseudo effective for $u>2$.}
%
If $u\in [0,1]$, then we have 
\begin{align*}
P(u)^2=5-2u-u^2 .%\quad P(u)F = 1+u
\end{align*}
If $u\in [1,2]$, then we have 
\begin{align*}
P(u)^2=2(2-u)^2 .%\quad P(u)F = (3-u)
\end{align*}
Therefore,  we get 
\begin{align*}
S(E_2)%=\frac{1}{5}\int_{0}^{2}\left(5-2u-\frac{u^2}{4}\right) du
%+ \frac{1}{5}\int_{1}^{2} 2(2-u)^2 du
\geq \frac{13}{15}
\end{align*}
by the definition of $S(E_2)$.
% and 
%\begin{align*}
%S(W_{\bullet,\bullet}^{E_{1}},p)
%=\frac{2}{5}\int_{0}^{1} \Rm{ord}_{p}(N(u)|_{F}) du
%&=%\frac{1}{5}\int_{0}^{1} (1+u)^2 du
%+ \frac{2}{5}\int_{0}^{2}\left(1+\frac{u}{4}\right) \Rm{ord}_{p}(N(u)|_{E_{1}})  du
%+\frac{1}{5}\int_{0}^{2}\left(1+\frac{u}{4}\right)^2 du\\
% &=\frac{4}{5}
% \begin{cases}
 % \frac{17}{15} & \text{if $p \in E_{1}\cap F_{j} $ for $j=1,2$},\\
 %                   & \\
 %  \frac{7}{15}   & \text{if $p \in E_{1} \setminus \bigcup_{j=1}^{2} F_{j}$}, \\
%  \end{cases}
%\end{align*}
%by \Red{Theorem~$S(W_{\bullet,\bullet}^{F},p)$の明示式}.
Hence we have 
\begin{align*}
\frac{15}{13}  \geq \delta_{p}(S) %\geq  \Rm{min}\left\{ \frac{1}{S(E_1)}, \frac{1}{S(W_{\bullet,\bullet}^{E_1},p)} \right\}
%= \frac{15}{16} 
\end{align*}
%from \Red{Theorem~Abban-Zhang}.
By \cite[\S 2]{FAND}, the (global) delta invariant of del Pezzo surfaces with the anti-canonical degree $5$ is $15/13$.
Hence we have $\delta_{p}(S)= 15/13$ for $p \in E_2$.
Since there exists the automorphisms in $\Rm{Aut}(S)$ which permutes the curves $E_2$, $E_4$, $E_6$ and $E_8$,
it holds $\delta_{p}(S)= 15/13$ for $p \in E_i$ for $i=2,4,6,8$.

%%%%%%%%%%%%%%%%%%%%% (3)

\noindent
(3)\;The case $p \in E_9$.

\noindent
We denote a divisor $D= a_9E_9 + a_1E_1 + a_5E_5 + a_{10}E_{10}  \in \Rm{Div}(S)$ ($a_i , b \in \Z$) by 
$D=(a_9,a_1,a_5, a_{10})$.
We give a lower bound of $S(E_{9})$.
Take $u \in \R_{\geq 0}$. 
Let $P(u)+N(u)$ be the Zariski decomposition of $-K_{S}-uE_{9}$.
%where $P(u)$ is the positive part and $N(u)$ is the negative part.
%
If $u\in [0,1]$, then we have
\begin{align*}
&P(u)=\left(2-u, 1,1,1\right),\\
%(2-u)E_1+ 3E_2 + \left(3-\frac{3}{4}u\right)F_1 + \left(4-\frac{u}{2}\right)F_2 + \left(2-\frac{u}{4}\right)F_3. \\
&N(u)=\left(0, 0 ,0,0\right).
%\frac{3}{4}uF_1 + \frac{u}{2}F_2 +  \frac{u}{4}F_3.
\end{align*}
If $u\in [1,2]$, then we have
\begin{align*}
&P(u)=\left(2-u, 2-u,2-u,2-u\right),\\
%(2-u)E_1+ 3E_2 + \left(3-\frac{3}{4}u\right)F_1 + \left(4-\frac{u}{2}\right)F_2 + \left(2-\frac{u}{4}\right)F_3. \\
&N(u)=\left(0, u-1 ,u-1,u-1\right).
%\frac{3}{4}uF_1 + \frac{u}{2}F_2 +  \frac{u}{4}F_3.
\end{align*}
%\Red{We note that $-K_{S}-uE_{1}$ is not pseudo effective for $u>2$.}
%
If $u\in [0,1]$, then we have 
\begin{align*}
P(u)^2=5-2u-u^2 .%\quad P(u)F = 1+u
\end{align*}
If $u\in [1,2]$, then we have 
\begin{align*}
P(u)^2=2(2-u)^2 .%\quad P(u)F = (3-u)
\end{align*}
Therefore,  we get 
\begin{align*}
S(E_9)%=\frac{1}{5}\int_{0}^{2}\left(5-2u-\frac{u^2}{4}\right) du
%+ \frac{1}{5}\int_{1}^{2} 2(2-u)^2 du
\geq \frac{13}{15}
\end{align*}
by the definition of $S(E_9)$.
% and 
%\begin{align*}
%S(W_{\bullet,\bullet}^{E_{1}},p)
%=\frac{2}{5}\int_{0}^{1} \Rm{ord}_{p}(N(u)|_{F}) du
%&=%\frac{1}{5}\int_{0}^{1} (1+u)^2 du
%+ \frac{2}{5}\int_{0}^{2}\left(1+\frac{u}{4}\right) \Rm{ord}_{p}(N(u)|_{E_{1}})  du
%+\frac{1}{5}\int_{0}^{2}\left(1+\frac{u}{4}\right)^2 du\\
% &=\frac{4}{5}
% \begin{cases}
 % \frac{17}{15} & \text{if $p \in E_{1}\cap F_{j} $ for $j=1,2$},\\
 %                   & \\
 %  \frac{7}{15}   & \text{if $p \in E_{1} \setminus \bigcup_{j=1}^{2} F_{j}$}, \\
%  \end{cases}
%\end{align*}
%by \Red{Theorem~$S(W_{\bullet,\bullet}^{F},p)$の明示式}.
Hence we have 
\begin{align*}
\frac{15}{13}  \geq \delta_{p}(S) %\geq  \Rm{min}\left\{ \frac{1}{S(E_1)}, \frac{1}{S(W_{\bullet,\bullet}^{E_1},p)} \right\}
%= \frac{15}{16} 
\end{align*}
%from \Red{Theorem~Abban-Zhang}.
By \cite[\S 2]{FAND}, the (global) delta invariant of del Pezzo surfaces with the anti-canonical degree $5$ is $15/13$.
Hence we have $\delta_{p}(S)= 15/13$ for $p \in E_2$.
Since there exists the automorphisms in $\Rm{Aut}(S)$ which permutes the curves $E_9$ and $E_{10}$,
it holds $\delta_{p}(S)= 15/13$ for $p \in E_i$ for $i=9,10$.

%%%%%%%%%%%%%%%%%%%%%%%%(4)
\noindent
(4)\;The case $p \in S \setminus \bigcup_{i} E_{i}$.

\noindent
%To apply \Red{Theorem~Abban-Zhang} for prime divisor $F$,
Consider a blowing up $\sigma:\Ti{S} \to S$ at $p$.
Take a conic $\Ov{C} \subset \Bb{P}^2$ passing through $q_1,q_2,q_3,q_4$ and $\rho(p)$.
Let $C$ be the proper transform of $\Ov{C}$,
$E$ the exceptional curve of $p$, 
$\Ti{E}_{i}$ the proper transform of $E_{i}$ for $i=1, \cdots,10$.
Put $L_i :=(\rho \sigma )^{-1}_{\ast}\Ov{\rho(p)q_i}$.
%We write a divisor $D=\alpha C + \beta_{1} \Ti{E}_{1} +  \beta_{3} \Ti{E}_{3}+  \beta_{5} \Ti{E}_{5} +  \beta_{7} \Ti{E}_{7} + \gamma E
 %\in \Rm{Div}(\Ti{S})$ ($\alpha , \beta_i, \gamma \in \Z$) as 
%$D=(\alpha, \beta_{1}, \beta_{3}, \beta_{5}, \beta_{7}, \gamma)$.
Then we have
\[
\sigma^{\ast}(-K_{S})-uE \sim
\frac{1}{2}(3C+\Ti{E}_1+\Ti{E}_3 + \Ti{E}_5 + \Ti{E}_7 + (3-2u)E).    
\]
We calculate $S(E)$ and $S(W_{\bullet,\bullet}^{E},p)$.
Take $u \in \R_{\geq 0}$. 
Let $\Ti{P}(u)+\Ti{N}(u)$ be the Zariski decomposition of $\sigma^{\ast}(-K_{S})-uE$,
where $\Ti{P}(u)$ is the positive part and $\Ti{N}(u)$ is the negative part.
If $u\in [0,2]$, then we have
\begin{align*}
&\Ti{P}(u)=\frac{1}{2}(3C+\Ti{E}_1+\Ti{E}_3 + \Ti{E}_5 + \Ti{E}_7 + (3-2u)E)  ,\\
%\frac{1}{2}(3,1,1,1,1,3-2u),\\
&\Ti{N}(u)=0.
\end{align*}
If $u\in [2,\frac{5}{2}]$, then we have
\begin{align*}
&\Ti{P}(u)=\frac{1}{2}\big((7-2u)C+\Ti{E}_1+\Ti{E}_3 + \Ti{E}_5 + \Ti{E}_7 + \sum_{i=1}^{4}(2-u)L_i + (3-2u)E\big)  ,\\
%\frac{1}{2}(3,1,1,1,1,3-2u),\\
&\Ti{N}(u)=(u-2)C+ \sum_{i=1}^{4}(u-2)L_i.
\end{align*}
We note that $\sigma^{\ast}(-K_{S})-uE$ is not pseudo effective for $u>5/2$.
If $u \in [ 0,2]$, then we have
\begin{align*}
\Ti{P}(u)^2=5-u^2,\quad \Ti{P}(u)E = u.
\end{align*}
If $u \in [ 2, \frac{5}{2}]$, then we have
\begin{align*}
\Ti{P}(u)^2=21-18u+4u^2,\quad \Ti{P}(u)E = 2(5-2u).
\end{align*}
Therefore, we get 
\begin{align*}
S(E) =\frac{31}{20}
\end{align*}
by the definition of $S(E)$ and 
\begin{align*}
S(W_{\bullet,\bullet}^{E},p)
%\frac{2}{5}\int_{0}^{1}2u \Rm{ord}_{p}(N(u)|_{F}) du
%\frac{1}{5}\int_{0}^{2}u^2 du 
& =
\begin{cases}
  \frac{21}{30} & \text{if $p \in E \cap (C \cup \bigcup_{i=1}^{4}L_i)$},\\
                    & \\
   \frac{2}{3}   & \text{if $p \in E \setminus (C \cup \bigcup_{i=1}^{4}L_i)$}, \\
     \end{cases}
\end{align*}
by Definition~\ref{S-inv}.
Hence we have 
\begin{align*}
\frac{40}{31}  \geq \delta_{p}(S) \geq  \Rm{min}\left\{ \frac{2}{S(E)}, \frac{1}{S(W_{\bullet,\bullet}^{E},p)} \right\}
= \frac{40}{31} 
\end{align*}
from Corollary~\ref{AZ}.
 Thus, we have $\delta_{p}(S)= 40/31$ in this case.

\end{proof}

\section{The case of the anti-canonical degree $6$}

It is known that there exist $6$ types of weak del Pezzo surfaces of the anti-canonical degree $6$
in terms of the configuration of negative curves (\cite{CT}).

\begin{prop}
Let $S$ be the anti-canonical degree $6$ weak del Pezzo surface such
that the dual graph of negative curves is
%The dual graph of negative curves of $S$ is 
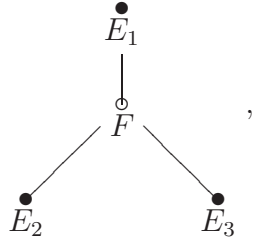
\[
\xygraph{
 %   {\circ}* +!U{F_{1}} - [r]
  %  {\circ}*+!U{F_{2}}  -[r]
 %   {\circ}* +!U{F_{1}} - [r]
   {\circ}* +!D{ }*+!U{F}
      (%- [ru] 
 %   {\circ}* +!U{F_{4}} 
   % {\bullet}*+!U{E_{4}}  
   %     - [r] \cdots
   %     - [r] \bullet 
   %     - [r] \bullet 
     %  - [rd] \bullet
              ,
                -[rd]
    {\bullet}*+!U{E_{3}} 
       ,
                -[ld]
    {\bullet}*+!U{E_{2}} 
     ,
                -[u]
    {\bullet}*+!U{E_{1}}
   % {\bullet}*+!U{E_{6}}
   %     - [r] \cdots
  %      - [r] \bullet 
  %      - [r] \bullet 
    %     - [ru] {\bullet}*+!U{E_{7}}
         },
        \]
where $E_{i}$ $(i=1,2,3)$ is a $(-1)$-curve and $F$ is a $(-2)$-curve.
Then, for a point $p \in S$, it holds that
\begin{equation}
\nonumber
  \delta_{p}(S)=
  \begin{cases}
    \frac{9}{10} & \text{if $p \in E_i \setminus F$,} \\
    \frac{3}{4}                 & \text{if $p \in F $,} \\
  \frac{6}{5} & \text{if $p \in S \setminus \left( \bigcup_{i}  E_{i} \cup  F\right)$.}\\
 %  \frac{15}{11} & \text{if $p \in S \setminus \left( \bigcup_{i,j}( E_{i} \cup F_{j})\right)$}
  \end{cases}
\end{equation}
\end{prop}

\begin{proof}
We can assume that we get $S$ from $\Bb{P}^2$ as follows.
Take three colinear points $q_1, q_2, q_3 \in \Bb{P}^2$
and the line $l$ passing through these points.
Then we have $\rho:S=\Rm{Bl}_{\{q_1, q_2, q_3\}}\Bb{P}^{2} \to \Bb{P}^2$.
Moreover, we have
$E_i := \rho^{-1}(q_i)$ ($i=1,2,3$) and $F=\rho^{-1}_{\ast}l$.
We denote $D=\sum_{i=1}^{3}a_i E_{i} + bF \in \Rm{Div}(S)$ ($a_i , b \in \Z$)
by $D=(a_1, a_2, a_3, b)$.
The intersection matrix of $\{E_1, E_2, E_3, F  \}$ is 
\[A:=
\left(
\begin{array}{ccc|c}
  -1 & 0 &  0  &  1       \\ 
  0 & -1 &  1  &  1       \\ 
  0 & 1 &  -1  &  1      \\ \hline
  1 & 1 &  1 &  -2       \\ 
%  0 & 0 &  0  &  1  &  -2    \\
\end{array}
\right).
\]
 We note that 
 $
 -K_{S} \sim 2E_1 + 2E_2 + 2E_3 + 3F = (2,2,2,3).
$
%%%%%%%%%%%%%%%%%%%%% (1)

\noindent
(1)\;The case $p \in E_1$.

\noindent
We calculate $S(E_{1})$ and $S(W_{\bullet,\bullet}^{E_{1}},p)$.
Take $u \in \R_{\geq 0}$. 
Let $P(u)+N(u)$ be the Zariski decomposition of $-K_{S}-uE_{1}$.
%where $P(u)$ is the positive part and $N(u)$ is the negative part.
%
If $u\in [0,2]$, then we have
\begin{align*}
&P(u)=\left(2-u,2,2, 3-\frac{u}{2}\right),\\
%(2-u)E_1+ 3E_2 + \left(3-\frac{3}{4}u\right)F_1 + \left(4-\frac{u}{2}\right)F_2 + \left(2-\frac{u}{4}\right)F_3. \\
&N(u)=\left(0, 0 ,0, \frac{u}{2}\right).
%\frac{3}{4}uF_1 + \frac{u}{2}F_2 +  \frac{u}{4}F_3.
\end{align*}
%If $u\in [1,2]$, it holds
%\begin{align*}
%&P(u)= (2-u)E_1 + 3(2-u)E_2 + 2(2-u)F_1 + (2-u)F_2, \\
%& N(u)=3(u-1)E_2+2(u-1)F_1 + (u-1)F_2.
%\end{align*}
We note that $-K_{S}-uE_{1}$ is not pseudo effective for $u>2$.
Hence we have 
\begin{align*}
P(u)^2=6-2u-\frac{u^2}{2},\quad P(u)E_1 = \frac{2+u}{2}.
\end{align*}
Therefore,  we get 
\begin{align*}
S(E_1)%=\frac{1}{4}\int_{0}^{2}(2-u)(2+u) du
%+ \frac{1}{5}\int_{1}^{2} 2(2-u)^2 du
=\frac{10}{9}, \quad 
S(W_{\bullet,\bullet}^{E_{1}},p)
%=\frac{2}{5}\int_{0}^{1} \Rm{ord}_{p}(N(u)|_{F}) du
%&=%\frac{1}{5}\int_{0}^{1} (1+u)^2 du
%+ \frac{2}{5}\int_{0}^{2}\frac{5}{4}u \cdot \Rm{ord}_{p}(N(u)|_{E_{1}})  du
%+\frac{1}{5}\int_{0}^{2}\frac{25}{16}u^2 du\\
 &=%\frac{7}{9}
 \begin{cases}
 \frac{7}{9} & \text{if $p \in E_{1} \setminus F $},\\
                    & \\
  \frac{4}{3}   & \text{if $p \in E_{1} \cap F $}. \\
  \end{cases}
\end{align*}
Hence we have 
\begin{align*}
\frac{9}{10}  \geq \delta_{p}(S) \geq  \Rm{min}\left\{ \frac{1}{S(E_1)}, \frac{1}{S(W_{\bullet,\bullet}^{E_1},p)} \right\}
= \frac{9}{10} 
\end{align*}
for a point $p \in E_{1} \setminus F$.
Thus, we have 
\begin{align*}
\delta_{p}(S)
 &
 \begin{cases}
 =\frac{9}{10} & \text{if $p \in E_{1} \setminus F $},\\
                    & \\
 \geq \frac{3}{4}   & \text{if $p \in E_{1} \cap F $}. \\
  \end{cases}
\end{align*}
For $i=2,3$, one can show 
\begin{align*}
\delta_{p}(S)
 &
 \begin{cases}
 =\frac{9}{10} & \text{if $p \in E_{i} \setminus F $},\\
                    & \\
 \geq \frac{3}{4}   & \text{if $p \in E_{i} \cap F $}, \\
  \end{cases}
\end{align*}
by the same calculation.

%%%%%%%%%%%%%%%%%%%%% (2)

\noindent
(2)\;The case $p \in F $.

\noindent
We calculate $S(F)$ and $S(W_{\bullet,\bullet}^{F},p)$.
Take $u \in \R_{\geq 0}$. 
Let $P(u)+N(u)$ be the Zariski decomposition of $-K_{S}-uF$.
%where $P(u)$ is the positive part and $N(u)$ is the negative part.
%
If $u\in [0,1]$, then we have
\begin{align*}
&P(u)=\left(2,2,2, 3-u\right),\\
%(2-u)E_1+ 3E_2 + \left(3-\frac{3}{4}u\right)F_1 + \left(4-\frac{u}{2}\right)F_2 + \left(2-\frac{u}{4}\right)F_3. \\
&N(u)=\left(0, 0 ,0, 0\right).
%\frac{3}{4}uF_1 + \frac{u}{2}F_2 +  \frac{u}{4}F_3.
\end{align*}
If $u\in [1,3]$, then we have
\begin{align*}
&P(u)= \left(3-u,3-u,3-u, 3-u\right), \\
& N(u)=\left(u-1, u-1 ,u-1, 0\right).
\end{align*}
We note that $-K_{S}-uF$ is not pseudo effective for $u>3$.
If $u\in [0,1]$, then we have 
\begin{align*}
P(u)^2=(6-2u^2),\quad P(u)F = 2u.
\end{align*}
If $u\in [1,3]$, then we have 
\begin{align*}
P(u)^2=(3-u)^2,\quad P(u)F = (3-u).
\end{align*}
Therefore,  we get 
\begin{align*}
S(F)%=\frac{1}{4}\int_{0}^{2}(2-u)(2+u) du
%+ \frac{1}{5}\int_{1}^{2} 2(2-u)^2 du
=\frac{4}{3}
\end{align*}
by the definition of $S(F)$.
Hence we get $3/4 \geq \delta_{p}(S)$ for any $p \in F$.
If $p \in F \cap \bigcup_{i=1,2,3} E_i$, then we have $\delta_{p}(S) \geq 3/4$ by (1). 
Hence we get $\delta_{p}(S) = 3/4$ at $p \in F \cap \bigcup_{i=1,2,3} E_i$.
If $p \in F \setminus \bigcup_{i=1,2,3} E_i$, then we have
\begin{align*}
S(W_{\bullet,\bullet}^{F},p)
%=\frac{2}{5}\int_{0}^{1} \Rm{ord}_{p}(N(u)|_{F}) du
%&=%\frac{1}{5}\int_{0}^{1} (1+u)^2 du
%+ \frac{2}{5}\int_{0}^{2}\frac{5}{4}u \cdot \Rm{ord}_{p}(N(u)|_{E_{1}})  du
%+\frac{1}{5}\int_{0}^{2}\frac{25}{16}u^2 du\\
 &=\frac{10}{9}.
% \begin{cases}
 % \frac{17}{15} & \text{if $p \in E_{1}\cap F_{j} $ for $j=1,2$},\\
 %                   & \\
 %  \frac{7}{15}   & \text{if $p \in E_{1} \setminus \bigcup_{j=1}^{2} F_{j}$}, \\
%  \end{cases}
\end{align*}
%by Theorem~\ref{S-inv}.
Hence we have 
\begin{align*}
\frac{3}{4}  \geq \delta_{p}(S) \geq  \Rm{min}\left\{ \frac{1}{S(F)}, \frac{1}{S(W_{\bullet,\bullet}^{F},p)} \right\}
= \frac{3}{4} 
\end{align*}
at a point $p \in F \setminus \bigcup_{i=1,2,3} E_i$.
Thus, we have $\delta_{p}(S)= 3/4$ for any $p \in F$.

%%%%%%%%%%%%%%%%%%%%% (3)

\noindent
(3)\;The case $p \in S \setminus \left( \bigcup_{i}  E_{i} \cup  F\right)$.

Consider a blowing up $\sigma:\Ti{S} \to S$ at $p$.
Let $E$ be the exceptional curve of $p$,
let $\Ti{E}_{i}$ and $\Ti{F}$ be the proper transform of $E_{i}$ and $F$, respectively.
Take three $(-1)$-curves $G_i := (\rho\sigma)_{\ast}^{-1}(\Ov{\rho\sigma(p)q_i})$ 
for $i=1,2,3$.
We note that $\sigma^{\ast}(-K_{S})\sim G_1 + G_2 +G_3 +3E$.
Hence we have 
\[
\sigma^{\ast}(-K_{S})-uE \sim
G_1 + G_2 +G_3 +(3-u)E
\]
We calculate $S(E)$ and $S(W_{\bullet,\bullet}^{E},p)$.
Take $u \in \R_{\geq 0}$. 
Let $\Ti{P}(u)+\Ti{N}(u)$ be the Zariski decomposition of $\sigma^{\ast}(-K_{S})-uE$,
where $\Ti{P}(u)$ is the positive part and $\Ti{N}(u)$ is the negative part.
If $u\in [0,2]$, then we have
\begin{align*}
&\Ti{P}(u)=G_1 + G_2 +G_3 +3E,\\
&\Ti{N}(u)=0.
\end{align*}
If $u\in [2,3]$, then we have
\begin{align*}
&\Ti{P}(u)=(3-u)\left(G_1 + G_2 +G_3 +E\right),\\
&\Ti{N}(u)=(u-2)\left(G_1 + G_2 +G_3 \right).
\end{align*}
We note that $\sigma^{\ast}(-K_{S})-uE$ is not pseudo effective for $u>3$.
If $u\in [0,2]$,then we have 
\begin{align*}
\Ti{P}(u)^2=(6-u^2),\quad \Ti{P}(u)E = u.
\end{align*}
If $u\in [2,3]$, then we have 
\begin{align*}
\Ti{P}(u)^2=2(3-u)^2, \quad \Ti{P}(u)E = 2(3-u).
\end{align*}

Therefore, we get 
\begin{align*}
S(E)%=\frac{1}{5}\int_{0}^{2} 5 -u^2 du 
=\frac{5}{3}, \quad
S(W_{\bullet,\bullet}^{E},p)
%=\frac{2}{5}\int_{0}^{1}2u \Rm{ord}_{p}(N(u)|_{F}) du
%\frac{1}{5}\int_{0}^{2}u^2 du 
= \frac{2}{3}.
\end{align*}
Hence we have 
\begin{align*}
\frac{6}{5}  \geq \delta_{p}(S) \geq  \Rm{min}\left\{ \frac{2}{S(E)}, \frac{1}{S(W_{\bullet,\bullet}^{E},p)} \right\}
= \frac{6}{5} 
\end{align*}
from Corollary~\ref{AZ}.
 Thus, we have $\delta_{p}(S)= 6/5$ in this case.
\end{proof}

\begin{prop}
Let $S$ be the anti-canonical degree $6$ weak del Pezzo surface such
that the dual graph of negative curves is
\[
\xygraph{
   {\bullet}* +!D{ }*+!U{E_1}
                    -[r]
    {\bullet}*+!U{E_{2}}  - [r]
    {\circ}*+!U{F}
   %     - [r] \cdots
%       - [r] {\circ}* +!U{F_{2}}
       - [r] {\bullet}* +!U{E_{3}}
         - [r] {\bullet}*+!U{E_{4}}
              },
        \]
where $E_{i}$ $(i=1,2,3,4)$ is a $(-1)$-curve and $F$ is a $(-2)$-curve.
Then, for a point $p \in S$, it holds that
\begin{equation}
\nonumber
  \delta_{p}(S)=
  \begin{cases}
    \frac{9}{10} & \text{if $p \in (E_1 \setminus E_2 )\cup( E_4 \setminus E_3)$,} \\
    \frac{9}{11}                 & \text{if $p \in E_2 \cup E_3 $,} \\
     \frac{9}{11} & \text{if $p \in F \setminus (E_2 \cup E_3)$,} \\
     \frac{9}{8} & \text{if $p \in S \setminus \left( \bigcup_{i}E_{i} \cup F)\right)$.}
  \end{cases}
\end{equation}
\end{prop}

\begin{proof}
We denote $D=\sum_{i=1}^{4} a_i E_{i} + bF \in \Rm{Div}(S)$ ($a_i , b \in \Z$)
by $D=(a_1, a_2, a_3, a_4,  b)$.
The intersection matrix of $\{E_1, E_2, E_3, E_4, F  \}$ is 
\[A:=
\left(
\begin{array}{cccc|c}
  -1 & 1 &  0  &  0   &  0    \\ 
  1 & -1 &  0  &  0   & 1   \\ 
  0 & 0 &  -1  &  1   & 1 \\ 
  0 & 0 &  1 &  -1    & 0 \\ \hline
  0 & 1 &  1  &  0    &  -2    \\
\end{array}
\right).
\]
 We note that 
 $
 -K_{S} \sim 2E_1 + 3E_2 + E_3 + 2F = (2,3,1,0,2).
$
%%%%%%%%%%%%%%%%%%%%% (1)

\noindent
(1)\;The case $p \in E_1 \setminus E_2$.

\noindent
We calculate $S(E_{1})$ and $S(W_{\bullet,\bullet}^{E_{1}},p)$.
Take $u \in \R_{\geq 0}$. 
Let $P(u)+N(u)$ be the Zariski decomposition of $-K_{S}-uE_{1}$.
%where $P(u)$ is the positive part and $N(u)$ is the negative part.
%
If $u\in [0,1]$, then we have
\begin{align*}
&P(u)=\left(2-u,3,1,0,2\right),\\
%(2-u)E_1+ 3E_2 + \left(3-\frac{3}{4}u\right)F_1 + \left(4-\frac{u}{2}\right)F_2 + \left(2-\frac{u}{4}\right)F_3. \\
&N(u)=\left(0, 0 ,0, 0,0\right).
%\frac{3}{4}uF_1 + \frac{u}{2}F_2 +  \frac{u}{4}F_3.
\end{align*}
If $u\in [1,2]$, then we have
\begin{align*}
&P(u)=\left(2-u,5-2u,1,0,3-u\right),\\
%(2-u)E_1+ 3E_2 + \left(3-\frac{3}{4}u\right)F_1 + \left(4-\frac{u}{2}\right)F_2 + \left(2-\frac{u}{4}\right)F_3. \\
&N(u)=\left(0, 2(u-1) , 0 , 0, u-1\right).
%\frac{3}{4}uF_1 + \frac{u}{2}F_2 +  \frac{u}{4}F_3.
\end{align*}
We note that $-K_{S}-uE_{1}$ is not pseudo effective for $u>2$.
If $u\in [0,1]$, then we have 
\begin{align*}
P(u)^2=(6-2u-u^2),\quad P(u)E_1 = 1+u.
\end{align*}
If $u\in [2,3]$, then we have 
\begin{align*}
P(u)^2=5-2u,\quad P(u)E_1 = (3-u).
\end{align*}

Therefore,  we get 
\begin{align*}
S(E_1)%=\frac{1}{4}\int_{0}^{2}(2-u)(2+u) du
%+ \frac{1}{5}\int_{1}^{2} 2(2-u)^2 du
=\frac{10}{9}, \quad 
S(W_{\bullet,\bullet}^{E_{1}},p)
%=\frac{2}{5}\int_{0}^{1} \Rm{ord}_{p}(N(u)|_{F}) du
%&=%\frac{1}{5}\int_{0}^{1} (1+u)^2 du
%+ \frac{2}{5}\int_{0}^{2}\frac{5}{4}u \cdot \Rm{ord}_{p}(N(u)|_{E_{1}})  du
%+\frac{1}{5}\int_{0}^{2}\frac{25}{16}u^2 du\\
 &=\frac{7}{9}.
% \begin{cases}
 % \frac{17}{15} & \text{if $p \in E_{1}\cap F_{j} $ for $j=1,2$},\\
 %                   & \\
 %  \frac{7}{15}   & \text{if $p \in E_{1} \setminus \bigcup_{j=1}^{2} F_{j}$}, \\
%  \end{cases}
\end{align*}
Hence we have 
\begin{align*}
\frac{9}{10}  \geq \delta_{p}(S) \geq  \Rm{min}\left\{ \frac{1}{S(E_1)}, \frac{1}{S(W_{\bullet,\bullet}^{E_1},p)} \right\}
= \frac{9}{10} 
\end{align*}
from Corollary~\ref{AZ}.
Thus, we have $\delta_{p}(S)= 9/10$ in this case.
We can check $\delta_{p}(S)= 9/10$ for $p \in E_4 \setminus E_3$ by the same calculation.

%%%%%%%%%%%%%%%%%%%%% (2)

\noindent
(2)\;The case $p \in E_2 $.

\noindent
We calculate $S(E_{2})$ and $S(W_{\bullet,\bullet}^{E_{2}},p)$.
Take $u \in \R_{\geq 0}$. 
Let $P(u)+N(u)$ be the Zariski decomposition of $-K_{S}-uE_{2}$.
%where $P(u)$ is the positive part and $N(u)$ is the negative part.
%
If $u\in [0,1]$, then we have
\begin{align*}
&P(u)=\left(2,3-u,1,0, 2- \frac{u}{2}\right),\\
%(2-u)E_1+ 3E_2 + \left(3-\frac{3}{4}u\right)F_1 + \left(4-\frac{u}{2}\right)F_2 + \left(2-\frac{u}{4}\right)F_3. \\
&N(u)=\left(0, 0 ,0, 0, \frac{u}{2} \right).
%\frac{3}{4}uF_1 + \frac{u}{2}F_2 +  \frac{u}{4}F_3.
\end{align*}
If $u\in [1,2]$, then we have
\begin{align*}
&P(u)=\left(3-u,3-u,1,0,2-\frac{u}{2}\right),\\
%(2-u)E_1+ 3E_2 + \left(3-\frac{3}{4}u\right)F_1 + \left(4-\frac{u}{2}\right)F_2 + \left(2-\frac{u}{4}\right)F_3. \\
&N(u)=\left(u-1, 0 ,0, 0, \frac{u}{2}\right).
%\frac{3}{4}uF_1 + \frac{u}{2}F_2 +  \frac{u}{4}F_3.
\end{align*}
If $u\in [2,3]$, then we have
\begin{align*}
&P(u)=\left(3-u,3-u,3-u,0, 3-u\right),\\
%(2-u)E_1+ 3E_2 + \left(3-\frac{3}{4}u\right)F_1 + \left(4-\frac{u}{2}\right)F_2 + \left(2-\frac{u}{4}\right)F_3. \\
&N(u)=\left(u-1, 0 , u-2 , 0, u-1 \right).
%\frac{3}{4}uF_1 + \frac{u}{2}F_2 +  \frac{u}{4}F_3.
\end{align*}
We note that $-K_{S}-uE_{2}$ is not pseudo effective for $u>3$.
If $u\in [0,1]$, then we have 
\begin{align*}
P(u)^2=6-2u-\frac{u^2}{2},\quad P(u)E_2 = 1+\frac{u}{2}.
\end{align*}
If $u\in [1,2]$, then we have 
\begin{align*}
P(u)^2=7-4u+\frac{u^2}{2},\quad P(u)E_2 = 2-\frac{u}{2}.
\end{align*}
If $u\in [2,3]$, then we have 
\begin{align*}
P(u)^2=(3-u)^2,\quad P(u)E_2 = (3-u).
\end{align*}

Therefore,  we get 
\begin{align*}
S(E_2)%=\frac{1}{4}\int_{0}^{2}(2-u)(2+u) du
%+ \frac{1}{5}\int_{1}^{2} 2(2-u)^2 du
=\frac{11}{9}, \quad 
S(W_{\bullet,\bullet}^{E_{2}},p)
%=\frac{2}{5}\int_{0}^{1} \Rm{ord}_{p}(N(u)|_{F}) du
%&=%\frac{1}{5}\int_{0}^{1} (1+u)^2 du
%+ \frac{2}{5}\int_{0}^{2}\frac{5}{4}u \cdot \Rm{ord}_{p}(N(u)|_{E_{1}})  du
%+\frac{1}{5}\int_{0}^{2}\frac{25}{16}u^2 du\\
 &=%\frac{7}{12}
 \begin{cases}
  1 & \text{if $p \in E_{2}\cap E_1 $},\\
                    & \\
   \frac{11}{9}   & \text{if $p \in E_{2} \cap F $}, \\
                    & \\
  \frac{7}{12}   & \text{if $p \in E_{2} \setminus (E_1 \cup F ) $}. \\
  \end{cases}
\end{align*}
Hence we have 
\begin{align*}
\frac{9}{11}  \geq \delta_{p}(S) \geq  \Rm{min}\left\{ \frac{1}{S(E_2)}, \frac{1}{S(W_{\bullet,\bullet}^{E_2},p)} \right\}
= \frac{9}{11} 
\end{align*}
from Corollary~\ref{AZ}.
Thus, we have $\delta_{p}(S)= 9/11$ in this case.
We can check $\delta_{p}(S)= 9/11$ for $p \in E_3 $ by the same calculation.

%%%%%%%%%%%%%%%%%%%%% (3)

\noindent
(3)\;The case $p \in F \setminus (E_2 \cup E_3)$.

\noindent
We calculate $S(F)$ and $S(W_{\bullet,\bullet}^{F},p)$.
Take $u \in \R_{\geq 0}$. 
Let $P(u)+N(u)$ be the Zariski decomposition of $-K_{S}-uF$.
%where $P(u)$ is the positive part and $N(u)$ is the negative part.
%
If $u\in [0,1]$, then we have
\begin{align*}
&P(u)=\left(2,3,1,0, 2-u\right),\\
%(2-u)E_1+ 3E_2 + \left(3-\frac{3}{4}u\right)F_1 + \left(4-\frac{u}{2}\right)F_2 + \left(2-\frac{u}{4}\right)F_3. \\
&N(u)=\left(0, 0 ,0, 0, 0 \right).
%\frac{3}{4}uF_1 + \frac{u}{2}F_2 +  \frac{u}{4}F_3.
\end{align*}
If $u\in [1,2]$, then we have
\begin{align*}
&P(u)=\left(2, 4-u, 2-u ,0, 2-u\right),\\
%(2-u)E_1+ 3E_2 + \left(3-\frac{3}{4}u\right)F_1 + \left(4-\frac{u}{2}\right)F_2 + \left(2-\frac{u}{4}\right)F_3. \\
&N(u)=\left(0, u-1 , u-1, 0 , 0\right).
%\frac{3}{4}uF_1 + \frac{u}{2}F_2 +  \frac{u}{4}F_3.
\end{align*}
%If $u\in [2,3]$, it holds
%\begin{align*}
%&P(u)=\left(3-u,3-u,3-u,0, 3-u\right)\\
%(2-u)E_1+ 3E_2 + \left(3-\frac{3}{4}u\right)F_1 + \left(4-\frac{u}{2}\right)F_2 + \left(2-\frac{u}{4}\right)F_3. \\
%&N(u)=\left(u-1, 0 , u-2 , 0, u-1 \right)
%\frac{3}{4}uF_1 + \frac{u}{2}F_2 +  \frac{u}{4}F_3.
%\end{align*}
We note that $-K_{S}-uF$ is not pseudo effective for $u>2$.
If $u\in [0,1]$, then we have 
\begin{align*}
P(u)^2=6-2u^2,\quad P(u)F = 2u.
\end{align*}
If $u\in [1,2]$, then we have 
\begin{align*}
P(u)^2=4(2-u),\quad P(u)F = 2.
\end{align*}
%f $u\in [2,3]$ we have 
%\begin{align*}
%P(u)^2=(3-u)^2,\quad P(u)F = (3-u)
%\end{align*}
Therefore,  we get 
\begin{align*}
S(F)%=\frac{1}{4}\int_{0}^{2}(2-u)(2+u) du
%+ \frac{1}{5}\int_{1}^{2} 2(2-u)^2 du
=\frac{11}{9}, \quad
S(W_{\bullet,\bullet}^{F},p)
%=\frac{2}{5}\int_{0}^{1} \Rm{ord}_{p}(N(u)|_{F}) du
%&=%\frac{1}{5}\int_{0}^{1} (1+u)^2 du
%+ \frac{2}{5}\int_{0}^{2}\frac{5}{4}u \cdot \Rm{ord}_{p}(N(u)|_{E_{1}})  du
%+\frac{1}{5}\int_{0}^{2}\frac{25}{16}u^2 du\\
 &=\frac{8}{9}.
% \begin{cases}
 % \frac{17}{15} & \text{if $p \in E_{1}\cap F_{j} $ for $j=1,2$},\\
 %                   & \\
 %  \frac{7}{15}   & \text{if $p \in E_{1} \setminus \bigcup_{j=1}^{2} F_{j}$}, \\
%  \end{cases}
\end{align*}
Hence we have 
\begin{align*}
\frac{9}{11}  \geq \delta_{p}(S) \geq  \Rm{min}\left\{ \frac{1}{S(F)}, \frac{1}{S(W_{\bullet,\bullet}^{F},p)} \right\}
= \frac{9}{11} 
\end{align*}
from Corollary~\ref{AZ}.
Thus, we have $\delta_{p}(S)= 9/11$ in this case.
%We can check $\delta_{p}(S)= 9/11$ for $p \in E_3 $ similarly.

%%%%%%%%%%%%%%%%%%%%% (4)

\noindent
(4)\;The case $p \in S \setminus \left( \bigcup_{i}E_{i} \cup F)\right)$.

\noindent
Let $L \in |E_1+E_2|$ be a smooth irreducible curve.
We calculate $S(L)$ and $S(W_{\bullet,\bullet}^{L},p)$.
Take $u \in \R_{\geq 0}$. 
Let $P(u)+N(u)$ be the Zariski decomposition of $-K_{S}-uL$.
%where $P(u)$ is the positive part and $N(u)$ is the negative part.
%
If $u\in [0,2]$, then we have
\begin{align*}
&P(u)=\left(2-u,3-u,1,0, 2-\frac{u}{2}\right),\\
%(2-u)E_1+ 3E_2 + \left(3-\frac{3}{4}u\right)F_1 + \left(4-\frac{u}{2}\right)F_2 + \left(2-\frac{u}{4}\right)F_3. \\
&N(u)=\left(0, 0 ,0, 0, \frac{u}{2} \right).
%\frac{3}{4}uF_1 + \frac{u}{2}F_2 +  \frac{u}{4}F_3.
\end{align*}
%If $u\in [1,2]$, it holds
%\begin{align*}
%&P(u)=\left(2, 4-u, 2-u ,0, 2-u\right)\\
%(2-u)E_1+ 3E_2 + \left(3-\frac{3}{4}u\right)F_1 + \left(4-\frac{u}{2}\right)F_2 + \left(2-\frac{u}{4}\right)F_3. \\
%&N(u)=\left(0, u-1 , u-1, 0 , 0\right)
%\frac{3}{4}uF_1 + \frac{u}{2}F_2 +  \frac{u}{4}F_3.
%\end{align*}
%If $u\in [2,3]$, it holds
%\begin{align*}
%&P(u)=\left(3-u,3-u,3-u,0, 3-u\right)\\
%(2-u)E_1+ 3E_2 + \left(3-\frac{3}{4}u\right)F_1 + \left(4-\frac{u}{2}\right)F_2 + \left(2-\frac{u}{4}\right)F_3. \\
%&N(u)=\left(u-1, 0 , u-2 , 0, u-1 \right)
%\frac{3}{4}uF_1 + \frac{u}{2}F_2 +  \frac{u}{4}F_3.
%\end{align*}
We note that $-K_{S}-uL$ is not pseudo effective for $u>2$.
If $u\in [0,2]$, then we have
\begin{align*}
P(u)^2=\frac{(u-2)(u-6)}{2},\quad P(u)L = \frac{4-u}{2}.
\end{align*}
%If $u\in [1,2]$ we have 
%\begin{align*}
%P(u)^2=4(2-u),\quad P(u)F = 2
%\end{align*}
%f $u\in [2,3]$ we have 
%\begin{align*}
%P(u)^2=(3-u)^2,\quad P(u)F = (3-u)
%\end{align*}
Therefore,  we get 
\begin{align*}
S(L)%=\frac{1}{4}\int_{0}^{2}(2-u)(2+u) du
%+ \frac{1}{5}\int_{1}^{2} 2(2-u)^2 du
=\frac{8}{9}, \quad S(W_{\bullet,\bullet}^{L},p)
%=\frac{2}{5}\int_{0}^{1} \Rm{ord}_{p}(N(u)|_{F}) du
%&=%\frac{1}{5}\int_{0}^{1} (1+u)^2 du
%+ \frac{2}{5}\int_{0}^{2}\frac{5}{4}u \cdot \Rm{ord}_{p}(N(u)|_{E_{1}})  du
%+\frac{1}{5}\int_{0}^{2}\frac{25}{16}u^2 du\\
 &=\frac{7}{9}.
% \begin{cases}
 % \frac{17}{15} & \text{if $p \in E_{1}\cap F_{j} $ for $j=1,2$},\\
 %                   & \\
 %  \frac{7}{15}   & \text{if $p \in E_{1} \setminus \bigcup_{j=1}^{2} F_{j}$}, \\
%  \end{cases}
\end{align*}
Hence we have 
\begin{align*}
\frac{9}{8}  \geq \delta_{p}(S) \geq  \Rm{min}\left\{ \frac{1}{S(L)}, \frac{1}{S(W_{\bullet,\bullet}^{L},p)} \right\}
= \frac{9}{8} 
\end{align*}
from Corollary~\ref{AZ}.
Thus, we have $\delta_{p}(S)= 9/8$ in this case.
%We can check $\delta_{p}(S)= 9/11$ for $p \in E_3 $ similarly.
\end{proof}

\begin{prop}
Let $S$ be the anti-canonical degree $6$ weak del Pezzo surface such
that the dual graph of negative curves is
\[
\xygraph{
   {\circ}* +!D{ }*+!U{F_1}
                    -[r]
    {\bullet}*+!U{E_{1}}  - [r]
    {\circ}*+!U{F_{2}}
   %     - [r] \cdots
%       - [r] {\circ}* +!U{F_{2}}
       - [r] {\bullet}* +!U{E_{2}}
       %  - [r] {\bullet}*+!U{E_{4}}
              },
        \]
where $E_{i}$ $(i=1,2,3,4)$ is a $(-1)$-curve and $F_j$ $(j=1,2)$ is a $(-2)$-curve.
Then, for a point $p \in S$, it holds that
\begin{equation}
\nonumber
  \delta_{p}(S)=
  \begin{cases}
    \frac{9}{11} & \text{if $p \in F_1 \setminus E_1$,} \\
    \frac{9}{14}                 & \text{if $p \in E_1 $,} \\
     \frac{3}{4}                 & \text{if $p \in F_2 \setminus E_1 $,} \\
     \frac{9}{10} & \text{if $p \in  E_2 \setminus F_2$,} \\
     \frac{9}{8} & \text{if $p \in S \setminus \left( E_1\cup E_2 \cup F_1 \cup F_2\right)$.}
  \end{cases}
\end{equation}
\end{prop}

\begin{proof}
We denote $D=\sum_{i=1,2} a_{i}E_{i} + \sum_{j=1,2}b_jF_{j} \in \Rm{Div}(S)$ ($a_i , b_j \in \Z$)
by $D=(a_1, a_2, b_1,b_2)$.
The intersection matrix of $\{E_1, E_2, F_1,F_2  \}$ is 
\[A:=
\left(
\begin{array}{cc|cc}
  -1 & 0 &  1  &  1       \\ 
  0 & -1 &  0  &  1      \\ \hline
  1 & 0 &  -2  &  0    \\ 
  1 & 1 &  0 &  -2     \\ 
%  0 & 1 &  1  &  0       \\
\end{array}
\right).
\]
 We note that 
 $
 -K_{S} \sim 4E_1 + 2E_2 + 2F_1 + 3F_2 = (4,2,2,3).
$
%%%%%%%%%%%%%%%%%%%%% (1)

\noindent
(1)\;The case $p \in F_1 \setminus E_1$.

\noindent
We calculate $S(F_{1})$ and $S(W_{\bullet,\bullet}^{F_{1}},p)$.
Take $u \in \R_{\geq 0}$. 
Let $P(u)+N(u)$ be the Zariski decomposition of $-K_{S}-uF_{1}$.
%where $P(u)$ is the positive part and $N(u)$ is the negative part.
%
If $u\in [0,1]$, then we have
\begin{align*}
&P(u)=\left(4,2,2-u,3\right),\\
%(2-u)E_1+ 3E_2 + \left(3-\frac{3}{4}u\right)F_1 + \left(4-\frac{u}{2}\right)F_2 + \left(2-\frac{u}{4}\right)F_3. \\
&N(u)=\left(0, 0 ,0, 0\right).
%\frac{3}{4}uF_1 + \frac{u}{2}F_2 +  \frac{u}{4}F_3.
\end{align*}
If $u\in [1,2]$, then we have
\begin{align*}
&P(u)=\left(2(3-u),2,2-u,4-u\right),\\
%(2-u)E_1+ 3E_2 + \left(3-\frac{3}{4}u\right)F_1 + \left(4-\frac{u}{2}\right)F_2 + \left(2-\frac{u}{4}\right)F_3. \\
&N(u)=\left(2(u-1), 0 , 0 , u-1\right).
%\frac{3}{4}uF_1 + \frac{u}{2}F_2 +  \frac{u}{4}F_3.
\end{align*}
We note that $-K_{S}-uF_{1}$ is not pseudo effective for $u>2$.
If $u\in [0,1]$, then we have 
\begin{align*}
P(u)^2=(6-2u^2),\quad P(u)F_1 = 2u.
\end{align*}
If $u\in [1,2]$, then we have 
\begin{align*}
P(u)^2=4(2-u),\quad P(u)F_1 = 2.
\end{align*}

Therefore,  we get 
\begin{align*}
S(F_1)%=\frac{1}{4}\int_{0}^{2}(2-u)(2+u) du
%+ \frac{1}{5}\int_{1}^{2} 2(2-u)^2 du
=\frac{11}{9}, \quad 
S(W_{\bullet,\bullet}^{F_{1}},p)
%=\frac{2}{5}\int_{0}^{1} \Rm{ord}_{p}(N(u)|_{F}) du
%&=%\frac{1}{5}\int_{0}^{1} (1+u)^2 du
%+ \frac{2}{5}\int_{0}^{2}\frac{5}{4}u \cdot \Rm{ord}_{p}(N(u)|_{E_{1}})  du
%+\frac{1}{5}\int_{0}^{2}\frac{25}{16}u^2 du\\
 &=\frac{8}{9}.
% \begin{cases}
 % \frac{17}{15} & \text{if $p \in E_{1}\cap F_{j} $ for $j=1,2$},\\
 %                   & \\
 %  \frac{7}{15}   & \text{if $p \in E_{1} \setminus \bigcup_{j=1}^{2} F_{j}$}, \\
%  \end{cases}
\end{align*}
Hence we have 
\begin{align*}
\frac{9}{11}  \geq \delta_{p}(S) \geq  \Rm{min}\left\{ \frac{1}{S(F_1)}, \frac{1}{S(W_{\bullet,\bullet}^{F_1},p)} \right\}
= \frac{9}{11} 
\end{align*}
from Corollary~\ref{AZ}.
Thus, we have $\delta_{p}(S)= 9/11$ in this case.

%%%%%%%%%%%%%%%%%%%%% (2)

\noindent
(2)\;The case $p \in E_1$.

\noindent
We calculate $S(E_{1})$ and $S(W_{\bullet,\bullet}^{E_{1}},p)$.
Take $u \in \R_{\geq 0}$. 
Let $P(u)+N(u)$ be the Zariski decomposition of $-K_{S}-uE_{1}$.
%where $P(u)$ is the positive part and $N(u)$ is the negative part.
%
If $u\in [0,2]$, then we have
\begin{align*}
&P(u)=\left(4-u,2,2-\frac{u}{2},3-\frac{u}{2}\right),\\
%(2-u)E_1+ 3E_2 + \left(3-\frac{3}{4}u\right)F_1 + \left(4-\frac{u}{2}\right)F_2 + \left(2-\frac{u}{4}\right)F_3. \\
&N(u)=\left(0, 0 ,\frac{u}{2} , \frac{u}{2} \right).
%\frac{3}{4}uF_1 + \frac{u}{2}F_2 +  \frac{u}{4}F_3.
\end{align*}
If $u\in [2,4]$, then we have
\begin{align*}
&P(u)=\left(4-u ,4-u,2-\frac{u}{2},4-u\right),\\
%(2-u)E_1+ 3E_2 + \left(3-\frac{3}{4}u\right)F_1 + \left(4-\frac{u}{2}\right)F_2 + \left(2-\frac{u}{4}\right)F_3. \\
&N(u)=\left(0, u-2 , \frac{u}{2} , u-1\right).
%\frac{3}{4}uF_1 + \frac{u}{2}F_2 +  \frac{u}{4}F_3.
\end{align*}
We note that $-K_{S}-uE_{1}$ is not pseudo effective for $u>4$.
If $u\in [0,2]$, then we have 
\begin{align*}
P(u)^2=6-2u,\quad P(u)E_1 = 1.
\end{align*}
If $u\in [2,4]$, then we have 
\begin{align*}
P(u)^2=\frac{(4-u)^2}{2},\quad P(u)E_1 = \frac{4-u}{2}.
\end{align*}

Therefore,  we get 
\begin{align*}
S(E_1)%=\frac{1}{4}\int_{0}^{2}(2-u)(2+u) du
%+ \frac{1}{5}\int_{1}^{2} 2(2-u)^2 du
=\frac{14}{9}, \quad 
S(W_{\bullet,\bullet}^{E_{1}},p)
%=\frac{2}{5}\int_{0}^{1} \Rm{ord}_{p}(N(u)|_{F}) du
%&=%\frac{1}{5}\int_{0}^{1} (1+u)^2 du
%+ \frac{2}{5}\int_{0}^{2}\frac{5}{4}u \cdot \Rm{ord}_{p}(N(u)|_{E_{1}})  du
%+\frac{1}{5}\int_{0}^{2}\frac{25}{16}u^2 du\\
 %=\frac{8}{9}
= \begin{cases}
 1 & \text{if $p \in F_1 $},\\
                    & \\
 \frac{10}{9} & \text{if $p \in F_2 $},\\
                    & \\
  \frac{4}{9}   & \text{if $p \in E_{1} \setminus \bigcup_{j=1}^{2} F_{j}$}. \\
  \end{cases}
\end{align*}
Hence we have 
\begin{align*}
\frac{9}{14}  \geq \delta_{p}(S) \geq  \Rm{min}\left\{ \frac{1}{S(E_1)}, \frac{1}{S(W_{\bullet,\bullet}^{E_1},p)} \right\}
= \frac{9}{14} 
\end{align*}
from Corollary~\ref{AZ}.
Thus, we have $\delta_{p}(S)= 9/14$ in this case.

%%%%%%%%%%%%%%%%%%%%% (3)

\noindent
(3)\;The case $p \in F_2 \setminus E_1$.

\noindent
We calculate $S(F_{2})$ and $S(W_{\bullet,\bullet}^{F_{2}},p)$.
Take $u \in \R_{\geq 0}$. 
Let $P(u)+N(u)$ be the Zariski decomposition of $-K_{S}-uF_{2}$.
%where $P(u)$ is the positive part and $N(u)$ is the negative part.
%
If $u\in [0,1]$, then we have
\begin{align*}
&P(u)=\left(4, 2 , 2, 3-u\right),\\
%(2-u)E_1+ 3E_2 + \left(3-\frac{3}{4}u\right)F_1 + \left(4-\frac{u}{2}\right)F_2 + \left(2-\frac{u}{4}\right)F_3. \\
&N(u)=\left(0, 0 ,0 , 0 \right).
%\frac{3}{4}uF_1 + \frac{u}{2}F_2 +  \frac{u}{4}F_3.
\end{align*}
If $u\in [1,3]$, then we have
\begin{align*}
&P(u)=\left(2(3-u) ,3-u, 3-u, 3-u \right),\\
%(2-u)E_1+ 3E_2 + \left(3-\frac{3}{4}u\right)F_1 + \left(4-\frac{u}{2}\right)F_2 + \left(2-\frac{u}{4}\right)F_3. \\
&N(u)=\left(2(u-1), u-1 , u-1 , 0 \right).
%\frac{3}{4}uF_1 + \frac{u}{2}F_2 +  \frac{u}{4}F_3.
\end{align*}
We note that $-K_{S}-uF_{2}$ is not pseudo effective for $u>3$.
If $u\in [0,1]$, then we have 
\begin{align*}
P(u)^2=6-2u^2,\quad P(u)F_2 = 2u.
\end{align*}
If $u\in [1,3]$, then we have 
\begin{align*}
P(u)^2=(3-u)^2,\quad P(u)F_2 = 3-u.
\end{align*}

Therefore,  we get 
\begin{align*}
S(F_2)%=\frac{1}{4}\int_{0}^{2}(2-u)(2+u) du
%+ \frac{1}{5}\int_{1}^{2} 2(2-u)^2 du
=\frac{4}{3}, \quad S(W_{\bullet,\bullet}^{F_{2}},p)
%=\frac{2}{5}\int_{0}^{1} \Rm{ord}_{p}(N(u)|_{F}) du
%&=%\frac{1}{5}\int_{0}^{1} (1+u)^2 du
%+ \frac{2}{5}\int_{0}^{2}\frac{5}{4}u \cdot \Rm{ord}_{p}(N(u)|_{E_{1}})  du
%+\frac{1}{5}\int_{0}^{2}\frac{25}{16}u^2 du\\
 %=\frac{8}{9}
= \begin{cases}
 \frac{20}{9} & \text{if $p \in F_2 \cap E_1 $},\\
                    & \\
 \frac{16}{9} & \text{if $p \in F_2 \cap E_2 $},\\
                    & \\
  \frac{4}{3}   & \text{if $p \in F_{2} \setminus \bigcup_{i=1}^{2} E_{i}$}. \\
  \end{cases}
\end{align*}
%\begin{align*}
%\frac{3}{4}  \geq \delta_{p}(S) \geq  \Rm{min}\left\{ \frac{1}{S(F_1)}, \frac{1}{S(W_{\bullet,\bullet}^{F_1},p)} \right\}
%= \frac{3}{4} 
%\end{align*}
%from \Red{Theorem~\ref{AZ}}.
Thus, we have 
\begin{align*}
\delta_{p}(S)
 \begin{cases}
 = \frac{3}{4} & \text{if $p \in F_{2} \setminus (E_1 \cup E_2)$},\\
                    & \\
% \frac{16}{9} & \text{if $p \in E_2 $},\\
 %                   & \\
 \leq \frac{3}{4}   & \text{if $\{p\} = F_{2}\cap E_{2}$}. \\
  \end{cases}
\end{align*}
%%%%%%%%%%%%%%%%%%%%% (4)

\noindent
(4)\;The case $p \in E_2 $.

\noindent
We calculate $S(E_{2})$ and $S(W_{\bullet,\bullet}^{E_{2}},p)$.
Take $u \in \R_{\geq 0}$. 
Let $P(u)+N(u)$ be the Zariski decomposition of $-K_{S}-uE_{2}$.
%where $P(u)$ is the positive part and $N(u)$ is the negative part.
%
If $u\in [0,2]$, then we have
\begin{align*}
&P(u)=\left(4, 2-u , 2, 3-\frac{u}{2}\right),\\
%(2-u)E_1+ 3E_2 + \left(3-\frac{3}{4}u\right)F_1 + \left(4-\frac{u}{2}\right)F_2 + \left(2-\frac{u}{4}\right)F_3. \\
&N(u)=\left(0, 0 ,0 , \frac{u}{2} \right).
%\frac{3}{4}uF_1 + \frac{u}{2}F_2 +  \frac{u}{4}F_3.
\end{align*}
%If $u\in [1,3]$, it holds
%\begin{align*}
%&P(u)=\left(2(3-u) ,3-u, 3-u, 3-u \right)\\
%(2-u)E_1+ 3E_2 + \left(3-\frac{3}{4}u\right)F_1 + \left(4-\frac{u}{2}\right)F_2 + \left(2-\frac{u}{4}\right)F_3. \\
%&N(u)=\left(2(u-1), u-1 , u-1 , 0 \right)
%\frac{3}{4}uF_1 + \frac{u}{2}F_2 +  \frac{u}{4}F_3.
%\end{align*}
We note that $-K_{S}-uE_{2}$ is not pseudo effective for $u>2$.
If $u\in [0,2]$, then we have
\begin{align*}
P(u)^2=6-2u-\frac{u^2}{2},\quad P(u)E_2 = 1+\frac{u}{2}.
\end{align*}
%If $u\in [1,3]$ we have 
%\begin{align*}
%P(u)^2=(3-u)^2,\quad P(u)F_2 = 3-u
%\end{align*}

Therefore,  we get 
\begin{align*}
S(E_2)%=\frac{1}{4}\int_{0}^{2}(2-u)(2+u) du
%+ \frac{1}{5}\int_{1}^{2} 2(2-u)^2 du
=\frac{10}{9}, \quad S(W_{\bullet,\bullet}^{E_{2}},p)
%=\frac{2}{5}\int_{0}^{1} \Rm{ord}_{p}(N(u)|_{F}) du
%&=%\frac{1}{5}\int_{0}^{1} (1+u)^2 du
%+ \frac{2}{5}\int_{0}^{2}\frac{5}{4}u \cdot \Rm{ord}_{p}(N(u)|_{E_{1}})  du
%+\frac{1}{5}\int_{0}^{2}\frac{25}{16}u^2 du\\
 %=\frac{8}{9}
= \begin{cases}
 \frac{4}{3} & \text{if $p \in E_2 \cap F_2 $},\\
                    & \\
% \frac{16}{9} & \text{if $p \in E_2 $},\\
 %                   & \\
  \frac{7}{9}   & \text{if $p \in E_{2} \setminus F_2$}. \\
  \end{cases}
\end{align*}
%by Theorem~\ref{S-inv}.
Thus, we have 
\begin{align*}
\delta_{p}(S) 
 \begin{cases}
 \geq \frac{3}{4}   & \text{if $\{p\} = F_{2}\cap E_{2}$} ,\\
                    & \\
% \frac{16}{9} & \text{if $p \in E_2 $},\\
 %                   & \\
= \frac{9}{10} & \text{if $p \in E_{2} \setminus F_2$}.  \\
  \end{cases}
\end{align*}
%\begin{align*}
%\frac{3}{4}  \geq \delta_{p}(S) \geq  \Rm{min}\left\{ \frac{1}{S(F_1)}, \frac{1}{S(W_{\bullet,\bullet}^{F_1},p)} \right\}
%= \frac{3}{4} 
%\end{align*}
%from \Red{Theorem~\ref{AZ}}.
By (3), we have $  3/4 \geq \delta_{p}(S)$ for $\{p\} = F_{2}\cap E_{2}$.
Therefore, we get $\delta_{p}(S) = 3/4$ for $\{p\} = F_{2}\cap E_{2}$.

%%%%%%%%%%%%%%%%%%%%% (5)

\noindent
(5)\;The case $p \in S \setminus \left( E_1\cup E_2 \cup F_1 \cup F_2\right)$.

\noindent
Let $L \in |E_1+E_2+F_2|$ be a smooth irreducible curve.
We calculate $S(L)$ and $S(W_{\bullet,\bullet}^{L},p)$.
Take $u \in \R_{\geq 0}$. 
Let $P(u)+N(u)$ be the Zariski decomposition of $-K_{S}-uL$.
%where $P(u)$ is the positive part and $N(u)$ is the negative part.
%
If $u\in [0,2]$, then we have
\begin{align*}
&P(u)=\left(4-u,2-u,2-\frac{u}{2},3-u \right),\\
%(2-u)E_1+ 3E_2 + \left(3-\frac{3}{4}u\right)F_1 + \left(4-\frac{u}{2}\right)F_2 + \left(2-\frac{u}{4}\right)F_3. \\
&N(u)=\left(0, 0 , \frac{u}{2},0 \right).
%\frac{3}{4}uF_1 + \frac{u}{2}F_2 +  \frac{u}{4}F_3.
\end{align*}
%If $u\in [1,2]$, it holds
%\begin{align*}
%&P(u)=\left(2, 4-u, 2-u ,0, 2-u\right)\\
%(2-u)E_1+ 3E_2 + \left(3-\frac{3}{4}u\right)F_1 + \left(4-\frac{u}{2}\right)F_2 + \left(2-\frac{u}{4}\right)F_3. \\
%&N(u)=\left(0, u-1 , u-1, 0 , 0\right)
%\frac{3}{4}uF_1 + \frac{u}{2}F_2 +  \frac{u}{4}F_3.
%\end{align*}
%If $u\in [2,3]$, it holds
%\begin{align*}
%&P(u)=\left(3-u,3-u,3-u,0, 3-u\right)\\
%(2-u)E_1+ 3E_2 + \left(3-\frac{3}{4}u\right)F_1 + \left(4-\frac{u}{2}\right)F_2 + \left(2-\frac{u}{4}\right)F_3. \\
%&N(u)=\left(u-1, 0 , u-2 , 0, u-1 \right)
%\frac{3}{4}uF_1 + \frac{u}{2}F_2 +  \frac{u}{4}F_3.
%\end{align*}
We note that $-K_{S}-uL$ is not pseudo effective for $u>2$.
If $u\in [0,2]$, then we have 
\begin{align*}
P(u)^2=\frac{(u-2)(u-6)}{2},\quad P(u)L = \frac{4-u}{2}.
\end{align*}
%If $u\in [1,2]$ we have 
%\begin{align*}
%P(u)^2=4(2-u),\quad P(u)F = 2
%\end{align*}
%f $u\in [2,3]$ we have 
%\begin{align*}
%P(u)^2=(3-u)^2,\quad P(u)F = (3-u)
%\end{align*}
Therefore,  we get 
\begin{align*}
S(L)%=\frac{1}{4}\int_{0}^{2}(2-u)(2+u) du
%+ \frac{1}{5}\int_{1}^{2} 2(2-u)^2 du
=\frac{8}{9}, \quad S(W_{\bullet,\bullet}^{L},p)
%=\frac{2}{5}\int_{0}^{1} \Rm{ord}_{p}(N(u)|_{F}) du
%&=%\frac{1}{5}\int_{0}^{1} (1+u)^2 du
%+ \frac{2}{5}\int_{0}^{2}\frac{5}{4}u \cdot \Rm{ord}_{p}(N(u)|_{E_{1}})  du
%+\frac{1}{5}\int_{0}^{2}\frac{25}{16}u^2 du\\
 &=\frac{7}{9}.
% \begin{cases}
 % \frac{17}{15} & \text{if $p \in E_{1}\cap F_{j} $ for $j=1,2$},\\
 %                   & \\
 %  \frac{7}{15}   & \text{if $p \in E_{1} \setminus \bigcup_{j=1}^{2} F_{j}$}, \\
%  \end{cases}
\end{align*}
%by Theorem~\ref{S-inv}.
Hence we have 
\begin{align*}
\frac{9}{8}  \geq \delta_{p}(S) \geq  \Rm{min}\left\{ \frac{1}{S(L)}, \frac{1}{S(W_{\bullet,\bullet}^{L},p)} \right\}
= \frac{9}{8} 
\end{align*}
from Corollary~\ref{AZ}.
Thus, we have $\delta_{p}(S)= 9/8$ in this case.
%We can check $\delta_{p}(S)= 9/11$ for $p \in E_3 $ similarly.
\end{proof}

\begin{prop}
Let $S$ be the anti-canonical degree $6$ weak del Pezzo surface such
that the dual graph of negative curves is
%The dual graph of negative curves of $S$ is 
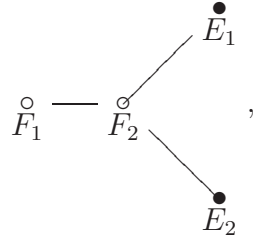
\[
\xygraph{
 %   {\circ}* +!U{F_{1}} - [r]
  %  {\circ}*+!U{F_{2}}  -[r]
 %   {\circ}* +!U{F_{1}} - [r]
   {\circ}* +!D{ }*+!U{F_2}
      (%- [ru] 
 %   {\circ}* +!U{F_{4}} 
   % {\bullet}*+!U{E_{4}}  
   %     - [r] \cdots
   %     - [r] \bullet 
   %     - [r] \bullet 
     %  - [rd] \bullet
              ,
                -[rd]
    {\bullet}*+!U{E_{2}} 
       ,
                -[l]
    {\circ}*+!U{F_{1}} 
     ,
                -[ru]
    {\bullet}*+!U{E_{1}}
   % {\bullet}*+!U{E_{6}}
   %     - [r] \cdots
  %      - [r] \bullet 
  %      - [r] \bullet 
    %     - [ru] {\bullet}*+!U{E_{7}}
         },
        \]
where $E_{i}$ $(i=1,2)$ is a $(-1)$-curve and $F_j$ $(j=1,2)$ is a $(-2)$-curve.
Then, for a point $p \in S$, it holds that
\begin{equation}
\nonumber
  \delta_{p}(S)=
  \begin{cases}
    \frac{3}{4} & \text{if $p \in F_1 \setminus F_2$,} \\
    \frac{3}{5}                 & \text{if $p \in F_2 $,} \\
    \frac{4}{5}                 & \text{if $p \in (E_1 \cup E_2) \setminus F_2$,} \\
%    \frac{3}{4}                 & \text{if $p \in F \setminus E_i$,} \\
  1 & \text{if $p \in S \setminus \left( \bigcup_{i}  E_{i} \cup  F\right)$.}\\
 %  \frac{15}{11} & \text{if $p \in S \setminus \left( \bigcup_{i,j}( E_{i} \cup F_{j})\right)$}
  \end{cases}
\end{equation}
\end{prop}

\begin{proof}
We denote $D=\sum_{i=1,2} a_{i}E_{i} + \sum_{j=1,2}b_jF_{j} \in \Rm{Div}(S)$ ($a_i , b_j \in \Z$)
by $D=(a_1, a_2, b_1,b_2)$.
The intersection matrix of $\{E_1, E_2, F_1,F_2  \}$ is 
\[A:=
\left(
\begin{array}{cc|cc}
  -1 & 0 &  0  &  1       \\ 
  0 & -1 &  0  &  1      \\ \hline
  0 & 0 &  -2  &  1    \\ 
  1 & 1 &  1 &  -2     \\ 
%  0 & 1 &  1  &  0       \\
\end{array}
\right).
\]
 We note that 
 $
 -K_{S} \sim 3E_1 + 3E_2 + 2F_1 + 4F_2 = (3,3,2,4).
$
%%%%%%%%%%%%%%%%%%%%% (1)

\noindent
(1)\;The case $p \in F_1 \setminus F_2$.

\noindent
We calculate $S(F_{1})$ and $S(W_{\bullet,\bullet}^{F_{1}},p)$.
Take $u \in \R_{\geq 0}$. 
Let $P(u)+N(u)$ be the Zariski decomposition of $-K_{S}-uF_{1}$.
%where $P(u)$ is the positive part and $N(u)$ is the negative part.
%
If $u\in [0,2]$, then we have
\begin{align*}
&P(u)=\left(3,3,2-u,4-\frac{u}{2}\right),\\
%(2-u)E_1+ 3E_2 + \left(3-\frac{3}{4}u\right)F_1 + \left(4-\frac{u}{2}\right)F_2 + \left(2-\frac{u}{4}\right)F_3. \\
&N(u)=\left(0, 0 ,0, \frac{u}{2} \right).
%\frac{3}{4}uF_1 + \frac{u}{2}F_2 +  \frac{u}{4}F_3.
\end{align*}
%If $u\in [1,2]$, it holds
%\begin{align*}
%&P(u)=\left(2(3-u),2,2-u,4-u\right)\\
%(2-u)E_1+ 3E_2 + \left(3-\frac{3}{4}u\right)F_1 + \left(4-\frac{u}{2}\right)F_2 + \left(2-\frac{u}{4}\right)F_3. \\
%&N(u)=\left(2(u-1), 0 , 0 , u-1\right)
%\frac{3}{4}uF_1 + \frac{u}{2}F_2 +  \frac{u}{4}F_3.
%\end{align*}
We note that $-K_{S}-uF_{1}$ is not pseudo effective for $u>2$.
If $u\in [0,2]$, then we have
\begin{align*}
P(u)^2=\frac{3(2-u)(2+u)}{2},\quad P(u)F_1 = \frac{3u}{2}.
\end{align*}
%If $u\in [1,2]$ we have 
%\begin{align*}
%P(u)^2=4(2-u),\quad P(u)F_1 = 2
%\end{align*}

Therefore,  we get 
\begin{align*}
S(F_1)%=\frac{1}{4}\int_{0}^{2}(2-u)(2+u) du
%+ \frac{1}{5}\int_{1}^{2} 2(2-u)^2 du
=\frac{4}{3}, \quad S(W_{\bullet,\bullet}^{F_{1}},p)
%=\frac{2}{5}\int_{0}^{1} \Rm{ord}_{p}(N(u)|_{F}) du
%&=%\frac{1}{5}\int_{0}^{1} (1+u)^2 du
%+ \frac{2}{5}\int_{0}^{2}\frac{5}{4}u \cdot \Rm{ord}_{p}(N(u)|_{E_{1}})  du
%+\frac{1}{5}\int_{0}^{2}\frac{25}{16}u^2 du\\
 &=1.
% \begin{cases}
 % \frac{17}{15} & \text{if $p \in E_{1}\cap F_{j} $ for $j=1,2$},\\
 %                   & \\
 %  \frac{7}{15}   & \text{if $p \in E_{1} \setminus \bigcup_{j=1}^{2} F_{j}$}, \\
%  \end{cases}
\end{align*}
Hence we have 
\begin{align*}
\frac{3}{4}  \geq \delta_{p}(S) \geq  \Rm{min}\left\{ \frac{1}{S(F_1)}, \frac{1}{S(W_{\bullet,\bullet}^{F_1},p)} \right\}
= \frac{3}{4} 
\end{align*}
from Corollary~\ref{AZ}.
Thus, we have $\delta_{p}(S)= 3/4$ in this case.

%%%%%%%%%%%%%%%%%%%%% (2)

\noindent
(2)\;The case $p \in F_2$.

\noindent
We calculate $S(F_{2})$ and $S(W_{\bullet,\bullet}^{F_{2}},p)$.
Take $u \in \R_{\geq 0}$. 
Let $P(u)+N(u)$ be the Zariski decomposition of $-K_{S}-uF_{2}$.
%where $P(u)$ is the positive part and $N(u)$ is the negative part.
%
If $u\in [0,1]$, then we have
\begin{align*}
&P(u)=\left(3,3,2-\frac{u}{2},4-u\right),\\
%(2-u)E_1+ 3E_2 + \left(3-\frac{3}{4}u\right)F_1 + \left(4-\frac{u}{2}\right)F_2 + \left(2-\frac{u}{4}\right)F_3. \\
&N(u)=\left(0, 0 , \frac{u}{2} ,0 \right).
%\frac{3}{4}uF_1 + \frac{u}{2}F_2 +  \frac{u}{4}F_3.
\end{align*}
If $u\in [1,4]$, then we have
\begin{align*}
&P(u)=\left(4-u , 4-u ,2- \frac{u}{2} ,4-u\right),\\
%(2-u)E_1+ 3E_2 + \left(3-\frac{3}{4}u\right)F_1 + \left(4-\frac{u}{2}\right)F_2 + \left(2-\frac{u}{4}\right)F_3. \\
&N(u)=\left(u-1, u-1 , \frac{u}{2} , 0\right).
%\frac{3}{4}uF_1 + \frac{u}{2}F_2 +  \frac{u}{4}F_3.
\end{align*}
We note that $-K_{S}-uF_{2}$ is not pseudo effective for $u>4$.
If $u\in [0,1]$, then we have 
\begin{align*}
P(u)^2=6-\frac{3u^2}{2},\quad P(u)F_2 = \frac{3u}{2}.
\end{align*}
If $u\in [1,4]$, then we have
\begin{align*}
P(u)^2=\frac{(4-u)^2}{2},\quad P(u)F_2 = 2-\frac{u}{2}.
\end{align*}

Therefore,  we get 
\begin{align*}
S(F_2)%=\frac{1}{4}\int_{0}^{2}(2-u)(2+u) du
%+ \frac{1}{5}\int_{1}^{2} 2(2-u)^2 du
=\frac{5}{3}, \quad 
S(W_{\bullet,\bullet}^{F_{2}},p)
%=\frac{2}{5}\int_{0}^{1} \Rm{ord}_{p}(N(u)|_{F}) du
%&=%\frac{1}{5}\int_{0}^{1} (1+u)^2 du
%+ \frac{2}{5}\int_{0}^{2}\frac{5}{4}u \cdot \Rm{ord}_{p}(N(u)|_{E_{1}})  du
%+\frac{1}{5}\int_{0}^{2}\frac{25}{16}u^2 du\\
 &=
 \begin{cases}
  \frac{4}{3} & \text{if $p \in  F_{2}\cap F_{1} $},\\
                    & \\
   \frac{5}{4}   & \text{if $p \in F_{2} \cap ( E_1 \cup E_2)$}, \\
              & \\
   \frac{1}{2}   & \text{if $p \in F_{2} \setminus  (F_{1}\cup E_1 \cup E_2)$}. \\
  \end{cases}
\end{align*}
Hence we have 
\begin{align*}
\frac{3}{5}  \geq \delta_{p}(S) \geq  \Rm{min}\left\{ \frac{1}{S(F_2)}, \frac{1}{S(W_{\bullet,\bullet}^{F_2},p)} \right\}
= \frac{3}{5} 
\end{align*}
from Corollary~\ref{AZ}.
Thus, we have $\delta_{p}(S)= 3/5$ in this case.

%%%%%%%%%%%%%%%%%%%%% (3)

\noindent
(3)\;The case $p \in E_1 \setminus F_2$.

\noindent
We calculate $S(E_{1})$ and $S(W_{\bullet,\bullet}^{E_{1}},p)$.
Take $u \in \R_{\geq 0}$. 
Let $P(u)+N(u)$ be the Zariski decomposition of $-K_{S}-uE_{1}$.
%where $P(u)$ is the positive part and $N(u)$ is the negative part.
%
If $u\in [0,\frac{3}{2}]$, then we have
\begin{align*}
&P(u)=\left(3-u,3,2-\frac{u}{3},4-\frac{2}{3}u \right),\\
%(2-u)E_1+ 3E_2 + \left(3-\frac{3}{4}u\right)F_1 + \left(4-\frac{u}{2}\right)F_2 + \left(2-\frac{u}{4}\right)F_3. \\
&N(u)=\left(0, 0 , \frac{u}{3} , \frac{2}{3}u \right).
%\frac{3}{4}uF_1 + \frac{u}{2}F_2 +  \frac{u}{4}F_3.
\end{align*}
If $u\in [\frac{3}{2},3]$, then we have
\begin{align*}
&P(u)=\left(3-u , 2(3-u) , 3-u  , 2(3-u)\right),\\
%(2-u)E_1+ 3E_2 + \left(3-\frac{3}{4}u\right)F_1 + \left(4-\frac{u}{2}\right)F_2 + \left(2-\frac{u}{4}\right)F_3. \\
&N(u)=\left(0, 2u-3 , u-1 , 2(u-1) \right).
%\frac{3}{4}uF_1 + \frac{u}{2}F_2 +  \frac{u}{4}F_3.
\end{align*}
We note that $-K_{S}-uE_{1}$ is not pseudo effective for $u>3$.
If $u\in [0,\frac{3}{2}]$, then we have 
\begin{align*}
P(u)^2=6-2u-\frac{u^2}{3},\quad P(u)E_1 =1+ \frac{u}{3}.
\end{align*}
If $u\in [\frac{3}{2},3]$, then we have
\begin{align*}
P(u)^2=(3-u)^2,\quad P(u)E_1 = 3-u.
\end{align*}

Therefore,  we get 
\begin{align*}
S(E_1)%=\frac{1}{4}\int_{0}^{2}(2-u)(2+u) du
%+ \frac{1}{5}\int_{1}^{2} 2(2-u)^2 du
=\frac{5}{4}, \quad
S(W_{\bullet,\bullet}^{E_1},p)
%=\frac{2}{5}\int_{0}^{1} \Rm{ord}_{p}(N(u)|_{F}) du
%&=%\frac{1}{5}\int_{0}^{1} (1+u)^2 du
%+ \frac{2}{5}\int_{0}^{2}\frac{5}{4}u \cdot \Rm{ord}_{p}(N(u)|_{E_{1}})  du
%+\frac{1}{5}\int_{0}^{2}\frac{25}{16}u^2 du\\
 &=\frac{7}{12}.
 %\begin{cases}
  %\frac{4}{3} & \text{if $p \in  F_{2}\cap F_{1} $},\\
    %                & \\
   %\frac{5}{4}   & \text{if $p \in F_{2} \cap ( E_1 \cup E_2)$}, \\
    %          & \\
  % \frac{1}{2}   & \text{if $p \in F_{2} \setminus  (F_{1}\cup E_1 \cup E_2)$}, \\
%  \end{cases}
\end{align*}
Hence we have 
\begin{align*}
\frac{4}{5}  \geq \delta_{p}(S) \geq  \Rm{min}\left\{ \frac{1}{S(E_1)}, \frac{1}{S(W_{\bullet,\bullet}^{E_1},p)} \right\}
= \frac{4}{5} 
\end{align*}
from Corollary~\ref{AZ}.
Thus, we have $\delta_{p}(S)= 4/5$ in this case.

%%%%%%%%%%%%%%%%%%%%% (4)

\noindent
(4)\;The case $p \in S \setminus \left( E_1\cup E_2 \cup F_1 \cup F_2\right)$.

\noindent
Let $L \in |E_1+E_2+F_2|$ be a smooth irreducible curve.
We calculate $S(L)$ and $S(W_{\bullet,\bullet}^{L},p)$.
Take $u \in \R_{\geq 0}$. 
Let $P(u)+N(u)$ be the Zariski decomposition of $-K_{S}-uL$.
%where $P(u)$ is the positive part and $N(u)$ is the negative part.
%
If $u\in [0,3]$, then we have
\begin{align*}
&P(u)=\left(3-u,3-u,\frac{2}{3}(3-u),\frac{4}{3}(3-u) \right),\\
%(2-u)E_1+ 3E_2 + \left(3-\frac{3}{4}u\right)F_1 + \left(4-\frac{u}{2}\right)F_2 + \left(2-\frac{u}{4}\right)F_3. \\
&N(u)=\left(0, 0 , \frac{2u}{3} , \frac{u}{3} \right).
%\frac{3}{4}uF_1 + \frac{u}{2}F_2 +  \frac{u}{4}F_3.
\end{align*}
%If $u\in [1,2]$, it holds
%\begin{align*}
%&P(u)=\left(2, 4-u, 2-u ,0, 2-u\right)\\
%(2-u)E_1+ 3E_2 + \left(3-\frac{3}{4}u\right)F_1 + \left(4-\frac{u}{2}\right)F_2 + \left(2-\frac{u}{4}\right)F_3. \\
%&N(u)=\left(0, u-1 , u-1, 0 , 0\right)
%\frac{3}{4}uF_1 + \frac{u}{2}F_2 +  \frac{u}{4}F_3.
%\end{align*}
%If $u\in [2,3]$, it holds
%\begin{align*}
%&P(u)=\left(3-u,3-u,3-u,0, 3-u\right)\\
%(2-u)E_1+ 3E_2 + \left(3-\frac{3}{4}u\right)F_1 + \left(4-\frac{u}{2}\right)F_2 + \left(2-\frac{u}{4}\right)F_3. \\
%&N(u)=\left(u-1, 0 , u-2 , 0, u-1 \right)
%\frac{3}{4}uF_1 + \frac{u}{2}F_2 +  \frac{u}{4}F_3.
%\end{align*}
We note that $-K_{S}-uL$ is not pseudo effective for $u>3$.
If $u\in [0,3]$, then we have
\begin{align*}
P(u)^2=\frac{2(3-u)^2}{3},\quad P(u)L = \frac{2(3-u)}{3}.
\end{align*}
%If $u\in [1,2]$ we have 
%\begin{align*}
%P(u)^2=4(2-u),\quad P(u)F = 2
%\end{align*}
%f $u\in [2,3]$ we have 
%\begin{align*}
%P(u)^2=(3-u)^2,\quad P(u)F = (3-u)
%\end{align*}
Therefore,  we get 
\begin{align*}
S(L)%=\frac{1}{4}\int_{0}^{2}(2-u)(2+u) du
%+ \frac{1}{5}\int_{1}^{2} 2(2-u)^2 du
=1, \quad
S(W_{\bullet,\bullet}^{L},p)
%=\frac{2}{5}\int_{0}^{1} \Rm{ord}_{p}(N(u)|_{F}) du
%&=%\frac{1}{5}\int_{0}^{1} (1+u)^2 du
%+ \frac{2}{5}\int_{0}^{2}\frac{5}{4}u \cdot \Rm{ord}_{p}(N(u)|_{E_{1}})  du
%+\frac{1}{5}\int_{0}^{2}\frac{25}{16}u^2 du\\
 &=\frac{2}{3}.
% \begin{cases}
 % \frac{17}{15} & \text{if $p \in E_{1}\cap F_{j} $ for $j=1,2$},\\
 %                   & \\
 %  \frac{7}{15}   & \text{if $p \in E_{1} \setminus \bigcup_{j=1}^{2} F_{j}$}, \\
%  \end{cases}
\end{align*}
Hence we have 
\begin{align*}
1 \geq \delta_{p}(S) \geq  \Rm{min}\left\{ \frac{1}{S(L)}, \frac{1}{S(W_{\bullet,\bullet}^{L},p)} \right\}
= 1
\end{align*}
from Corollary~\ref{AZ}.
Thus, we have $\delta_{p}(S)= 1$ in this case.
%We can check $\delta_{p}(S)= 9/11$ for $p \in E_3 $ similarly.
\end{proof}

\begin{prop}
Let $S$ be the anti-canonical degree $6$ weak del Pezzo surface such
that the dual graph of negative curves is
\[
\xygraph{
   {\circ}* +!D{ }*+!U{F_1}
                    -[r]
    {\circ}*+!U{F_{2}}  - [r]
    {\bullet}*+!U{E}
   %     - [r] \cdots
%       - [r] {\circ}* +!U{F_{2}}
       - [r] {\circ}* +!U{F_{3}}
       %  - [r] {\bullet}*+!U{E_{4}}
              },
        \]
where $E$ is a $(-1)$-curve and $F_j$ $(j=1,2,3)$ is a $(-2)$-curve.
Then, for a point $p \in S$, it holds that
\begin{equation}
\nonumber
  \delta_{p}(S)=
  \begin{cases}
    \frac{3}{4} & \text{if $p \in F_1 \setminus F_2$,} \\
    \frac{3}{5}                 & \text{if $ p \in F_2 \setminus E $,} \\
     \frac{1}{2}                 & \text{if $p \in E $,} \\
     \frac{3}{4} & \text{if $p \in  F_3 \setminus E $,} \\
     1 & \text{if $p \in S \setminus \left( E \cup F_1 \cup F_2 \cup F_3 \right)$.}
  \end{cases}
\end{equation}
\end{prop}

\begin{proof}
We denote $D= a E + \sum_{j=1,2,3}b_jF_{j} \in \Rm{Div}(S)$ ($a , b_j \in \Z$)
by $D=(a,  b_1,b_2, b_3)$.
The intersection matrix of $\{E, F_1,F_2, F_3  \}$ is 
\[A:=
\left(
\begin{array}{c|ccc}
  -1 & 0 &  1  &  1       \\ \hline
  0 & -2 &  1  &  0      \\ 
  1 & 1 &  -2  &  0    \\ 
  1 & 0 &  0 &  -2     \\ 
%  0 & 1 &  1  &  0       \\
\end{array}
\right).
\]
 We note that 
 $
 -K_{S} \sim 6E + 2F_1 + 4F_2 + 3F_3 = (6,2,4,3).
$
%%%%%%%%%%%%%%%%%%%%% (1)

\noindent
(1)\;The case $p \in F_1 \setminus F_2$.

\noindent
We calculate $S(F_{1})$ and $S(W_{\bullet,\bullet}^{F_{1}},p)$.
Take $u \in \R_{\geq 0}$. 
Let $P(u)+N(u)$ be the Zariski decomposition of $-K_{S}-uF_{1}$.
%where $P(u)$ is the positive part and $N(u)$ is the negative part.
%
If $u\in [0,2]$, then we have
\begin{align*}
&P(u)=\left(6,2-u,4-\frac{u}{2},3\right),\\
%(2-u)E_1+ 3E_2 + \left(3-\frac{3}{4}u\right)F_1 + \left(4-\frac{u}{2}\right)F_2 + \left(2-\frac{u}{4}\right)F_3. \\
&N(u)=\left(0, 0 , \frac{u}{2} , 0\right).
%\frac{3}{4}uF_1 + \frac{u}{2}F_2 +  \frac{u}{4}F_3.
\end{align*}
%If $u\in [1,2]$, it holds
%\begin{align*}
%&P(u)=\left(2(3-u),2,2-u,4-u\right)\\
%(2-u)E_1+ 3E_2 + \left(3-\frac{3}{4}u\right)F_1 + \left(4-\frac{u}{2}\right)F_2 + \left(2-\frac{u}{4}\right)F_3. \\
%&N(u)=\left(2(u-1), 0 , 0 , u-1\right)
%\frac{3}{4}uF_1 + \frac{u}{2}F_2 +  \frac{u}{4}F_3.
%\end{align*}
We note that $-K_{S}-uF_{1}$ is not pseudo effective for $u>2$.
If $u\in [0,2]$, then we have
\begin{align*}
P(u)^2=\frac{3(4-u^2)}{2},\quad P(u)F_1 = \frac{3u}{2}.
\end{align*}
%If $u\in [1,2]$ we have 
%\begin{align*}
%P(u)^2=4(2-u),\quad P(u)F_1 = 2
%\end{align*}

Therefore, we get 
\begin{align*}
S(F_1)%=\frac{1}{4}\int_{0}^{2}(2-u)(2+u) du
%+ \frac{1}{5}\int_{1}^{2} 2(2-u)^2 du
=\frac{4}{3}, \quad 
S(W_{\bullet,\bullet}^{F_{1}},p)
%=\frac{2}{5}\int_{0}^{1} \Rm{ord}_{p}(N(u)|_{F}) du
%&=%\frac{1}{5}\int_{0}^{1} (1+u)^2 du
%+ \frac{2}{5}\int_{0}^{2}\frac{5}{4}u \cdot \Rm{ord}_{p}(N(u)|_{E_{1}})  du
%+\frac{1}{5}\int_{0}^{2}\frac{25}{16}u^2 du\\
 &=1.
% \begin{cases}
 % \frac{17}{15} & \text{if $p \in E_{1}\cap F_{j} $ for $j=1,2$},\\
 %                   & \\
 %  \frac{7}{15}   & \text{if $p \in E_{1} \setminus \bigcup_{j=1}^{2} F_{j}$}, \\
%  \end{cases}
\end{align*}
Hence we have 
\begin{align*}
\frac{3}{4}  \geq \delta_{p}(S) \geq  \Rm{min}\left\{ \frac{1}{S(F_1)}, \frac{1}{S(W_{\bullet,\bullet}^{F_1},p)} \right\}
= \frac{3}{4} 
\end{align*}
from Corollary~\ref{AZ}.
Thus, we have $\delta_{p}(S)= 3/4$ in this case.

%%%%%%%%%%%%%%%%%%%%% (2)

\noindent
(2)\;The case $p \in F_2 \setminus E$.

\noindent
We calculate $S(F_{2})$ and $S(W_{\bullet,\bullet}^{F_{2}},p)$.
Take $u \in \R_{\geq 0}$. 
Let $P(u)+N(u)$ be the Zariski decomposition of $-K_{S}-uF_{2}$.
%where $P(u)$ is the positive part and $N(u)$ is the negative part.
%
If $u\in [0,1]$, then we have
\begin{align*}
&P(u)=\left(6,2-\frac{u}{2},4-u,3\right),\\
%(2-u)E_1+ 3E_2 + \left(3-\frac{3}{4}u\right)F_1 + \left(4-\frac{u}{2}\right)F_2 + \left(2-\frac{u}{4}\right)F_3. \\
&N(u)=\left(0, \frac{u}{2},0 ,0 \right).
%\frac{3}{4}uF_1 + \frac{u}{2}F_2 +  \frac{u}{4}F_3.
\end{align*}
If $u\in [1,4]$, then we have
\begin{align*}
&P(u)=\left(2(4-u),2-\frac{u}{2},4-u,4-u\right),\\
%(2-u)E_1+ 3E_2 + \left(3-\frac{3}{4}u\right)F_1 + \left(4-\frac{u}{2}\right)F_2 + \left(2-\frac{u}{4}\right)F_3. \\
&N(u)=\left(2(u-1), \frac{u}{2},0 ,u-1 \right).
%\frac{3}{4}uF_1 + \frac{u}{2}F_2 +  \frac{u}{4}F_3.
\end{align*}
We note that $-K_{S}-uF_{2}$ is not pseudo effective for $u>4$.
If $u\in [0,1]$, then we have
\begin{align*}
P(u)^2=6-\frac{3u^2}{2},\quad P(u)F_2 = \frac{3u}{2}.
\end{align*}
If $u\in [1,4]$, then we have
\begin{align*}
P(u)^2=\frac{(4-u)^2}{2},\quad P(u)F_2 = 2-\frac{u}{2}.
\end{align*}

Therefore,  we get 
\begin{align*}
S(F_2)%=\frac{1}{4}\int_{0}^{2}(2-u)(2+u) du
%+ \frac{1}{5}\int_{1}^{2} 2(2-u)^2 du
=\frac{5}{3}, \quad 
S(W_{\bullet,\bullet}^{F_{2}},p)
%=\frac{2}{5}\int_{0}^{1} \Rm{ord}_{p}(N(u)|_{F}) du
%&=%\frac{1}{5}\int_{0}^{1} (1+u)^2 du
%+ \frac{2}{5}\int_{0}^{2}\frac{5}{4}u \cdot \Rm{ord}_{p}(N(u)|_{E_{1}})  du
%+\frac{1}{5}\int_{0}^{2}\frac{25}{16}u^2 du\\
 &=
 \begin{cases}
  \frac{4}{3} & \text{if $p \in  F_{2}\cap F_{1} $},\\
                    & \\
   \frac{1}{2}   & \text{if $p \in F_{2} \setminus  (F_{1}\cup E)$}. \\
  \end{cases}
\end{align*}
Hence we have 
\begin{align*}
\frac{3}{5}  \geq \delta_{p}(S) \geq  \Rm{min}\left\{ \frac{1}{S(F_2)}, \frac{1}{S(W_{\bullet,\bullet}^{F_2},p)} \right\}
= \frac{3}{5} 
\end{align*}
from Corollary~\ref{AZ}.
Thus, we have $\delta_{p}(S)= 3/5$ in this case.

%%%%%%%%%%%%%%%%%%%%% (3)

\noindent
(3)\;The case $p \in E$.

\noindent
We calculate $S(E)$ and $S(W_{\bullet,\bullet}^{E},p)$.
Take $u \in \R_{\geq 0}$. 
Let $P(u)+N(u)$ be the Zariski decomposition of $-K_{S}-uE$.
%where $P(u)$ is the positive part and $N(u)$ is the negative part.
%
If $u\in [0,6]$, then we have
\begin{align*}
&P(u)=\left(6-u,\frac{1}{3}(6-u),\frac{2}{3}(6-u),\frac{1}{2}(6-u)\right),\\
%(2-u)E_1+ 3E_2 + \left(3-\frac{3}{4}u\right)F_1 + \left(4-\frac{u}{2}\right)F_2 + \left(2-\frac{u}{4}\right)F_3. \\
&N(u)=\left(0, \frac{u}{3} , \frac{2u}{3} , \frac{u}{2}\right).
%\frac{3}{4}uF_1 + \frac{u}{2}F_2 +  \frac{u}{4}F_3.
\end{align*}
%If $u\in [1,2]$, it holds
%\begin{align*}
%&P(u)=\left(2(3-u),2,2-u,4-u\right)\\
%(2-u)E_1+ 3E_2 + \left(3-\frac{3}{4}u\right)F_1 + \left(4-\frac{u}{2}\right)F_2 + \left(2-\frac{u}{4}\right)F_3. \\
%&N(u)=\left(2(u-1), 0 , 0 , u-1\right)
%\frac{3}{4}uF_1 + \frac{u}{2}F_2 +  \frac{u}{4}F_3.
%\end{align*}
We note that $-K_{S}-uE$ is not pseudo effective for $u>6$.
If $u\in [0,6]$, then we have 
\begin{align*}
P(u)^2=\frac{1}{6}(6-u^2),\quad P(u)E = \frac{1}{6}(6-u).
\end{align*}
%If $u\in [1,2]$ we have 
%\begin{align*}
%P(u)^2=4(2-u),\quad P(u)F_1 = 2
%\end{align*}

Therefore,  we get 
\begin{align*}
S(E)%=\frac{1}{4}\int_{0}^{2}(2-u)(2+u) du
%+ \frac{1}{5}\int_{1}^{2} 2(2-u)^2 du
=2, \quad 
S(W_{\bullet,\bullet}^{E},p)
%=\frac{2}{5}\int_{0}^{1} \Rm{ord}_{p}(N(u)|_{F}) du
%&=%\frac{1}{5}\int_{0}^{1} (1+u)^2 du
%+ \frac{2}{5}\int_{0}^{2}\frac{5}{4}u \cdot \Rm{ord}_{p}(N(u)|_{E_{1}})  du
%+\frac{1}{5}\int_{0}^{2}\frac{25}{16}u^2 du\\
 &=
 \begin{cases}
  \frac{5}{3} & \text{if $p \in E \cap F_{2}$},\\
                    & \\
   \frac{4}{3} & \text{if $p \in E \cap F_{3} $},\\
                    & \\
   \frac{1}{3}   & \text{if $p \in E \setminus (F_2 \cup F_3)$}. \\
  \end{cases}
\end{align*}
Hence we have 
\begin{align*}
\frac{1}{2}  \geq \delta_{p}(S) \geq  \Rm{min}\left\{ \frac{1}{S(E)}, \frac{1}{S(W_{\bullet,\bullet}^{E},p)} \right\}
= \frac{1}{2} 
\end{align*}
from Corollary~\ref{AZ}.
Thus, we have $\delta_{p}(S)= 1/2$ in this case.

%%%%%%%%%%%%%%%%%%%%% (4)

\noindent
(4)\;The case $p \in F_3 \setminus E$.

\noindent
We calculate $S(F_{3})$ and $S(W_{\bullet,\bullet}^{F_{3}},p)$.
Take $u \in \R_{\geq 0}$. 
Let $P(u)+N(u)$ be the Zariski decomposition of $-K_{S}-uF_{3}$.
%where $P(u)$ is the positive part and $N(u)$ is the negative part.
%
If $u\in [0,1]$, then we have
\begin{align*}
&P(u)=\left(6,2,4,3-u \right),\\
%(2-u)E_1+ 3E_2 + \left(3-\frac{3}{4}u\right)F_1 + \left(4-\frac{u}{2}\right)F_2 + \left(2-\frac{u}{4}\right)F_3. \\
&N(u)=\left(0,0 ,0 ,0 \right).
%\frac{3}{4}uF_1 + \frac{u}{2}F_2 +  \frac{u}{4}F_3.
\end{align*}
If $u\in [1,3]$, then we have
\begin{align*}
&P(u)=\left(3(3-u), 3-u , 2(3-u)  ,3-u\right),\\
%(2-u)E_1+ 3E_2 + \left(3-\frac{3}{4}u\right)F_1 + \left(4-\frac{u}{2}\right)F_2 + \left(2-\frac{u}{4}\right)F_3. \\
&N(u)=\left(3(u-1), u-1 , 2(u-1), 0 \right).
%\frac{3}{4}uF_1 + \frac{u}{2}F_2 +  \frac{u}{4}F_3.
\end{align*}
We note that $-K_{S}-uF_{3}$ is not pseudo effective for $u>3$.
If $u\in [0,1]$, then we have
\begin{align*}
P(u)^2=6-2u^2 ,\quad P(u)F_3 = 2u.
\end{align*}
If $u\in [1,3]$, then we have
\begin{align*}
P(u)^2=(3-u)^2,\quad P(u)F_3 = 3-u .
\end{align*}

Therefore,  we get 
\begin{align*}
S(F_3)%=\frac{1}{4}\int_{0}^{2}(2-u)(2+u) du
%+ \frac{1}{5}\int_{1}^{2} 2(2-u)^2 du
=\frac{4}{3}, \quad 
S(W_{\bullet,\bullet}^{F_{3}},p)
%=\frac{2}{5}\int_{0}^{1} \Rm{ord}_{p}(N(u)|_{F}) du
%&=%\frac{1}{5}\int_{0}^{1} (1+u)^2 du
%+ \frac{2}{5}\int_{0}^{2}\frac{5}{4}u \cdot \Rm{ord}_{p}(N(u)|_{E_{1}})  du
%+\frac{1}{5}\int_{0}^{2}\frac{25}{16}u^2 du\\
 &=\frac{2}{3}.
% \begin{cases}
 % \frac{4}{3} & \text{if $p \in  F_{2}\cap F_{1} $},\\
  %                  & \\
  % \frac{1}{2}   & \text{if $p \in F_{2} \setminus  (F_{1}\cup E)$}, \\
 % \end{cases}
\end{align*}
Hence we have 
\begin{align*}
\frac{3}{4}  \geq \delta_{p}(S) \geq  \Rm{min}\left\{ \frac{1}{S(F_3)}, \frac{1}{S(W_{\bullet,\bullet}^{F_3},p)} \right\}
= \frac{3}{4} 
\end{align*}
from Corollary~\ref{AZ}.
Thus, we have $\delta_{p}(S)= 3/4$ in this case.

%%%%%%%%%%%%%%%%%%%%% (5)

\noindent
(5)\;The case $p \in S \setminus \left( E \cup F_1 \cup F_2 \cup F_3\right)$.

\noindent
Let $L \in |2E+F_2+F_3|$ be a smooth irreducible curve.
We calculate $S(L)$ and $S(W_{\bullet,\bullet}^{L},p)$.
Take $u \in \R_{\geq 0}$. 
Let $P(u)+N(u)$ be the Zariski decomposition of $-K_{S}-uL$.
%where $P(u)$ is the positive part and $N(u)$ is the negative part.
%
If $u\in [0,3]$, then we have
\begin{align*}
&P(u)=\left(2(3-u) ,\frac{2}{3}(3-u)  , \frac{4}{3}(3-u) , 3-u \right),\\
%(2-u)E_1+ 3E_2 + \left(3-\frac{3}{4}u\right)F_1 + \left(4-\frac{u}{2}\right)F_2 + \left(2-\frac{u}{4}\right)F_3. \\
&N(u)=\left(0, \frac{2u}{3} , \frac{u}{3}, 0 \right).
%\frac{3}{4}uF_1 + \frac{u}{2}F_2 +  \frac{u}{4}F_3.
\end{align*}
%If $u\in [1,2]$, it holds
%\begin{align*}
%&P(u)=\left(2, 4-u, 2-u ,0, 2-u\right)\\
%(2-u)E_1+ 3E_2 + \left(3-\frac{3}{4}u\right)F_1 + \left(4-\frac{u}{2}\right)F_2 + \left(2-\frac{u}{4}\right)F_3. \\
%&N(u)=\left(0, u-1 , u-1, 0 , 0\right)
%\frac{3}{4}uF_1 + \frac{u}{2}F_2 +  \frac{u}{4}F_3.
%\end{align*}
%If $u\in [2,3]$, it holds
%\begin{align*}
%&P(u)=\left(3-u,3-u,3-u,0, 3-u\right)\\
%(2-u)E_1+ 3E_2 + \left(3-\frac{3}{4}u\right)F_1 + \left(4-\frac{u}{2}\right)F_2 + \left(2-\frac{u}{4}\right)F_3. \\
%&N(u)=\left(u-1, 0 , u-2 , 0, u-1 \right)
%\frac{3}{4}uF_1 + \frac{u}{2}F_2 +  \frac{u}{4}F_3.
%\end{align*}
We note that $-K_{S}-uL$ is not pseudo effective for $u>3$.
If $u\in [0,3]$, then we have
\begin{align*}
P(u)^2=\frac{2(3-u)^2}{3},\quad P(u)L = \frac{2(3-u)}{3}.
\end{align*}
%If $u\in [1,2]$ we have 
%\begin{align*}
%P(u)^2=4(2-u),\quad P(u)F = 2
%\end{align*}
%f $u\in [2,3]$ we have 
%\begin{align*}
%P(u)^2=(3-u)^2,\quad P(u)F = (3-u)
%\end{align*}
Therefore,  we get 
\begin{align*}
S(L)%=\frac{1}{4}\int_{0}^{2}(2-u)(2+u) du
%+ \frac{1}{5}\int_{1}^{2} 2(2-u)^2 du
=1, \quad 
S(W_{\bullet,\bullet}^{L},p)
%=\frac{2}{5}\int_{0}^{1} \Rm{ord}_{p}(N(u)|_{F}) du
%&=%\frac{1}{5}\int_{0}^{1} (1+u)^2 du
%+ \frac{2}{5}\int_{0}^{2}\frac{5}{4}u \cdot \Rm{ord}_{p}(N(u)|_{E_{1}})  du
%+\frac{1}{5}\int_{0}^{2}\frac{25}{16}u^2 du\\
 &=\frac{2}{3}.
% \begin{cases}
 % \frac{17}{15} & \text{if $p \in E_{1}\cap F_{j} $ for $j=1,2$},\\
 %                   & \\
 %  \frac{7}{15}   & \text{if $p \in E_{1} \setminus \bigcup_{j=1}^{2} F_{j}$}, \\
%  \end{cases}
\end{align*}
Hence we have 
\begin{align*}
1 \geq \delta_{p}(S) \geq  \Rm{min}\left\{ \frac{1}{S(L)}, \frac{1}{S(W_{\bullet,\bullet}^{L},p)} \right\}
= 1
\end{align*}
from Corollary~\ref{AZ}.
Thus, we have $\delta_{p}(S)= 1$ in this case.
\end{proof}
\begin{prop}
Let $S$ be the del Pezzo surface with the anti-canonical degree $6$.
The dual graph of negative curves of $S$ is 
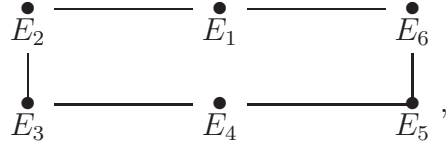
\[
\xygraph{
   {\bullet}* +!D{ }*+!U{E_{3}}
      (- [u] 
    {\bullet}* +!U{E_{2}} - [rr]
    {\bullet}*+!U{E_{1}}  -[rr]
       {\bullet}*+!U{E_6}  
   %     - [r] \cdots
   %     - [r] \bullet 
   %     - [r] \bullet 
       - [d] {\bullet}*+!U{E_{5}}
    ,
                -[rr]
    {\bullet}*+!U{E_{4}} - [rr]
%    {\bullet}*+!U{E_{4}}
   %     - [r] \cdots
  %      - [r] \bullet 
  %      - [r] \bullet 
     %    - [ru] \bullet
         },
        \]
where $E_{i}$ $(i=1, \cdots, 6)$ is a $(-1)$-curve.
Then, for a point $p \in S$, it holds that
%\begin{table}[htb]
 % \begin{tabular}{|ll|c|c|c|} \hline
 % & & Order & Length of non-trivial orbits &\\ \hline 
 % $\Z_{l}$ &Cyclic group &$ l $  & $1$, $1$ &\\ \hline
 % \end{tabular}
%\end{table}
\begin{equation}
\nonumber
  \delta_{p}(S)=
  \begin{cases}
    1 & \text{if $p \in E_i$ $(i=1, \cdots 6)$,} \\
    \frac{6}{5}                 & \text{if $p \in S \setminus \bigcup E_i$.} \\
  \end{cases}
\end{equation}
\end{prop}

\begin{proof}
We can assume that we get $S$ from $\Bb{P}^2$ as follows.
Take three non-colinear points $q_1, q_2, q_3 \in \Bb{P}^2$.
%and a line $l$ passing through these points.
Then we have $\rho:S=\Rm{Bl}_{\{q_1, q_2, q_3\}}\Bb{P}^{2} \to \Bb{P}^2$.
Moreover, we have $E_{2i-1} := \rho^{-1}(q_i)$ ($i=1,2,3$), 
$E_2=\rho^{-1}_{\ast}\Ov{q_1q_2}$, $E_4=\rho^{-1}_{\ast}\Ov{q_2q_3}$, 
and $E_6=\rho^{-1}_{\ast}\Ov{q_3q_1}$.
We denote $D=  \sum_{i=1}^{6}a_{i}E_{i} \in \Rm{Div}(S)$ ($a_i \in \Z$)
by $D=(a_1, a_2, a_3, a_4, a_5, a_6)$.
The intersection matrix of $\{E_1, E_2, E_3, E_4, E_5, E_6  \}$ is 
\[A:=
\left(
\begin{array}{cccccc}
  -1 & 1 & 0  & 0  &  0  &  1     \\
  1 & -1 & 1  & 0  &  0  &  0     \\
  0 & 1 & -1  & 1  &  0 &   0    \\ 
  0 & 0 & 1 & -1  &  1 &    0     \\ 
  0 & 0 & 0  & 1 &  -1  &   1     \\ 
  1 & 0 & 0  & 0  &  1  &  -1   \\
\end{array}
\right).
\]
We note that 
 $
 -K_{S} \sim  (2,2,1,0,0,1).
$
%%%%%%%%%%%%%%%%%%%%% (1)

\noindent
(1)\;The case $p \in E_1$.

\noindent
We calculate $S(E_{1})$.
Take $u \in \R_{\geq 0}$. 
Let $P(u)+N(u)$ be the Zariski decomposition of $-K_{S}-uE_{1}$.
%where $P(u)$ is the positive part and $N(u)$ is the negative part.
%
If $u\in [0,1]$, then we have
\begin{align*}
&P(u)=\left(2-u,2,1,0,0,1\right),\\
%(2-u)E_1+ 3E_2 + \left(3-\frac{3}{4}u\right)F_1 + \left(4-\frac{u}{2}\right)F_2 + \left(2-\frac{u}{4}\right)F_3. \\
&N(u)=\left(0, 0 ,0, 0,0, 0\right).
%\frac{3}{4}uF_1 + \frac{u}{2}F_2 +  \frac{u}{4}F_3.
\end{align*}
If $u\in [1,2]$, then we have
\begin{align*}
&P(u)=\left(2-u,3-u,1,0,0,2-u\right),\\
%(2-u)E_1+ 3E_2 + \left(3-\frac{3}{4}u\right)F_1 + \left(4-\frac{u}{2}\right)F_2 + \left(2-\frac{u}{4}\right)F_3. \\
&N(u)=\left(0, u-1 , 0, 0 , 0, u-1\right).
%\frac{3}{4}uF_1 + \frac{u}{2}F_2 +  \frac{u}{4}F_3.
\end{align*}
%We note that $-K_{S}-uE_{1}$ is not pseudo effective for $u>2$.
If $u\in [0,1]$, then we have
\begin{align*}
P(u)^2=(6-2u-u^2).%\quad P(u)F = 1+u.
\end{align*}
If $u\in [1,2]$, then we have
\begin{align*}
P(u)^2=8-6u+u^{2}.%\quad P(u)F = 3-u.
\end{align*}
Therefore,  we get 
\begin{align*}
S(E_1)%=\frac{1}{4}\int_{0}^{2}(2-u)(2+u) du
%+ \frac{1}{5}\int_{1}^{2} 2(2-u)^2 du
\geq 1
\end{align*}
by the definition of $S(E_1)$. 
%\begin{align*}
%S(W_{\bullet,\bullet}^{E_{1}},p)
%=\frac{2}{5}\int_{0}^{1} \Rm{ord}_{p}(N(u)|_{F}) du
%&=%\frac{1}{5}\int_{0}^{1} (1+u)^2 du
%+ \frac{2}{5}\int_{0}^{2}\frac{5}{4}u \cdot \Rm{ord}_{p}(N(u)|_{E_{1}})  du
%+\frac{1}{5}\int_{0}^{2}\frac{25}{16}u^2 du\\
 %&=\frac{7}{9}
% \begin{cases}
 % \frac{17}{15} & \text{if $p \in E_{1}\cap F_{j} $ for $j=1,2$},\\
 %                   & \\
 %  \frac{7}{15}   & \text{if $p \in E_{1} \setminus \bigcup_{j=1}^{2} F_{j}$}, \\
%  \end{cases}
%\end{align*}
%by \Red{Theorem~$S(W_{\bullet,\bullet}^{F},p)$の明示式}. 
By \cite[\S 2]{FAND}, the (global) delta invariant of del Pezzo surfaces with the anti-canonical degree $6$ is $1$.
Hence we have $\delta_{p}(S)= 1$ for $p \in E_1$.
It also holds $\delta_{p}(S)= 1$ for $p \in E_i$ for $i=1, \cdots, 6$.

%%%%%%%%%%%%%%%%%%%%% (2)

\noindent
(2)\;The case $p \in S \setminus \bigcup_{i} E_i$.

Consider a blowing up $\sigma:\Ti{S} \to S$ at $p$.
Let $E$ be the exceptional curve of $p$,
 $\Ti{E}_{i}$ and $\Ti{F}$ be the proper transform of $E_{i}$ and $F$, respectively.
Take three $(-1)$-curves $G_i := (\rho\sigma)_{\ast}^{-1}(\Ov{\rho\sigma(p)q_i})$ 
for $i=1,2,3$.
We note that $\sigma^{\ast}(-K_{S})\sim G_1 + G_2 +G_3 +3E$.
Hence we have 
\[
\sigma^{\ast}(-K_{S})-uE \sim
G_1 + G_2 +G_3 +(3-u)E.
\]
We calculate $S(E)$ and $S(W_{\bullet,\bullet}^{E},p)$.
Take $u \in \R_{\geq 0}$. 
Let $\Ti{P}(u)+\Ti{N}(u)$ be the Zariski decomposition of $\sigma^{\ast}(-K_{S})-uE$,
where $\Ti{P}(u)$ is the positive part and $\Ti{N}(u)$ is the negative part.
If $u\in [0,2]$, then we have
\begin{align*}
&\Ti{P}(u)=G_1 + G_2 +G_3 +(3-u)E,\\
&\Ti{N}(u)=0.
\end{align*}
If $u\in [2,3]$, then we have
\begin{align*}
&\Ti{P}(u)=(3-u)\left(G_1 + G_2 +G_3 +E\right),\\
&\Ti{N}(u)=(u-2)\left(G_1 + G_2 +G_3 \right).
\end{align*}
We note that $\sigma^{\ast}(-K_{S})-uE$ is not pseudo effective for $u>3$.
If $u\in [0,2]$, then we have
\begin{align*}
\Ti{P}(u)^2=(6-u^2),\quad \Ti{P}(u)E = u.
\end{align*}
If $u\in [2,3]$, then we have
\begin{align*}
\Ti{P}(u)^2=2(3-u)^2,\quad \Ti{P}(u)E = 2(3-u).
\end{align*}

Therefore, we get 
\begin{align*}
S(E)%=\frac{1}{5}\int_{0}^{2} 5 -u^2 du 
=\frac{5}{3}, \quad
S(W_{\bullet,\bullet}^{E},p)
%=\frac{2}{5}\int_{0}^{1}2u \Rm{ord}_{p}(N(u)|_{F}) du
%\frac{1}{5}\int_{0}^{2}u^2 du 
= \frac{2}{3}.
\end{align*}
Hence we have 
\begin{align*}
\frac{6}{5}  \geq \delta_{p}(S) \geq  \Rm{min}\left\{ \frac{2}{S(E)}, \frac{1}{S(W_{\bullet,\bullet}^{E},p)} \right\}
= \frac{6}{5} 
\end{align*}
from Corollary~\ref{AZ}.
 Thus, we have $\delta_{p}(S)= 6/5$ in this case.

\end{proof}

\section{The case of the anti-canonical degree $7$}
It is known that there exist $2$ types of weak del Pezzo surfaces of the anti-canonical degree $7$
in terms of the configuration of negative curves (\cite{CT},\cite[\S 8.4]{Dol}).

\begin{prop}
Let $S$ be the anti-canonical degree $7$ weak del Pezzo surface such
that the dual graph of negative curves is
\[
\xygraph{
   {\bullet}* +!D{ }*+!U{E_1}
                    -[r]
    {\bullet}*+!U{E_{2}}  - [r]
    {\circ}*+!U{F}
   %     - [r] \cdots
%       - [r] {\circ}* +!U{F_{2}}
     %  - [r] {\circ}* +!U{F_{3}}
       %  - [r] {\bullet}*+!U{E_{4}}
              },
        \]
where $E_i$ $(i=1,2)$ is a $(-1)$-curve and $F$ is a $(-2)$-curve.
Then, for a point $p \in S$, it holds that
\begin{equation}
\nonumber
  \delta_{p}(S)=
  \begin{cases}
    \frac{21}{25} & \text{if $p \in E_1 \setminus E_2$,} \\
   % \frac{21}{31}                 & \text{if $ p \in F_2 \setminus E $,} \\
     \frac{21}{31}                 & \text{if $p \in E_2 $,} \\
     \frac{7}{9} & \text{if $p \in  F \setminus E_2 $,} \\
     \frac{21}{23} & \text{if $p \in S \setminus \left( E_1 \cup E_2 \cup F \right)$.}
  \end{cases}
\end{equation}
\end{prop}

\begin{proof}
We denote $D= \sum_{i=1,2} a_i E_i + F \in \Rm{Div}(S)$ ($a_i, b \in \Z$)
by $D=(a_1, a_2, b)$.
The intersection matrix of $\{E_1, E_2, F  \}$ is 
\[A:=
\left(
\begin{array}{cc|c}
  -1 & 1 &  0        \\
  1 & -2 &  1        \\  \hline
  0 & 1 &  -2     \\ 
 % 1 & 0 &  0     \\ 
%  0 & 1 &  1  &  0       \\
\end{array}
\right).
\]
 We note that 
 $
 -K_{S} \sim 3E_1 + 4E_2 + 2F = (3,4,2).
$
%%%%%%%%%%%%%%%%%%%%% (1)

\noindent
(1)\;The case $p \in E_1 \setminus E_2$.

\noindent
We calculate $S(E_{1})$ and $S(W_{\bullet,\bullet}^{E_{1}},p)$.
Take $u \in \R_{\geq 0}$. 
Let $P(u)+N(u)$ be the Zariski decomposition of $-K_{S}-uE_{1}$.
%where $P(u)$ is the positive part and $N(u)$ is the negative part.
%
If $u\in [0,1]$, then we have
\begin{align*}
&P(u)=\left(3-u, 4, 2 \right),\\
%(2-u)E_1+ 3E_2 + \left(3-\frac{3}{4}u\right)F_1 + \left(4-\frac{u}{2}\right)F_2 + \left(2-\frac{u}{4}\right)F_3. \\
&N(u)=\left(0, 0 , 0 \right).
%\frac{3}{4}uF_1 + \frac{u}{2}F_2 +  \frac{u}{4}F_3.
\end{align*}
If $u\in [1,3]$, then we have
\begin{align*}
&P(u)=\left(3-u , 2(3-u) , 3-u \right),\\
%(2-u)E_1+ 3E_2 + \left(3-\frac{3}{4}u\right)F_1 + \left(4-\frac{u}{2}\right)F_2 + \left(2-\frac{u}{4}\right)F_3. \\
&N(u)=\left(0, 2(u-1) , u-1  \right).
%\frac{3}{4}uF_1 + \frac{u}{2}F_2 +  \frac{u}{4}F_3.
\end{align*}
We note that $-K_{S}-uE_{1}$ is not pseudo effective for $u>3$.
If $u\in [0,1]$, then we have 
\begin{align*}
P(u)^2=7-2u-u^2, \quad P(u)E_1 =1+u .
\end{align*}
If $u\in [1,3]$, then we have 
\begin{align*}
P(u)^2=(3-u)^2,\quad P(u)E_1 = 3-u.
\end{align*}

Therefore,  we get 
\begin{align*}
S(E_1)%=\frac{1}{4}\int_{0}^{2}(2-u)(2+u) du
%+ \frac{1}{5}\int_{1}^{2} 2(2-u)^2 du
=\frac{25}{21}, \quad
S(W_{\bullet,\bullet}^{E_1},p)
%=\frac{2}{5}\int_{0}^{1} \Rm{ord}_{p}(N(u)|_{F}) du
%&=%\frac{1}{5}\int_{0}^{1} (1+u)^2 du
%+ \frac{2}{5}\int_{0}^{2}\frac{5}{4}u \cdot \Rm{ord}_{p}(N(u)|_{E_{1}})  du
%+\frac{1}{5}\int_{0}^{2}\frac{25}{16}u^2 du\\
 &=\frac{15}{21}.
 %\begin{cases}
  %\frac{4}{3} & \text{if $p \in  F_{2}\cap F_{1} $},\\
    %                & \\
   %\frac{5}{4}   & \text{if $p \in F_{2} \cap ( E_1 \cup E_2)$}, \\
    %          & \\
  % \frac{1}{2}   & \text{if $p \in F_{2} \setminus  (F_{1}\cup E_1 \cup E_2)$}, \\
%  \end{cases}
\end{align*}
Hence we have 
\begin{align*}
\frac{21}{25}  \geq \delta_{p}(S) \geq  \Rm{min}\left\{ \frac{1}{S(E_1)}, \frac{1}{S(W_{\bullet,\bullet}^{E_1},p)} \right\}
= \frac{21}{25} 
\end{align*}
from Corollary~\ref{AZ}.
Thus, we have $\delta_{p}(S)= 21/25$ in this case.

%%%%%%%%%%%%%%%%%%%%% (2)

\noindent
(2)\;The case $p \in E_2 $.

\noindent
We calculate $S(E_{2})$ and $S(W_{\bullet,\bullet}^{E_{2}},p)$.
Take $u \in \R_{\geq 0}$. 
Let $P(u)+N(u)$ be the Zariski decomposition of $-K_{S}-uE_{2}$.
%where $P(u)$ is the positive part and $N(u)$ is the negative part.
%
If $u\in [0,1]$, then we have
\begin{align*}
&P(u)=\left(3,4-u, 2-\frac{u}{2} \right),\\
%(2-u)E_1+ 3E_2 + \left(3-\frac{3}{4}u\right)F_1 + \left(4-\frac{u}{2}\right)F_2 + \left(2-\frac{u}{4}\right)F_3. \\
&N(u)=\left(0, 0 , \frac{u}{2} \right).
%\frac{3}{4}uF_1 + \frac{u}{2}F_2 +  \frac{u}{4}F_3.
\end{align*}
If $u\in [1,4]$, then we have
\begin{align*}
&P(u)=\left(4-u, 4-u, 2-\frac{u}{2} \right),\\
%(2-u)E_1+ 3E_2 + \left(3-\frac{3}{4}u\right)F_1 + \left(4-\frac{u}{2}\right)F_2 + \left(2-\frac{u}{4}\right)F_3. \\
&N(u)=\left(u-1, 0, \frac{u}{2} \right).
%\frac{3}{4}uF_1 + \frac{u}{2}F_2 +  \frac{u}{4}F_3.
\end{align*}
We note that $-K_{S}-uE_{2}$ is not pseudo effective for $u>4$.
If $u\in [0,1]$, then we have 
\begin{align*}
P(u)^2=7-2u-\frac{u^2}{2},\quad P(u)E_2 = 1+\frac{u}{2}.
\end{align*}
If $u\in [1,4]$, then we have 
\begin{align*}
P(u)^2=2\left(2-\frac{u}{2}\right)^2,\quad P(u)E_2 = 2-\frac{u}{2}.
\end{align*}

Therefore,  we get 
\begin{align*}
S(E_2)%=\frac{1}{4}\int_{0}^{2}(2-u)(2+u) du
%+ \frac{1}{5}\int_{1}^{2} 2(2-u)^2 du
=\frac{31}{21}, \quad 
S(W_{\bullet,\bullet}^{E_{2}},p)
%=\frac{2}{5}\int_{0}^{1} \Rm{ord}_{p}(N(u)|_{F}) du
%&=%\frac{1}{5}\int_{0}^{1} (1+u)^2 du
%+ \frac{2}{5}\int_{0}^{2}\frac{5}{4}u \cdot \Rm{ord}_{p}(N(u)|_{E_{1}})  du
%+\frac{1}{5}\int_{0}^{2}\frac{25}{16}u^2 du\\
 %=\frac{8}{9}
= \begin{cases}
 \frac{25}{21} & \text{if $p \in  E_2 \cap E_1  $},\\
                    & \\
 \frac{9}{7} & \text{if $p \in E_2 \cap F$},\\
                   & \\
  \frac{23}{42}   & \text{if $p \in E_{2} \setminus (F \cup E_1)$}. \\
  \end{cases}
\end{align*}
%Hence we have 
%\begin{align*}
%\frac{3}{4}  \geq \delta_{p}(S) \geq  \Rm{min}\left\{ \frac{1}{S(F_1)}, \frac{1}{S(W_{\bullet,\bullet}^{F_1},p)} \right\}
%= \frac{3}{4} 
%\end{align*}
%from \Red{Theorem~\ref{AZ}}.
Thus, we have $\delta_{p}(S)= 21/31$ for $p \in E_{2}$.

%%%%%%%%%%%%%%%%%%%%% (3)

\noindent
(3)\;The case $p \in F \setminus E_2$ .

\noindent
We calculate $S(F)$ and $S(W_{\bullet,\bullet}^{F},p)$.
Take $u \in \R_{\geq 0}$. 
Let $P(u)+N(u)$ be the Zariski decomposition of $-K_{S}-uF$.
%where $P(u)$ is the positive part and $N(u)$ is the negative part.
%
If $u\in [0,1]$, then we have
\begin{align*}
&P(u)=\left(3,4, 2-u \right),\\
%(2-u)E_1+ 3E_2 + \left(3-\frac{3}{4}u\right)F_1 + \left(4-\frac{u}{2}\right)F_2 + \left(2-\frac{u}{4}\right)F_3. \\
&N(u)=\left(0, 0 , 0 \right).
%\frac{3}{4}uF_1 + \frac{u}{2}F_2 +  \frac{u}{4}F_3.
\end{align*}
If $u\in [1,2]$, then we have
\begin{align*}
&P(u)=\left(3, 5-u, 2-u \right),\\
%(2-u)E_1+ 3E_2 + \left(3-\frac{3}{4}u\right)F_1 + \left(4-\frac{u}{2}\right)F_2 + \left(2-\frac{u}{4}\right)F_3. \\
&N(u)=\left(0, u-1 , 0 \right).
%\frac{3}{4}uF_1 + \frac{u}{2}F_2 +  \frac{u}{4}F_3.
\end{align*}
We note that $-K_{S}-uF$ is not pseudo effective for $u>2$.
If $u\in [0,1]$, then we have 
\begin{align*}
P(u)^2=(7-2u^2),\quad P(u)F = 2u.
\end{align*}
If $u\in [1,2]$, then we have 
\begin{align*}
P(u)^2=8-2u-u^2,\quad P(u)F = 1+u.
\end{align*}
Therefore,  we get 
\begin{align*}
S(F)%=\frac{1}{4}\int_{0}^{2}(2-u)(2+u) du
%+ \frac{1}{5}\int_{1}^{2} 2(2-u)^2 du
=\frac{9}{7}, \quad
S(W_{\bullet,\bullet}^{F},p)
%=\frac{2}{5}\int_{0}^{1} \Rm{ord}_{p}(N(u)|_{F}) du
%&=%\frac{1}{5}\int_{0}^{1} (1+u)^2 du
%+ \frac{2}{5}\int_{0}^{2}\frac{5}{4}u \cdot \Rm{ord}_{p}(N(u)|_{E_{1}})  du
%+\frac{1}{5}\int_{0}^{2}\frac{25}{16}u^2 du\\
 &=\frac{23}{21}.
% \begin{cases}
 % \frac{17}{15} & \text{if $p \in E_{1}\cap F_{j} $ for $j=1,2$},\\
 %                   & \\
 %  \frac{7}{15}   & \text{if $p \in E_{1} \setminus \bigcup_{j=1}^{2} F_{j}$}, \\
%  \end{cases}
\end{align*}
Hence we have 
\begin{align*}
\frac{7}{9}  \geq \delta_{p}(S) \geq  \Rm{min}\left\{ \frac{1}{S(F)}, \frac{1}{S(W_{\bullet,\bullet}^{F},p)} \right\}
= \frac{7}{9} 
\end{align*}
from Corollary~\ref{AZ}.
Thus, we have $\delta_{p}(S)= 7/9$ in this case.

%%%%%%%%%%%%%%%%%%%%% (4)

\noindent
(4)\;The case $p \in S \setminus \left( E_1 \cup E_2 \cup F \right)$.

\noindent
Let $L \in |E_1+E_2|$ be a smooth irreducible curve.
We calculate $S(L)$ and $S(W_{\bullet,\bullet}^{L},p)$.
Take $u \in \R_{\geq 0}$. 
Let $P(u)+N(u)$ be the Zariski decomposition of $-K_{S}-uL$.
%where $P(u)$ is the positive part and $N(u)$ is the negative part.
%
If $u\in [0,2]$, then we have
\begin{align*}
&P(u)=\left(3-u , 4-u  , 2-\frac{u}{2}\right),\\
%(2-u)E_1+ 3E_2 + \left(3-\frac{3}{4}u\right)F_1 + \left(4-\frac{u}{2}\right)F_2 + \left(2-\frac{u}{4}\right)F_3. \\
&N(u)=\left(0, 0 , \frac{u}{2} \right).
%\frac{3}{4}uF_1 + \frac{u}{2}F_2 +  \frac{u}{4}F_3.
\end{align*}
If $u\in [2,3]$, then we have
\begin{align*}
&P(u)=\left(3-u, 2(3-u), 3-u\right),\\
%(2-u)E_1+ 3E_2 + \left(3-\frac{3}{4}u\right)F_1 + \left(4-\frac{u}{2}\right)F_2 + \left(2-\frac{u}{4}\right)F_3. \\
&N(u)=\left(0, u-2 , u-1\right).
%\frac{3}{4}uF_1 + \frac{u}{2}F_2 +  \frac{u}{4}F_3.
\end{align*}
%If $u\in [2,3]$, it holds
%\begin{align*}
%&P(u)=\left(3-u,3-u,3-u,0, 3-u\right)\\
%(2-u)E_1+ 3E_2 + \left(3-\frac{3}{4}u\right)F_1 + \left(4-\frac{u}{2}\right)F_2 + \left(2-\frac{u}{4}\right)F_3. \\
%&N(u)=\left(u-1, 0 , u-2 , 0, u-1 \right)
%\frac{3}{4}uF_1 + \frac{u}{2}F_2 +  \frac{u}{4}F_3.
%\end{align*}
We note that $-K_{S}-uL$ is not pseudo effective for $u>3$.
If $u\in [0,2]$, then we have 
\begin{align*}
P(u)^2=7-4u+\frac{u^2}{2},\quad P(u)L = 2-\frac{u}{2}.
\end{align*}
If $u\in [2,3]$, then we have 
\begin{align*}
P(u)^2=(3-u)^2,\quad P(u)L = 3-u.
\end{align*}
%f $u\in [2,3]$ we have 
%\begin{align*}
%P(u)^2=(3-u)^2,\quad P(u)F = (3-u)
%\end{align*}
Therefore,  we get 
\begin{align*}
S(L)%=\frac{1}{4}\int_{0}^{2}(2-u)(2+u) du
%+ \frac{1}{5}\int_{1}^{2} 2(2-u)^2 du
=\frac{23}{21}, \quad 
S(W_{\bullet,\bullet}^{L},p)
%=\frac{2}{5}\int_{0}^{1} \Rm{ord}_{p}(N(u)|_{F}) du
%&=%\frac{1}{5}\int_{0}^{1} (1+u)^2 du
%+ \frac{2}{5}\int_{0}^{2}\frac{5}{4}u \cdot \Rm{ord}_{p}(N(u)|_{E_{1}})  du
%+\frac{1}{5}\int_{0}^{2}\frac{25}{16}u^2 du\\
 &=\frac{15}{21}.
% \begin{cases}
 % \frac{17}{15} & \text{if $p \in E_{1}\cap F_{j} $ for $j=1,2$},\\
 %                   & \\
 %  \frac{7}{15}   & \text{if $p \in E_{1} \setminus \bigcup_{j=1}^{2} F_{j}$}, \\
%  \end{cases}
\end{align*}
Hence we have 
\begin{align*}
\frac{21}{23} \geq \delta_{p}(S) \geq  \Rm{min}\left\{ \frac{1}{S(L)}, \frac{1}{S(W_{\bullet,\bullet}^{L},p)} \right\}
= \frac{21}{23}
\end{align*}
from Corollary~\ref{AZ}.
Thus, we have $\delta_{p}(S)= 21/23$ in this case.
\end{proof}

\begin{prop}
Let $S$ be the anti-canonical degree $6$ weak del Pezzo surface such
that the dual graph of negative curves is
\[
\xygraph{
   {\bullet}* +!D{ }*+!U{E_1}
                    -[r]
    {\bullet}*+!U{E_{2}}  - [r]
    {\bullet}*+!U{E_3}
   %     - [r] \cdots
%       - [r] {\circ}* +!U{F_{2}}
     %  - [r] {\circ}* +!U{F_{3}}
       %  - [r] {\bullet}*+!U{E_{4}}
              },
        \]
where $E_i$ $(i=1,2,3)$ is a $(-1)$-curve.
% $F$ means a $(-2)$-curve.
Then, for a point $p \in S$, it holds that
\begin{equation}
\nonumber
  \delta_{p}(S)=
  \begin{cases}
    \frac{21}{23} & \text{if $p \in E_1 \setminus E_2$,} \\
   % \frac{21}{31}                 & \text{if $ p \in F_2 \setminus E $,} \\
     %\frac{21}{31}                 & \text{if $p \in E_2 $,} \\
     \frac{21}{25} & \text{if $p \in  E_2 $,} \\
     \frac{21}{22} & \text{if $p \in S \setminus \left( E_1 \cup E_2 \cup E_3 \right)$.}
  \end{cases}
\end{equation}
\end{prop}

\begin{proof}
We denote $D= \sum_{i=1,2,3} a_i E_i  \in \Rm{Div}(S)$ ($a_i \in \Z$)
by $D=(a_1, a_2, a_3)$.
The intersection matrix of $\{E_1, E_2, E_3 \}$ is 
\[A:=
\left(
\begin{array}{ccc}
  -1 & 1 &  0        \\
  1 & -1 &  1        \\  
  0 & 1 &  -1     \\ 
 % 1 & 0 &  0     \\ 
%  0 & 1 &  1  &  0       \\
\end{array}
\right).
\]
 We note that 
 $
 -K_{S} \sim 2E_1 + 3E_2 + 2E_3 = (2,3,2).
$
%%%%%%%%%%%%%%%%%%%%% (1)

\noindent
(1)\;The case $p \in E_1 \setminus E_2$.

\noindent
We calculate $S(E_{1})$ and $S(W_{\bullet,\bullet}^{E_{1}},p)$.
Take $u \in \R_{\geq 0}$. 
Let $P(u)+N(u)$ be the Zariski decomposition of $-K_{S}-uE_{1}$.
%where $P(u)$ is the positive part and $N(u)$ is the negative part.
%
If $u\in [0,1]$, then we have
\begin{align*}
&P(u)=\left(2-u, 3, 2 \right),\\
%(2-u)E_1+ 3E_2 + \left(3-\frac{3}{4}u\right)F_1 + \left(4-\frac{u}{2}\right)F_2 + \left(2-\frac{u}{4}\right)F_3. \\
&N(u)=\left(0, 0 , 0 \right).
%\frac{3}{4}uF_1 + \frac{u}{2}F_2 +  \frac{u}{4}F_3.
\end{align*}
If $u\in [1,2]$, then we have
\begin{align*}
&P(u)=\left(2-u , 4-u , 2 \right),\\
%(2-u)E_1+ 3E_2 + \left(3-\frac{3}{4}u\right)F_1 + \left(4-\frac{u}{2}\right)F_2 + \left(2-\frac{u}{4}\right)F_3. \\
&N(u)=\left(0, u-1 , 0  \right).
%\frac{3}{4}uF_1 + \frac{u}{2}F_2 +  \frac{u}{4}F_3.
\end{align*}
We note that $-K_{S}-uE_{1}$ is not pseudo effective for $u>2$.
If $u\in [0,1]$, then we have 
\begin{align*}
P(u)^2=7-2u-u^2, \quad P(u)E_1 =1+u.
\end{align*}
If $u\in [1,2]$, then we have 
\begin{align*}
P(u)^2=4(2-u),\quad P(u)E_1 = 2.
\end{align*}

Therefore,  we get 
\begin{align*}
S(E_1)%=\frac{1}{4}\int_{0}^{2}(2-u)(2+u) du
%+ \frac{1}{5}\int_{1}^{2} 2(2-u)^2 du
=\frac{23}{21}, \quad
S(W_{\bullet,\bullet}^{E_1},p)
%=\frac{2}{5}\int_{0}^{1} \Rm{ord}_{p}(N(u)|_{F}) du
%&=%\frac{1}{5}\int_{0}^{1} (1+u)^2 du
%+ \frac{2}{5}\int_{0}^{2}\frac{5}{4}u \cdot \Rm{ord}_{p}(N(u)|_{E_{1}})  du
%+\frac{1}{5}\int_{0}^{2}\frac{25}{16}u^2 du\\
 &=\frac{19}{21}.
 %\begin{cases}
  %\frac{4}{3} & \text{if $p \in  F_{2}\cap F_{1} $},\\
    %                & \\
   %\frac{5}{4}   & \text{if $p \in F_{2} \cap ( E_1 \cup E_2)$}, \\
    %          & \\
  % \frac{1}{2}   & \text{if $p \in F_{2} \setminus  (F_{1}\cup E_1 \cup E_2)$}, \\
%  \end{cases}
\end{align*}
Hence we have 
\begin{align*}
\frac{21}{23}  \geq \delta_{p}(S) \geq  \Rm{min}\left\{ \frac{1}{S(E_1)}, \frac{1}{S(W_{\bullet,\bullet}^{E_1},p)} \right\}
= \frac{21}{23} 
\end{align*}
from Corollary~\ref{AZ}.
Thus, we have $\delta_{p}(S)= 21/23$ in this case.

%%%%%%%%%%%%%%%%%%%%%%%%%%%%%%%(2)
\noindent
(2)\;The case $p \in E_2$.

\noindent
We calculate $S(E_{2})$ and $S(W_{\bullet,\bullet}^{E_{2}},p)$.
Take $u \in \R_{\geq 0}$. 
Let $P(u)+N(u)$ be the Zariski decomposition of $-K_{S}-uE_{2}$.
%where $P(u)$ is the positive part and $N(u)$ is the negative part.
%
If $u\in [0,1]$, then we have
\begin{align*}
&P(u)=\left(2, 3-u, 2 \right),\\
%(2-u)E_1+ 3E_2 + \left(3-\frac{3}{4}u\right)F_1 + \left(4-\frac{u}{2}\right)F_2 + \left(2-\frac{u}{4}\right)F_3. \\
&N(u)=\left(0, 0 , 0 \right).
%\frac{3}{4}uF_1 + \frac{u}{2}F_2 +  \frac{u}{4}F_3.
\end{align*}
If $u\in [1,3]$, then we have
\begin{align*}
&P(u)=\left(3-u , 3-u , 3-u \right),\\
%(2-u)E_1+ 3E_2 + \left(3-\frac{3}{4}u\right)F_1 + \left(4-\frac{u}{2}\right)F_2 + \left(2-\frac{u}{4}\right)F_3. \\
&N(u)=\left(u-1, 0 , u-1  \right).
%\frac{3}{4}uF_1 + \frac{u}{2}F_2 +  \frac{u}{4}F_3.
\end{align*}
We note that $-K_{S}-uE_{2}$ is not pseudo effective for $u>3$.
If $u\in [0,1]$, then we have 
\begin{align*}
P(u)^2=7-2u-u^2.% \quad P(u)E_1 =1+u .
\end{align*}
If $u\in [1,3]$, then we have 
\begin{align*}
P(u)^2=(3-u)^2.%\quad P(u)E_1 = 3-u.
\end{align*}

Therefore,  we get 
\begin{align*}
S(E_2)%=\frac{1}{4}\int_{0}^{2}(2-u)(2+u) du
%+ \frac{1}{5}\int_{1}^{2} 2(2-u)^2 du
=\frac{25}{21}
\end{align*}
by the definition of $S(E_2)$.
%\begin{align*}
%S(W_{\bullet,\bullet}^{E_2},p)
%=\frac{2}{5}\int_{0}^{1} \Rm{ord}_{p}(N(u)|_{F}) du
%&=%\frac{1}{5}\int_{0}^{1} (1+u)^2 du
%+ \frac{2}{5}\int_{0}^{2}\frac{5}{4}u \cdot \Rm{ord}_{p}(N(u)|_{E_{1}})  du
%+\frac{1}{5}\int_{0}^{2}\frac{25}{16}u^2 du\\
% &=\frac{19}{21}
 %\begin{cases}
  %\frac{4}{3} & \text{if $p \in  F_{2}\cap F_{1} $},\\
    %                & \\
   %\frac{5}{4}   & \text{if $p \in F_{2} \cap ( E_1 \cup E_2)$}, \\
    %          & \\
  % \frac{1}{2}   & \text{if $p \in F_{2} \setminus  (F_{1}\cup E_1 \cup E_2)$}, \\
%  \end{cases}
%\end{align*}
%by \Red{Theorem~$S(W_{\bullet,\bullet}^{F},p)$の明示式}.
By \cite[\S 2]{FAND}, the (global) delta invariant of del Pezzo surfaces with the anti-canonical degree $7$ is $21/25$.
Hence we have $\delta_{p}(S)= 21/25$ for $p \in E_2$.

%%%%%%%%%%%%%%%%%%%%% (3)

\noindent
(3)\;The case $p \in S \setminus \left( E_1 \cup E_2 \cup E_3 \right)$.

\noindent
Let $L \in |E_2+E_3|$ be a smooth irreducible curve.
We calculate $S(L)$ and $S(W_{\bullet,\bullet}^{L},p)$.
Take $u \in \R_{\geq 0}$. 
Let $P(u)+N(u)$ be the Zariski decomposition of $-K_{S}-uL$.
%where $P(u)$ is the positive part and $N(u)$ is the negative part.
%
If $u\in [0,1]$, then we have
\begin{align*}
&P(u)=\left(2 , 3-u  , 2-u\right),\\
%(2-u)E_1+ 3E_2 + \left(3-\frac{3}{4}u\right)F_1 + \left(4-\frac{u}{2}\right)F_2 + \left(2-\frac{u}{4}\right)F_3. \\
&N(u)=\left(0, 0 , 0 \right).
%\frac{3}{4}uF_1 + \frac{u}{2}F_2 +  \frac{u}{4}F_3.
\end{align*}
If $u\in [1,2]$, then we have
\begin{align*}
&P(u)=\left(3-u, 3-u, 2-u\right),\\
%(2-u)E_1+ 3E_2 + \left(3-\frac{3}{4}u\right)F_1 + \left(4-\frac{u}{2}\right)F_2 + \left(2-\frac{u}{4}\right)F_3. \\
&N(u)=\left(u-1, 0 , 0\right).
%\frac{3}{4}uF_1 + \frac{u}{2}F_2 +  \frac{u}{4}F_3.
\end{align*}
%If $u\in [2,3]$, it holds
%\begin{align*}
%&P(u)=\left(3-u,3-u,3-u,0, 3-u\right)\\
%(2-u)E_1+ 3E_2 + \left(3-\frac{3}{4}u\right)F_1 + \left(4-\frac{u}{2}\right)F_2 + \left(2-\frac{u}{4}\right)F_3. \\
%&N(u)=\left(u-1, 0 , u-2 , 0, u-1 \right)
%\frac{3}{4}uF_1 + \frac{u}{2}F_2 +  \frac{u}{4}F_3.
%\end{align*}
We note that $-K_{S}-uL$ is not pseudo effective for $u>2$.
If $u\in [0,1]$, then we have 
\begin{align*}
P(u)^2=7-2u,\quad P(u)L = 2.
\end{align*}
If $u\in [1,2]$, then we have 
\begin{align*}
P(u)^2=(2-u)(4-u),\quad P(u)L = 3-u.
\end{align*}
%f $u\in [2,3]$ we have 
%\begin{align*}
%P(u)^2=(3-u)^2,\quad P(u)F = (3-u)
%\end{align*}
Therefore,  we get 
\begin{align*}
S(L)%=\frac{1}{4}\int_{0}^{2}(2-u)(2+u) du
%+ \frac{1}{5}\int_{1}^{2} 2(2-u)^2 du
=\frac{22}{21}, \quad 
S(W_{\bullet,\bullet}^{L},p)
%=\frac{2}{5}\int_{0}^{1} \Rm{ord}_{p}(N(u)|_{F}) du
%&=%\frac{1}{5}\int_{0}^{1} (1+u)^2 du
%+ \frac{2}{5}\int_{0}^{2}\frac{5}{4}u \cdot \Rm{ord}_{p}(N(u)|_{E_{1}})  du
%+\frac{1}{5}\int_{0}^{2}\frac{25}{16}u^2 du\\
 &=\frac{19}{21}.
% \begin{cases}
 % \frac{17}{15} & \text{if $p \in E_{1}\cap F_{j} $ for $j=1,2$},\\
 %                   & \\
 %  \frac{7}{15}   & \text{if $p \in E_{1} \setminus \bigcup_{j=1}^{2} F_{j}$}, \\
%  \end{cases}
\end{align*}
Hence we have 
\begin{align*}
\frac{21}{22} \geq \delta_{p}(S) \geq  \Rm{min}\left\{ \frac{1}{S(L)}, \frac{1}{S(W_{\bullet,\bullet}^{L},p)} \right\}
= \frac{21}{22}
\end{align*}
from Corollary~\ref{AZ}.
Thus, we have $\delta_{p}(S)= 21/22$ in this case.

\end{proof}

\section{The case of the anti-canonical degree $8$}
We denote by $\pi:\Sigma_n \to \Bb{P}^1$  the $n$-th Hirzebruch surface.
Let $C_0$ be the section of $\pi$ with $C_0^2=-n$ and $\Gamma$ the fiber of $\pi$.
It is known that a weak del Pezzo surface of the anti-canonical degree $8$
is either $\Sigma_0$, $\Sigma_1$ or $\Sigma_2$ (\cite{CT},\cite[\S 8.4]{Dol}).

\begin{prop}
If $S = \Sigma_2$, then for any point $p \in S$, it holds that
\begin{equation}
\nonumber
  \delta_{p}(S)=\frac{3}{4}.
\end{equation}
\end{prop}

\begin{proof}

We denote $D= aC_0 + b\Gamma \in \Rm{Div}(S)$ ($a, b \in \Z$)
by $D=(a, b)$.
The intersection matrix of $\{C_0, \Gamma  \}$ is 
\[A:=
\left(
\begin{array}{cc}
  -2 & 1           \\ 
  1 & 0           \\  
 % 0 & 1      \\ 
 % 1 & 0 &  0     \\ 
%  0 & 1 &  1  &  0       \\
\end{array}
\right).
\]
 We note that 
 $
 -K_{S} \sim 2C_0 + 4\Gamma = (2,4).
$
%%%%%%%%%%%%%%%%%%%%% (1)

\noindent
(1)\;The case $p \in C_0$.

\noindent
We calculate $S(C_0)$ and $S(W_{\bullet,\bullet}^{C_0},p)$.
Take $u \in \R_{\geq 0}$. 
Let $P(u)+N(u)$ be the Zariski decomposition of $-K_{S}-uC_0$.
%where $P(u)$ is the positive part and $N(u)$ is the negative part.
%
If $u\in [0,2]$, then we have
\begin{align*}
&P(u)=\left(2-u, 4\right),\\
%(2-u)E_1+ 3E_2 + \left(3-\frac{3}{4}u\right)F_1 + \left(4-\frac{u}{2}\right)F_2 + \left(2-\frac{u}{4}\right)F_3. \\
&N(u)=\left(0, 0  \right).
%\frac{3}{4}uF_1 + \frac{u}{2}F_2 +  \frac{u}{4}F_3.
\end{align*}
%If $u\in [1,3]$, it holds
%\begin{align*}
%&P(u)=\left(3-u , 2(3-u) , 3-u \right)\\
%(2-u)E_1+ 3E_2 + \left(3-\frac{3}{4}u\right)F_1 + \left(4-\frac{u}{2}\right)F_2 + \left(2-\frac{u}{4}\right)F_3. \\
%&N(u)=\left(0, 2(u-1) , u-1  \right)
%\frac{3}{4}uF_1 + \frac{u}{2}F_2 +  \frac{u}{4}F_3.
%\end{align*}
We note that $-K_{S}-uC_0$ is not pseudo effective for $u>2$.
If $u\in [0,2]$, then we have
\begin{align*}
P(u)^2=2(4-u^2), \quad P(u)C_0 =2u .
\end{align*}
%If $u\in [1,3]$ we have 
%\begin{align*}
%P(u)^2=(3-u)^2,\quad P(u)C_0 = 3-u.
%\end{align*}

Therefore,  we get 
\begin{align*}
S(C_0)%=\frac{1}{4}\int_{0}^{2}(2-u)(2+u) du
%+ \frac{1}{5}\int_{1}^{2} 2(2-u)^2 du
=\frac{4}{3}, \quad
S(W_{\bullet,\bullet}^{C_0},p)
%=\frac{2}{5}\int_{0}^{1} \Rm{ord}_{p}(N(u)|_{F}) du
%&=%\frac{1}{5}\int_{0}^{1} (1+u)^2 du
%+ \frac{2}{5}\int_{0}^{2}\frac{5}{4}u \cdot \Rm{ord}_{p}(N(u)|_{E_{1}})  du
%+\frac{1}{5}\int_{0}^{2}\frac{25}{16}u^2 du\\
 &=\frac{4}{3}.
 %\begin{cases}
  %\frac{4}{3} & \text{if $p \in  F_{2}\cap F_{1} $},\\
    %                & \\
   %\frac{5}{4}   & \text{if $p \in F_{2} \cap ( E_1 \cup E_2)$}, \\
    %          & \\
  % \frac{1}{2}   & \text{if $p \in F_{2} \setminus  (F_{1}\cup E_1 \cup E_2)$}, \\
%  \end{cases}
\end{align*}
Hence we have 
\begin{align*}
\frac{3}{4}  \geq \delta_{p}(S) \geq  \Rm{min}\left\{ \frac{1}{S(C_0)}, \frac{1}{S(W_{\bullet,\bullet}^{C_0},p)} \right\}
= \frac{3}{4} 
\end{align*}
from Corollary~\ref{AZ}.
Thus, we have $\delta_{p}(S)= 3/4$ in this case.

%%%%%%%%%%%%%%%%%%%%% (2)

\noindent
(2)\;The case $p \in S\setminus C_0$.

\noindent
Let $\Gamma$ be the fiber of $\pi$ passing through $p$.
We calculate $S(\Gamma)$ and $S(W_{\bullet,\bullet}^{\Gamma},p)$.
Take $u \in \R_{\geq 0}$. 
Let $P(u)+N(u)$ be the Zariski decomposition of $-K_{S}-u\Gamma$.
%where $P(u)$ is the positive part and $N(u)$ is the negative part.
%
If $u\in [0,4]$, then we have
\begin{align*}
&P(u)=\left(2-\frac{u}{2}, 4-u\right),\\
%(2-u)E_1+ 3E_2 + \left(3-\frac{3}{4}u\right)F_1 + \left(4-\frac{u}{2}\right)F_2 + \left(2-\frac{u}{4}\right)F_3. \\
&N(u)=\left(\frac{u}{2}, 0  \right).
%\frac{3}{4}uF_1 + \frac{u}{2}F_2 +  \frac{u}{4}F_3.
\end{align*}
%If $u\in [1,3]$, it holds
%\begin{align*}
%&P(u)=\left(3-u , 2(3-u) , 3-u \right)\\
%(2-u)E_1+ 3E_2 + \left(3-\frac{3}{4}u\right)F_1 + \left(4-\frac{u}{2}\right)F_2 + \left(2-\frac{u}{4}\right)F_3. \\
%&N(u)=\left(0, 2(u-1) , u-1  \right)
%\frac{3}{4}uF_1 + \frac{u}{2}F_2 +  \frac{u}{4}F_3.
%\end{align*}
We note that $-K_{S}-u\Gamma$ is not pseudo effective for $u>4$.
If $u\in [0,4]$, then we have
\begin{align*}
P(u)^2=2(2-\frac{u}{2})^2, \quad P(u)\Gamma =2-\frac{u}{2} .
\end{align*}
%If $u\in [1,3]$ we have 
%\begin{align*}
%P(u)^2=(3-u)^2,\quad P(u)C_0 = 3-u.
%\end{align*}

Therefore,  we get 
\begin{align*}
S(\Gamma)%=\frac{1}{4}\int_{0}^{2}(2-u)(2+u) du
%+ \frac{1}{5}\int_{1}^{2} 2(2-u)^2 du
=\frac{4}{3}, \quad
S(W_{\bullet,\bullet}^{\Gamma},p)
%=\frac{2}{5}\int_{0}^{1} \Rm{ord}_{p}(N(u)|_{F}) du
%&=%\frac{1}{5}\int_{0}^{1} (1+u)^2 du
%+ \frac{2}{5}\int_{0}^{2}\frac{5}{4}u \cdot \Rm{ord}_{p}(N(u)|_{E_{1}})  du
%+\frac{1}{5}\int_{0}^{2}\frac{25}{16}u^2 du\\
 &=\frac{2}{3}.
 %\begin{cases}
  %\frac{4}{3} & \text{if $p \in  F_{2}\cap F_{1} $},\\
    %                & \\
   %\frac{5}{4}   & \text{if $p \in F_{2} \cap ( E_1 \cup E_2)$}, \\
    %          & \\
  % \frac{1}{2}   & \text{if $p \in F_{2} \setminus  (F_{1}\cup E_1 \cup E_2)$}, \\
%  \end{cases}
\end{align*}
%by Theorem~\ref{S-inv}.
Hence we have 
\begin{align*}
\frac{3}{4}  \geq \delta_{p}(S) \geq  \Rm{min}\left\{ \frac{1}{S(\Gamma)}, \frac{1}{S(W_{\bullet,\bullet}^{\Gamma},p)} \right\}
= \frac{3}{4} 
\end{align*}
from Corollary~\ref{AZ}.
Thus, we have $\delta_{p}(S)= 3/4$ in this case.
\end{proof}

\begin{prop}
If $S = \Sigma_1$, then for any point $p \in S$, it holds that
\begin{equation}
\nonumber
  \delta_{p}(S)=
  \begin{cases}
    \frac{6}{7} & \text{if $p \in C_0$,} \\
   % \frac{21}{31}                 & \text{if $ p \in F_2 \setminus E $,} \\
     %\frac{21}{31}                 & \text{if $p \in E_2 $,} \\
    % \frac{21}{25} & \text{if $p \in  E_2 $} \\
     \frac{12}{13} & \text{if $p \in S \setminus C_0 $.}
  \end{cases}
\end{equation}
\end{prop}

\begin{proof}

We denote $D= aC_0 + b\Gamma \in \Rm{Div}(S)$ ($a, b \in \Z$)
by $D=(a, b)$.
The intersection matrix of $\{C_0, \Gamma  \}$ is 
\[A:=
\left(
\begin{array}{cc}
  -1 & 1           \\ 
  1 & 0           \\  
 % 0 & 1      \\ 
 % 1 & 0 &  0     \\ 
%  0 & 1 &  1  &  0       \\
\end{array}
\right).
\]
 We note that 
 $
 -K_{S} \sim 2C_0 + 3\Gamma = (2,3).
$
%%%%%%%%%%%%%%%%%%%%% (1)

\noindent
(1)\;The case $p \in C_0$.

\noindent
We calculate $S(C_0)$ and $S(W_{\bullet,\bullet}^{C_0},p)$.
Take $u \in \R_{\geq 0}$. 
Let $P(u)+N(u)$ be the Zariski decomposition of $-K_{S}-uC_0$.
%where $P(u)$ is the positive part and $N(u)$ is the negative part.
%
If $u\in [0,2]$, then we have
\begin{align*}
&P(u)=\left(2-u, 3\right),\\
%(2-u)E_1+ 3E_2 + \left(3-\frac{3}{4}u\right)F_1 + \left(4-\frac{u}{2}\right)F_2 + \left(2-\frac{u}{4}\right)F_3. \\
&N(u)=\left(0, 0  \right).
%\frac{3}{4}uF_1 + \frac{u}{2}F_2 +  \frac{u}{4}F_3.
\end{align*}
%If $u\in [1,3]$, it holds
%\begin{align*}
%&P(u)=\left(3-u , 2(3-u) , 3-u \right)\\
%(2-u)E_1+ 3E_2 + \left(3-\frac{3}{4}u\right)F_1 + \left(4-\frac{u}{2}\right)F_2 + \left(2-\frac{u}{4}\right)F_3. \\
%&N(u)=\left(0, 2(u-1) , u-1  \right)
%\frac{3}{4}uF_1 + \frac{u}{2}F_2 +  \frac{u}{4}F_3.
%\end{align*}
We note that $-K_{S}-uC_0$ is not pseudo effective for $u>2$.
If $u\in [0,2]$, then we have
\begin{align*}
P(u)^2=8-2u-u^2, \quad P(u)C_0 =u+1 .
\end{align*}
%If $u\in [1,3]$ we have 
%\begin{align*}
%P(u)^2=(3-u)^2,\quad P(u)C_0 = 3-u.
%\end{align*}

Therefore,  we get 
\begin{align*}
S(C_0)%=\frac{1}{4}\int_{0}^{2}(2-u)(2+u) du
%+ \frac{1}{5}\int_{1}^{2} 2(2-u)^2 du
=\frac{7}{6}, \quad
S(W_{\bullet,\bullet}^{C_0},p)
%=\frac{2}{5}\int_{0}^{1} \Rm{ord}_{p}(N(u)|_{F}) du
%&=%\frac{1}{5}\int_{0}^{1} (1+u)^2 du
%+ \frac{2}{5}\int_{0}^{2}\frac{5}{4}u \cdot \Rm{ord}_{p}(N(u)|_{E_{1}})  du
%+\frac{1}{5}\int_{0}^{2}\frac{25}{16}u^2 du\\
 &=\frac{13}{12}.
 %\begin{cases}
  %\frac{4}{3} & \text{if $p \in  F_{2}\cap F_{1} $},\\
    %                & \\
   %\frac{5}{4}   & \text{if $p \in F_{2} \cap ( E_1 \cup E_2)$}, \\
    %          & \\
  % \frac{1}{2}   & \text{if $p \in F_{2} \setminus  (F_{1}\cup E_1 \cup E_2)$}, \\
%  \end{cases}
\end{align*}
Hence we have 
\begin{align*}
\frac{6}{7}  \geq \delta_{p}(S) \geq  \Rm{min}\left\{ \frac{1}{S(C_0)}, \frac{1}{S(W_{\bullet,\bullet}^{C_0},p)} \right\}
= \frac{6}{7} 
\end{align*}
from Corollary~\ref{AZ}.
Thus, we have $\delta_{p}(S)= 6/7$ in this case.

%%%%%%%%%%%%%%%%%%%%% (2)

\noindent
(2)\;The case $p \in S\setminus C_0$.

\noindent
Let $\Gamma$ be the fiber of $\pi$ passing through $p$.
We calculate $S(\Gamma)$ and $S(W_{\bullet,\bullet}^{\Gamma},p)$.
Take $u \in \R_{\geq 0}$. 
Let $P(u)+N(u)$ be the Zariski decomposition of $-K_{S}-u\Gamma$.
%where $P(u)$ is the positive part and $N(u)$ is the negative part.
%
If $u\in [0,1]$, then we have
\begin{align*}
&P(u)=\left(2, 3-u\right),\\
%(2-u)E_1+ 3E_2 + \left(3-\frac{3}{4}u\right)F_1 + \left(4-\frac{u}{2}\right)F_2 + \left(2-\frac{u}{4}\right)F_3. \\
&N(u)=\left(0 , 0  \right).
%\frac{3}{4}uF_1 + \frac{u}{2}F_2 +  \frac{u}{4}F_3.
\end{align*}
If $u\in [1,3]$, then we have
\begin{align*}
&P(u)=\left(3-u , 3-u \right),\\
%(2-u)E_1+ 3E_2 + \left(3-\frac{3}{4}u\right)F_1 + \left(4-\frac{u}{2}\right)F_2 + \left(2-\frac{u}{4}\right)F_3. \\
&N(u)=\left( u-1 , 0  \right).
%\frac{3}{4}uF_1 + \frac{u}{2}F_2 +  \frac{u}{4}F_3.
\end{align*}
We note that $-K_{S}-u\Gamma$ is not pseudo effective for $u>3$.
If $u\in [0,1]$, then we have
\begin{align*}
P(u)^2=8-4u, \quad P(u)\Gamma =2 .
\end{align*}
If $u\in [1,3]$, then we have
\begin{align*}
P(u)^2=(3-u)^2,\quad P(u)\Gamma = 3-u.
\end{align*}

Therefore,  we get 
\begin{align*}
S(\Gamma)%=\frac{1}{4}\int_{0}^{2}(2-u)(2+u) du
%+ \frac{1}{5}\int_{1}^{2} 2(2-u)^2 du
=\frac{13}{12}, \quad
S(W_{\bullet,\bullet}^{\Gamma},p)
%=\frac{2}{5}\int_{0}^{1} \Rm{ord}_{p}(N(u)|_{F}) du
%&=%\frac{1}{5}\int_{0}^{1} (1+u)^2 du
%+ \frac{2}{5}\int_{0}^{2}\frac{5}{4}u \cdot \Rm{ord}_{p}(N(u)|_{E_{1}})  du
%+\frac{1}{5}\int_{0}^{2}\frac{25}{16}u^2 du\\
 &=\frac{5}{6}.
 %\begin{cases}
  %\frac{4}{3} & \text{if $p \in  F_{2}\cap F_{1} $},\\
    %                & \\
   %\frac{5}{4}   & \text{if $p \in F_{2} \cap ( E_1 \cup E_2)$}, \\
    %          & \\
  % \frac{1}{2}   & \text{if $p \in F_{2} \setminus  (F_{1}\cup E_1 \cup E_2)$}, \\
%  \end{cases}
\end{align*}
Hence we have 
\begin{align*}
\frac{12}{13}  \geq \delta_{p}(S) \geq  \Rm{min}\left\{ \frac{1}{S(\Gamma)}, \frac{1}{S(W_{\bullet,\bullet}^{\Gamma},p)} \right\}
= \frac{12}{13} 
\end{align*}
from Corollary~\ref{AZ}.
Thus, we have $\delta_{p}(S)= 12/13$ in this case.
\end{proof}

\begin{prop}
If $S = \Sigma_0=\Bb{P}^1 \times \Bb{P}^1$, then for any point $p \in S$, it holds that
\begin{equation}
\nonumber
  \delta_{p}(S)=1.
\end{equation}
\end{prop}

\begin{proof}

We denote $D= aC_0 + b\Gamma \in \Rm{Div}(S)$ ($a, b \in \Z$)
by $D=(a, b)$.
The intersection matrix of $\{C_0, \Gamma  \}$ is 
\[A:=
\left(
\begin{array}{cc}
  0 & 1           \\ 
  1 & 0           \\  
 % 0 & 1      \\ 
 % 1 & 0 &  0     \\ 
%  0 & 1 &  1  &  0       \\
\end{array}
\right).
\]
 We note that 
 $
 -K_{S} \sim 2C_0 + 2\Gamma = (2,2).
$
%%%%%%%%%%%%%%%%%%%%% (1)

For any point $p \in S$, let $C_0$ be the section that passes through $p$.
We calculate $S(C_0)$ and $S(W_{\bullet,\bullet}^{C_0},p)$.
Take $u \in \R_{\geq 0}$. 
Let $P(u)+N(u)$ be the Zariski decomposition of $-K_{S}-uC_0$.
%where $P(u)$ is the positive part and $N(u)$ is the negative part.
%
If $u\in [0,2]$, then we have
\begin{align*}
&P(u)=\left(2-u, 2\right),\\
%(2-u)E_1+ 3E_2 + \left(3-\frac{3}{4}u\right)F_1 + \left(4-\frac{u}{2}\right)F_2 + \left(2-\frac{u}{4}\right)F_3. \\
&N(u)=\left(0, 0  \right).
%\frac{3}{4}uF_1 + \frac{u}{2}F_2 +  \frac{u}{4}F_3.
\end{align*}
%If $u\in [1,3]$, it holds
%\begin{align*}
%&P(u)=\left(3-u , 2(3-u) , 3-u \right)\\
%(2-u)E_1+ 3E_2 + \left(3-\frac{3}{4}u\right)F_1 + \left(4-\frac{u}{2}\right)F_2 + \left(2-\frac{u}{4}\right)F_3. \\
%&N(u)=\left(0, 2(u-1) , u-1  \right)
%\frac{3}{4}uF_1 + \frac{u}{2}F_2 +  \frac{u}{4}F_3.
%\end{align*}
We note that $-K_{S}-uC_0$ is not pseudo effective for $u>2$.
If $u\in [0,2]$, then we have 
\begin{align*}
P(u)^2=4(2-u), \quad P(u)C_0 =2 .
\end{align*}
%If $u\in [1,3]$ we have 
%\begin{align*}
%P(u)^2=(3-u)^2,\quad P(u)C_0 = 3-u.
%\end{align*}

Therefore,  we get 
\begin{align*}
S(C_0)%=\frac{1}{4}\int_{0}^{2}(2-u)(2+u) du
%+ \frac{1}{5}\int_{1}^{2} 2(2-u)^2 du
=1, \quad
S(W_{\bullet,\bullet}^{C_0},p)
%=\frac{2}{5}\int_{0}^{1} \Rm{ord}_{p}(N(u)|_{F}) du
%&=%\frac{1}{5}\int_{0}^{1} (1+u)^2 du
%+ \frac{2}{5}\int_{0}^{2}\frac{5}{4}u \cdot \Rm{ord}_{p}(N(u)|_{E_{1}})  du
%+\frac{1}{5}\int_{0}^{2}\frac{25}{16}u^2 du\\
 &=1.
 %\begin{cases}
  %\frac{4}{3} & \text{if $p \in  F_{2}\cap F_{1} $},\\
    %                & \\
   %\frac{5}{4}   & \text{if $p \in F_{2} \cap ( E_1 \cup E_2)$}, \\
    %          & \\
  % \frac{1}{2}   & \text{if $p \in F_{2} \setminus  (F_{1}\cup E_1 \cup E_2)$}, \\
%  \end{cases}
\end{align*}
Hence we have 
\begin{align*}
1  \geq \delta_{p}(S) \geq  \Rm{min}\left\{ \frac{1}{S(C_0)}, \frac{1}{S(W_{\bullet,\bullet}^{C_0},p)} \right\}
= 1
\end{align*}
from Corollary~\ref{AZ}.
Thus, we have $\delta_{p}(S)= 1$ in this case.
\end{proof}

\vspace{2\baselineskip}
{}
\bigskip
\bigskip

Hiroto Akaike

Department of Mathematics, 
Kyoto University,

Yoshida-honmachi, Sakyo-ku, Kyoto,
606-8501, Japan

e-mail: akaike.hiroto.6f@kyoto-u.ac.jp

\begin{thebibliography}{9}
\bibitem{AZ} H.\ Abban, Z.\ Zhuang, K-stability of Fano varieties via admissible flags, Forum.\ of.\ Mathematics.\ Pi.\ \textbf{10} (2022), 1--43.
\bibitem{FAND} C.\ Araujo, A.-M.\ Castravet, I.\ Cheltsov, K.\ Fujita, A.-S.\ Kaloghiros, J.\ Martinez-Garcia, C.\ Shramov, H.\ Suss, N.\ Viswanathan, The Calabi problem for Fano threefolds, https://www.maths.ed.ac.uk/cheltsov/pdf/Fanos.pdf
%\bibitem{Ber16} R.\ J.\ Berman, K-polystability of $\Bb{Q}$-Fano varieties admitting K\"ahler-Einstein metrics, Invent.\ Math.,
\textbf{203} (2016), 973--1025.
\bibitem{BJ} H.\ Blum, M.\ Jonsson, Thresholds, valuations, and K-stability, Adv.\  Math.\ \textbf{365} (2020), 107062.
\bibitem{BX}H.\ Blum, C.\ Xu, Uniqueness of K-polystable degenerations of Fano varieties, Ann.\ Math.\ \textbf{190}
(2019), 609--656.
\bibitem{CFKO}  I.\ Cheltsov, K.\ Fujita, T.\ Kishimoto, T.\ Okada, K-stable divisors in $\Bb{P}^1 \times \Bb{P}^1 \times \Bb{P}^1$ of degree $(1,1,2)$, arXiv:2206.08539v2 [math.AG]
%\bibitem{CDS} X.\ -X.\ Chen, S.\ Donaldson, S.\ Sun, K\"ahler-Einstein metrics on Fano manifolds. I, II, III, J.\ Am.\
Math.\ Soc.\ \textbf{28} (2015), no. 1, 183--197, 199--234, 235--278.
\bibitem{CT} D.\ Coray, M.\ Tsfasman, Arithmetic on singular Del Pezzo surfaces, Proc.\ LMS \textbf{57} (1988), 25--87.
\bibitem{Dol}I.\ Dolgachev, Classical algebraic geometry.\ A modern view, Cambridge University Press, (2012).
\bibitem{FO} K.\ Fujita, Y.\ Odaka, On the K-stability of Fano varieties and anticanonical divisors, Tohoku Math.\ J. \textbf{70} (2018), 511--521.
\bibitem{Fjt19}K.\ Fujita, A valuative criterion for uniform K-stability of Q-Fano varieties, J.\ Reine Angew.\ Math.\ \textbf{751} (2019), 309--338.
\bibitem{Li17}C.\ Li, K-semistability is equivariant volume minimization, Duke Math.\ J.\ \textbf{166} (2017), 3147--3218.
\bibitem{LXZ}Y.\ Liu, C.\ Xu, Z.\ Zhuang, Finite generation for valuations computing stability thresholds and applications to K-stability, Ann.\ of Math.\ \textbf{196} (2022) 507--566.
%\bibitem{Ti15}G.\ Tian, K-stability and K\"ahler-Einstein metrics, Comm.\ Pure Appl.\ Math.\ \textbf{68} (2015), no. 7, 1085--1156.
%\bibitem{Aka} H.\ Akaike, Bounds for the order of automorphism groups of cyclic covering fibrations of a ruled surface, arxiv:2011.03000v3[math.AG].
%\bibitem{Ara} T.\ Arakawa, Bounds for the order of automorphism groups of hyperelliptic fibrations, T\^ohoku.\ Math.\ J. \textbf{50} (1998), 317--323.
%\bibitem{Bar} M.\ A.\ Barja, On the slope of bielliptic fibrations, Proc.\ Amer.\ Math.\ Soc.\ \textbf{129} (2001), 1899--1906.  
% \bibitem{Bar1} M.\ A.\ Barja, V.\ G.\ Alonso, and J.\ C.\ Naranjo, Xiao's Conjecture for general fibred surfaces, 
% J.\ reine angew Math.\ \textbf{739} (2018), 297--308.
% \bibitem{Chen} Z.\ Chen, Best bounds of automorphism groups of hyperelliptic fibrations, T\^ohoku.\ Math.\ J. \textbf{50} (1998), 469--489.
 %\bibitem{Barth1}W.Barth, K.Hulek, C.Weters and A.Van de Ven, Compact Complex 
 %Surfaces, 2nd Edition, Springer-Verlag, Berlin (2004).
% \bibitem{Eno} M.\ Enokizono, Slopes of fibered surfaces with a finite cyclic automorphism, 
% Michigan Math.\ J.\ \textbf{66} (2017), 125--154.
% \bibitem{Eno2} M.\ Enokizono, Upper bounds on the slope of certain fibered surfaces, arXiv:1604.06919v1[math.AG]
%\bibitem{Kod} K.\ Kodaira, On compact analytic surfaces I, Ann.\ of Math. \textbf{71} (1960),
%111--152; II, ibid.\ \textbf{77} (1963), 563--626;
%III, ibid.\ \textbf{78} (1963), 1--40.
%\bibitem{Kon} K.\ Konno, A lower bound of the slope of trigonal fibrations, Intern.\ J.\ Math. \textbf{7} (1996), 19--27.
%\bibitem{Mir} R.\ Miranda and U.\ Persson, Torsion groups of elliptic surfaces, Compositio Math \textbf{72} (1989), 249--267.
\end{thebibliography}
\end{document}